\documentclass[a4paper,11pt,oneside]{article}

\makeindex
\usepackage{hyperref}
\usepackage{graphicx}
\usepackage{enumerate}
\hypersetup{colorlinks, linkcolor=blue, urlcolor=red, filecolor=green, citecolor=blue}
\usepackage[toc,page]{appendix}
\usepackage{subcaption,siunitx} 
\usepackage[shortlabels]{enumitem}

\usepackage{amsmath}
\usepackage{amsthm}
\usepackage{amssymb}
\usepackage{authblk}
\usepackage{enumerate}
\usepackage{fullpage}
\usepackage{esint}
\usepackage{dsfont} 
\usepackage{colonequals} 

\usepackage{pgfplots}
\pgfplotsset{compat=1.15}
\usepackage{mathrsfs}
\usetikzlibrary{arrows}

\theoremstyle{plain}
\newtheorem{theorem}{Theorem}[section]
\newtheorem{corollary}[theorem]{Corollary}
\newtheorem{lemma}[theorem]{Lemma}
\newtheorem{proposition}[theorem]{Proposition}
\newtheorem{algorithm}{Algorithm}[section]

\newtheorem{definition}[theorem]{Definition}
\newtheorem{assumption}[theorem]{Assumption}

\theoremstyle{remark}
\newtheorem{remark}[theorem]{Remark}
\theoremstyle{definition}
\newtheorem{example}[theorem]{Example}

\newcommand{\N}{\mathbb{N}}
\newcommand{\R}{\mathbb{R}}

\newcommand\e{\varepsilon}

\newcommand\dist{\operatorname{dist}}
\newcommand\sym{\operatorname{sym}}
\newcommand\trace{\operatorname{tr}}
\newcommand\interieur{\operatorname{int}}
\newcommand{\SO}[1]{\operatorname{SO}(#1)}

\newcommand\id{{I}}
\newcommand\Id{I}
\newcommand\II{{\operatorname{I\!I}}}
\newcommand{\step}[1]{\medskip\noindent\textit{Step #1. }}

\newcommand{\wto}{\rightharpoonup}

\def\iso{{\mathrm{ iso}}}

\def\OmegaT{\Omega_{\mathrm{top}}}

\usepackage{cleveref} 

\begin{document}

\begin{center}
  \LARGE
  A nonlinear bending theory for nematic LCE plates
  \bigskip

  \normalsize
  S\"oren Bartels\footnote{bartels@mathematik.uni-freiburg.de}\textsuperscript{*}, Max Griehl\footnote{max.griehl@tu-dresden.de}\textsuperscript{\$},
  Stefan Neukamm\footnote{stefan.neukamm@tu-dresden.de}\textsuperscript{\$}, David Padilla-Garza\footnote{david.padilla-garza@tu-dresden.de}\textsuperscript{\$} and
  Christian Palus\footnote{christian.palus@mathematik.uni-freiburg.de}\textsuperscript{*} \par \bigskip

  \textsuperscript{*}Department of Applied Mathematics, University of Freiburg \par
  \textsuperscript{\$}Faculty of Mathematics, Technische Universit\"at Dresden\par \bigskip

  \today
\end{center}

\begin{abstract}
  In this paper, we study an elastic bilayer plate composed of a nematic liquid crystal elastomer in the top layer and a nonlinearly elastic material in the bottom layer.
  While the bottom layer is assumed to be stress-free in the flat reference configuration, the top layer features an eigenstrain that depends on the local liquid crystal orientation.
  As a consequence, the plate shows non-flat deformations in equilibrium with a geometry that non-trivially depends on the relative thickness and shape of the plate, material parameters, boundary conditions for the deformation, and anchorings of the liquid crystal orientation.
  We focus on thin plates in the bending regime and derive a two-dimensional bending model that combines a nonlinear bending energy for the deformation, with a surface Oseen-Frank energy for the director field that describes the local orientation of the liquid crystal elastomer.
  Both energies are nonlinearly coupled by means of a spontaneous curvature term that effectively describes the nematic-elastic coupling.
  We rigorously derive this model as a $\Gamma$-limit from three-dimensional, nonlinear elasticity.
  We also devise a new numerical algorithm to compute stationary points of the two-dimensional model.
  We conduct numerical experiments and present simulation results that illustrate the practical properties of the proposed scheme as well as the rich mechanical behavior of the system.
  \smallskip


  \textbf{MSC-2020:} 74B20; 76A15; 74K20; 65N30; 74-10.
\end{abstract}

\tableofcontents


\section{Introduction}
\begin{figure}
\centering
\begin{subfigure}{0.75\textwidth}
  \centering
  \includegraphics[width=\textwidth]{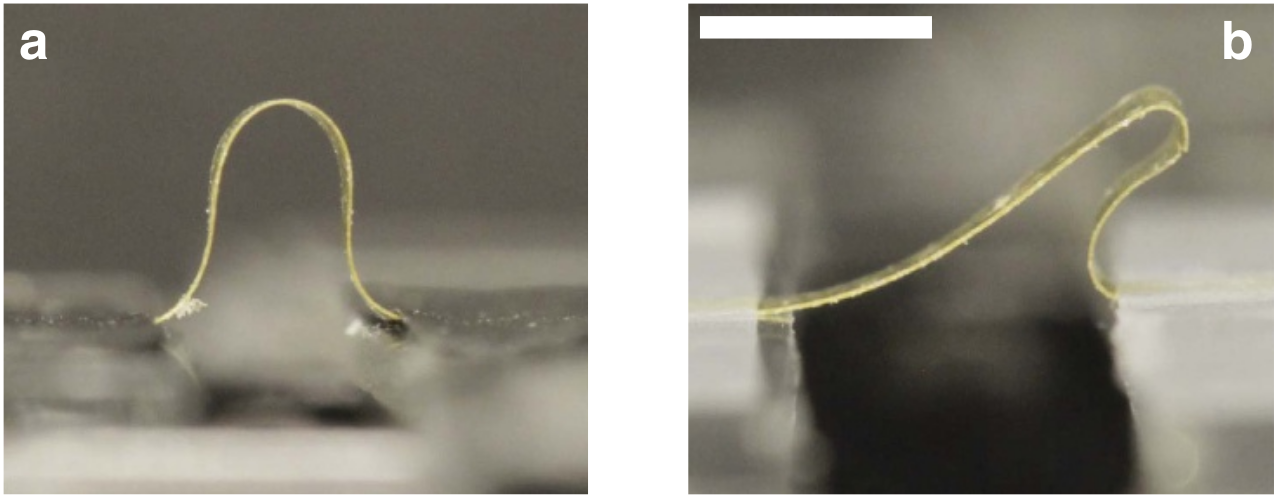}
  \caption{Bending/buckling behavior of a rectangular film with uniaxially aligned (left) vs. patterned (right) director. Source: Ware et al.~(~\cite{ware2016localized}), license:~CC BY 4.0 (\url{https://creativecommons.org/licenses/by/4.0/}).}
\end{subfigure}
\\[5mm]
\begin{subfigure}{0.75\textwidth}
  \centering
  \includegraphics[width=\textwidth]{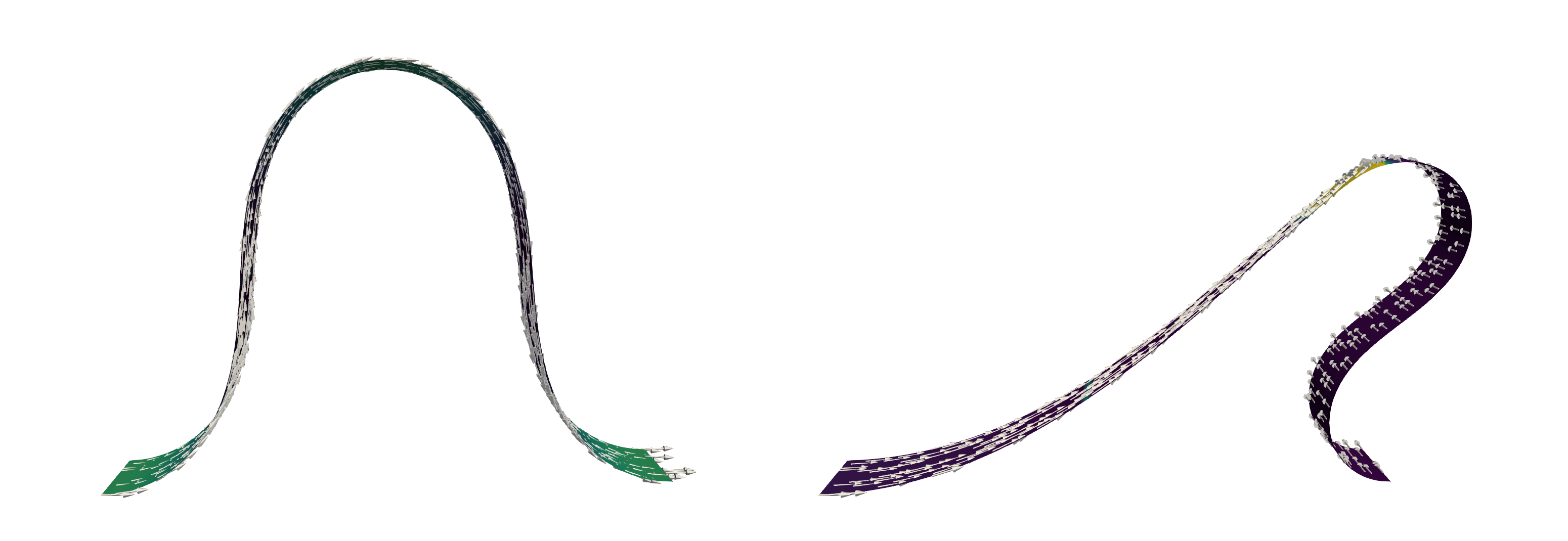}
  \caption{Simulation of an LC bilayer with either uniaxial (left) or patterned (right) prescribed director field. Parameters used:~$\mu=10^{-3}$, $\lambda=10^3$, $\bar{r}=1.5$, $\bar{\e}=8\times10^{-2}$ (cf.~Lemma~\ref{L:representations:hom:iso} and Section~\ref{S:modeling}).}
\end{subfigure}
\caption{Comparing findings from physical experiments with LCE simulations: Director anchoring may be used to control the buckling behavior.}
\label{fig:introfig}
\end{figure}

Liquid crystal elastomers (LCE) are solids made of liquid crystals (LC) incorporated into a polymer network.
In an isotropic phase (at high temperature) the LC are randomly oriented, while in the nematic phase (at low temperature) the LC show an orientational order and the material features a coupling between the entropic elasticity of the polymer network and the LC orientation.
LCE exhibit various interesting physical properties (e.g., soft elasticity \cite{B_Finkelman91, B_Kundler95}, or thermo-mechanical coupling \cite{B_broer11, B_Warner07}), and are considered in the design of active thin sheets that show a complex change of shape upon thermo-mechanical (or photo-mechanical) actuation, see \cite{B_White15} for a recent review and \cite{ware2016localized, MaJaZa18} for existing and future applications, cf.~Figure~\ref{fig:introfig}.
Understanding the mechanical properties of such thin films is the subject of current research in physics, mechanics, and mathematics; e.g.,~see \cite{B_conti2002soft, B_conti2002semisoft, B_Plucinsky18, B_plucinsky2018patterning, B_cesana2015effective,B_Agostiniani17plate, B_Agostiniani17platesoft, B_Agostiniani17platehetero, B_Agostiniani17ribbon} for a selection of contributions that focus on mathematical modeling and analysis of nematic LCE films.
\smallskip

In our paper, we consider a thin bilayer plate with a top layer consisting of a nematic LCE and a bottom layer consisting of a usual nonlinearly elastic material.
In view of the nematic-elastic coupling, the top layer features a non-vanishing eigenstrain, while the bottom layer is stress-free in flat reference configurations.
As a result, non-flat equilibrium shapes emerge with a geometry that is the result of a non-trivial interplay between the orientation of the LC and the deformation of the plate.
Various shapes have been observed experimentally, e.g.,~see~\cite{B_White15, SUTDT10,greco2012bending}.
In particular, films have been studied that ``bend'' in order to achieve equilibrium, e.g.,~\cite{SUTDT10}.
\smallskip

The goal of the present paper is to develop and analyze an effective, two-dimensional model for the nematic LCE bilayer plate, which is able to predict non-trivial bending deformations as energy minimizers and stationary points.
For this purpose, we consider a three-dimensional, geometrically nonlinear elasticity model for the bilayer plate that invokes both --- the elastic energy of the elastomer and the Frank energy of the LCE.
As a main result, we rigorously derive a two-dimensional bending model from the three-dimensional model as a $\Gamma$-limit; here, we focus on a regime where the thickness of the plate is comparable with the relaxation length-scale of the LC directors and we consider an energy scaling that is compatible with bending deformations.
Secondly, we devise a numerical algorithm to compute stationary points of the two-dimensional model and present simulation results that illustrate the rich mechanical behavior of the system.
\smallskip

Our simulations are able to reproduce characteristic material traits that are observed in physical experiments as shown in Figure~\ref{fig:introfig}, where the compression of an LCE strip leads to classical Euler buckling or non-Euler buckling.
In the physical experiment from~\cite{ware2016localized} as well as in the numerical simulation illustrated in the figure, the buckling behavior is determined by the alignment of the director, which is either uniaxial along the long direction of the strip, or is patterned, i.\,e. it is constantly parallel to the long direction of the strip in its left half, undergoes a rotation of $\pi/2$ around the surface normal in a small transition region close to the center of the strip, and then is constant again in most of its right half.
\smallskip

In contrast to earlier works, our model is not restricted to situations, where the orientation of the LC is prescribed.
We also analyze settings, where the orientation of the LC is not prescribed and thus allows for a soft elastic response.
We also consider partial constraints (e.g.,~restrictions to the tangential plane) for the LC orientation.
\smallskip

Let us present the main contributions of our paper with more detail:
We denote by $S\subseteq\R^2$ the mid-surface of the plate and by $0<h\ll 1$ its thickness.
We denote the rescaled domain of the plate by $\Omega\colonequals S\times(-\frac12,\frac12)$ and its top layer by $\OmegaT$.
As we shall explain in Section~\ref{S:modeling} an appropriate model for the bilayer plate is given by the (rescaled and non-dimensionalized) energy functional,
\begin{equation}\label{eq:energy}
  \begin{aligned}
    \mathcal E_h(y_h,n_h)\colonequals& \,\frac1{h^2}\int_{\Omega\setminus\OmegaT} W\big(x_3,\nabla\!_h y_h(x)\big)\,\mathrm{d} x+\frac1{h^2}\int_{\OmegaT}W\big(x_3,(L_h(n_h(x)))^{-\frac12}\nabla\!_hy_h(x)\big)\,\mathrm{d} x\\
    &\;+\bar\e^2\int_{\OmegaT}|\nabla\!_hn_h(x)\big(\nabla\!_h y_h(x)\big)^{-1}|^2\det(\nabla\!_h y_h(x))\,\mathrm{d} x,
  \end{aligned}
\end{equation}
where $y_h:\Omega\to\R^3$ denotes the deformation of the plate and $n_h:\OmegaT\to\mathbb S^2$ (with $\mathbb S^2$ being the unit sphere in $\R^3$) the director field (pulled back to the reference domain).
It describes the orientation of the LC.
Above, $\nabla\!_h\colonequals (\nabla',\frac1h\partial_3)\colonequals (\partial_1,\partial_2,\frac1h\partial_3)$ denotes the scaled gradient that comes from upscaling the domain to unit thickness.
The third integral describes the Oseen-Frank energy of the nematic LCE (in the simple form of a one-constant approximation). Here, $\bar\e^2> 0$ is a parameter that is related to the modulus of Frank elasticity, the shear modulus of the elastic material, and the physical thickness of the plate; see Section~\ref{S:modeling} for details.
The first and second integral in \eqref{eq:energy} describe the elastic energy stored in the deformed plate.
Here, $W$ denotes a stored energy function that describes a frame-indifferent, non-degenerate material with a stress-free reference state (see Assumption~\ref{ass:W} below for the precise assumptions), and
\begin{equation}\label{eq:steplength}
  L_h(n)\colonequals r_h^{-\frac13}(I+(r_h-1)n\otimes n),\qquad n\in\mathbb S^2,
\end{equation}
denotes the step-length tensor introduced by Bladon, Terentjev and Warner \cite{bladon1993transitions, B_Warner07} to model the nematic-elastic coupling. We note that the coupling is active in the nematic phase (when the  order parameter satisfies $r_h>1$) and enforces stretching (resp. compression) in the direction parallel to $n$ (resp.~in directions orthogonal to $n$). In the isotropic phase, which corresponds to $r_h=1$, the coupling is inactive.
Since we are interested in the bending behavior of the plate, and since bending modes have an energy per unit volume, that scales quadratically with the thickness of the plate, we consider the scaling factor $h^{-2}$ in front of the first two integrals in \eqref{eq:energy}, and we make the assumption that $r_h=1+h\bar r$ for some fixed constant $\bar r\in\R$.
\medskip

Theorem~\ref{T1}, which is the analytic main result of the paper, shows that in the zero-thickness limit $h\to 0$, the energy functional $\mathcal E_h$ $\Gamma$-converges to a two-dimensional model.
It takes the form of an energy $\mathcal E(y,n)$
defined for bending deformations
\begin{equation*}
  y\in H^2_{\iso}(S;\R^3)\colonequals \Big\{y\in H^2(S;\R^3)\,:\,(\nabla'y)^\top\nabla'y=I_{2\times 2}\Big\},
\end{equation*}
and director fields $n\in H^1(S;\mathbb S^2)$.
Here and below, $I_{d\times d}$ denotes the unit matrix in $\R^{d\times d}$.
The energy $\mathcal E$ couples a nonlinear plate model (quadratic in curvature) with a surface Oseen-Frank energy.
In the special case of a homogeneous material, the derived model takes the form
\begin{equation*}
\mathcal E(y,n)\colonequals \int_S Q_\mathrm{el}\left(\II_y+\frac34 \bar r\left(\frac 13 I-
\hat n'
\otimes
\hat n'
\right)\right)
+ \frac{3}{16}\bar r^2 Q_\mathrm{el}
\left(
	\frac 13 I -
	\hat n'
	\otimes
	\hat n'
\right)
+\frac{\bar \varepsilon^2}{2} |\nabla' n|^2\mathrm{d}x',
\end{equation*}
where
\begin{equation*}
  \II_y\colonequals \nabla'y^\top\;\nabla'(\partial_1 y\wedge\partial_2 y)
\quad
\text{and}
\quad
\hat n'\colonequals \nabla' y^\top n
\end{equation*}
denote the second fundamental form (in local coordinates) of the isometric immersion $y$ and the in-plane component of the director $n$ in local coordinates, respectively. (See Figure~\ref{Fig:2D} for the geometrical meaning of $n$ and $\hat n'$.)
Above, $Q_\mathrm{el}$ is a quadratic form obtained from  $W$ by a relaxation formula, which in the special case of a homogeneous and isotropic material can be explicitly written down, see Lemma~\ref{L:representations:hom:iso}.

The elastic energy of the bottom layer, the preferred curvature resulting from an inhomogeneous mismatch between the layers, and the Oseen-Frank term describing the crystalline order of the LCE layer yield competing contributions to the dimensionally reduced energy functional.
Characteristic properties of the three-dimensional model such as non-vanishing eigenstrain in the LCE layer are preserved in the two-dimensional setting, which itself exhibits its own characteristics: in many cases, the preferred curvature tensor $\frac{3}{4}\bar{r}(\frac{1}{3}I - \hat{n}'\otimes\hat{n}')$ cannot be attained by the second fundamental form of an isometric immersion of a flat reference configuration.
This leads to non-vanishing bending stresses.
In order to explore the model behavior by means of numerical simulations in Section~\ref{Sec:Freiburg}, we discretize the two-dimensional energy using discrete Kirchhoff triangular (DKT) elements for the deformation $y$ and Lagrange $P_1$ finite elements for the director $n$.
Functions in the corresponding finite element spaces are required to satisfy the isometry and unit-length constraints node-wise up to some small tolerance.
We approximate stationary points of the discretized energy via alternating between descent steps resulting from discrete $H^2$- and $H^1$-gradient flows for the deformation and the director field, respectively, using linearizations of the corresponding constraints.
Our numerical experiments confirm the rich behavior of the model: We investigate typical properties of the modeled LCE plates with an emphasis on the effects induced by prescribing various boundary conditions for the director, which in some situations have a significant influence on the energy landscape.
Most noteworthy, for a rectangular LCE plate clamped on two opposing sides, with certain combinations of coupling and order parameters $\bar{r}$ and $\bar{\varepsilon}$, imposing director boundary conditions may be used to select between a notable energy response or almost no response with respect to compression, effectively rendering the plate softer or more rigid.

\smallskip

We conclude the introduction by comparing our findings with previous results of the literature.
As all rigorous derivations of nonlinear bending models from 3d, our $\Gamma$-convergence result for $\mathcal E_h$ is based on the theory developed by Friesecke, James, M\"uller \cite{FJM02}.
In particular, if we consider $y_h\mapsto \mathcal E_h(y_h,n_h)$ with the director field $n_h$ being fixed, and if we neglect the nematic-elastic coupling, i.e.,~$r_h=1$, then we recover the celebrated result of \cite{FJM02}.
Furthermore, if we consider $y_h\mapsto \mathcal E_h(y_h,n_h)$ (with $n_h$ being still fixed) and $r_h=1+h\bar r$ with $\bar r\neq 0$, we end up with a model similar to $y_h\mapsto \int_\Omega W(\nabla\!_hy_h (I+hB)^{-1})\,\mathrm d x$ with a prestrain $(I+hB)$ that is of the order of the plate's thickness.
Problems of that type (with applications beyond nematic LCE plates) have been extensively studied in recent years, e.g., see \cite{Schm07a,SCHMIDT07,neukamm2010homogenization,BLS16,DPG20,LL20,BNS, BNPS}.
In this flavor, also bending models for nematic LCE plates with a prescribed and a fixed director $n$ have been derived, \cite{B_Agostiniani17plate,B_Agostiniani17platesoft,B_Agostiniani17platehetero,B_Aharoni14}.
In contrast to these works, in our model the director field $n$ is a free variable of the model. This leads to the presence of additional nonlinearities in our model; namely, the coupling term $L_h(n_h)$ in the second integral in \eqref{eq:energy} and the factor $(\nabla\!_hy_h)^{-1}\sqrt{\det\nabla\!_hy_h}$ in the third integral in \eqref{eq:energy}.
In our proof we shall first obtain an a priori estimate for $\nabla\!_hn_h$ in $L^p$ for $1<p<2$ that eventually implies that $\frac 1h (L_h(n_h)^{-\frac12}-I)$ strongly converges. In order to obtain the a priori estimate on $\nabla\!_hn_h$ (see Lemma~\ref{L:pf:apriori}) we need to assume a sufficiently strong (yet physical) growth of $W(F)$ for $\det F\to 0$ and for $|F|\rightarrow\infty$ (see Assumption~\ref{ass:W} below).
\smallskip

We also treat the case of clamped, affine boundary conditions for $y$ on flat parts of the lateral boundary $\partial S\times(-\frac12,\frac12)$, see Section~\ref{S:BC}.
We note that such boundary conditions have already been treated in the purely elastic case in \cite{FJM02} based on a Lipschitz truncation argument, which, however, does not extend to situations with prestrain.
Instead, we introduce a different argument that invokes an approximation result of independent interest:
we prove that the well-known density of $H^2_{\iso}(S;\R^3)\cap C^\infty(\overline S;\R^3)$ in $H^2_{\iso}(S;\R^3)$ (see~\cite{Pakz04,HornungApprox}) extends to the case with clamped, affine boundary conditions, see Proposition~\ref{P:approx}.
In Section~\ref{S:anchoring} we also analyze in detail the emergence of 2d-anchorings for the director field from 3d-anchorings.
This is partly motivated by \cite{Sternberg15}, where dimension reduction is carried out for non-deformable nematic plates in a Q-Tensor setting.
We consider both strong and weak anchorings on lateral parts of the boundary, as well as on the mid-surface.
In particular, the analysis of boundary- and anchoring conditions for $y$ and $n$ justifies the model problems that we consider in our numerical simulations in Section~\ref{subsec:numexp}.

In Section~\ref{subsec:discret} we devise a discretization of the two-dimensional energy based on discrete Kirchhoff triangular (DKT) elements for the deformation and Lagrange $P_1$ finite elements for the director that requires the unit-length and isometry constraints to be satisfied up to a small tolerance in the nodes of a triangulation.
Consequently, in Section~\ref{subsec:gradflow} we propose a discrete gradient flow scheme using linearizations of the respective constraints for approximating stationary points of the discretized energy.
The approach, which alternates between $H^1$- and $H^2$-gradient descent steps for director and deformation, essentially combines techniques developed in \cite{Bart05,Bart16} for harmonic maps and in \cite{Bart13a,Bart13b,BaBoNo17,BBMN18,Bart19-pre,BaPa21} for the isometric bending of prestrained bilayer plates.
Numerical experiments in Section~\ref{subsec:numexp} serve to illustrate both the behavior of the two-dimensional model as well as practical properties of the numerical scheme which we plan to analyze in a follow-up paper.
For isometric plate bending a similar $H^2$-gradient flow approach based on a discontinuous Galerkin discretization has been analyzed in~\cite{BoNoNt19}.
An application of that method to the case of (normal elastic) bilayers has been investigated experimentally in~\cite{BoNoNt20}.
For harmonic maps, an alternative approach that relies on a conforming discretization of the target manifold in the form of geodesic finite elements is presented in~\cite{Sand12,Sand16}.
Finite element methods for the Ericksen model of nematic LC with variable degree of orientation, i.\,e., in the case that $\bar{\varepsilon} = \bar{\varepsilon}(x)$ is not constant, have been considered in~\cite{NoWaZh17,Walk20}.

\paragraph{Structure of the paper.}
The paper is structured as follows:
In Section~\ref{S:setup} we introduce the precise setup, present the definition of the two-dimensional limiting functional and state the main $\Gamma$-convergence result.
Section~\ref{S:BC} is devoted to the discussion of boundary conditions for the deformation.
Anchoring for $n$ is discussed in Section~\ref{S:anchoring}.
Section~\ref{S:modeling} discusses the modeling and non-dimensionalization of the energy. All proofs are postponed to Section~\ref{S:proofs}.
In Section~\ref{Sec:Freiburg} we formulate the numerical algorithm and present simulation results.

\paragraph{Notation.}
In the paper we shall mostly use standard notation whose meaning is evident from the context. Furthermore, the following less standard notation is used:
\begin{itemize}
\item For $x\in\R^3$ we write $x=(x',x_3)$, and call $x'$ the in-plane components and $x_3$ the out-of-plane component of $x$.
\item $\nabla'\colonequals (\partial_1,\partial_2)$ denotes the in-plane gradient, and $\nabla\!_h\colonequals (\nabla',\frac1h\partial_3)$ denotes the scaled gradient.
\item Both of the above conventions are dropped in Section~\ref{Sec:Freiburg} where, to simplify notation, 2d-coordinates and gradients are denoted with $x$ and $\nabla$, respectively.
\item $I$ denotes the identity matrix on $\mathbb{R}^{d \times d}$; sometimes we write $I_{d \times d}$ to avoid ambiguities.
\item $H^2_\mathrm{iso}(S;\mathbb{R}^3)$ denotes the space of Sobolev isometries, i.e.,~$y\in H^2(S;\R^3)$ satisfying the pointwise isometry constraint $\nabla'y^\top\nabla'y=I$.
For $y\in H^2_{\iso}(S;\R^3)$ we write $b_y\colonequals \partial_1 y\wedge \partial_2 y$, $R_y\colonequals (\partial_1 y,\partial_2y,b_y)$ and $\II_y\colonequals \nabla'y^\top\nabla'b_y$.
The latter is the negative of the second fundamental form of the surface parametrized by $y$ represented in local coordinates.
\item $E_h(y)=\frac{\sqrt{(\nabla\!_hy)^\top\nabla\!_h y}-I}{h}$ denotes the scaled nonlinear strain associated with the deformation $y\in H^1(\Omega;\R^3)$.
\item $|\cdot|$ denotes the Euclidean norm on $\R^d$ and the Frobenius norm on matrix spaces.
For $A,B\in\mathbb R^{d\times d}$, we consdier the scalar product $A\cdot B\colonequals \trace(A^\top B)$.
\item For $A\in\R^{2\times 2}$ we denote by $\iota(A)$ the $3\!\!\times\!\!3$-matrix that is zero except for the upper-left $2\!\!\times\!\!2$-block, which is given by $A$.
\item For $A\in\R^{2\times 2}$ the inequality $A\geq 0$ has to be understood in the sense of quadratic forms, i.e.,~$Aa\cdot a\geq |a|^2$ for all $a\in\R^2$.
\end{itemize}

\section{Setup and statement of analytical main results}\label{S:setup}

\subsection{The rescaled three-dimensional plate model and assumptions}
\label{Sec:3d:assumptions}
We consider a bilayer plate with (rescaled) reference configuration $\Omega\colonequals S\times(-\frac12,\frac12)$ where
\begin{equation}\label{ass:domain}
  S\subseteq\R^2\text{ is an open, bounded and convex domain with piece-wise } C^1\text{-boundary}.
\end{equation}
Convexity is assumed for simplicity and the extension to non-convex $S$ is possible, see~Remark~\ref{R:nonconvex} below.
We assume that the top layer of the plate  $\OmegaT\colonequals S\times(0,\frac12)$ is occupied by a nematic LCE, while the bottom layer of the plate $\Omega\setminus\OmegaT$ consists of a standard nonlinearly elastic material.
As we shall explain in Section~\ref{S:modeling} we can model such a plate with thickness $h>0$ by the rescaled, non-dimensionalized energy functional
\begin{equation*}
  \mathcal E_h:H^1(\Omega;\R^3)\times H^1(\OmegaT;\mathbb S^2)\to[0,\infty],
\end{equation*}
with  $\mathcal E_h(y_h,n_h)$ being defined by \eqref{eq:energy}.
As already mentioned in the introduction, throughout the paper we assume $\bar\varepsilon>0$ and $r_h\colonequals 1+h\bar r$ (with $\bar r\in\R$ fixed), where $\bar \varepsilon$ denotes the parameter in front of the Oseen-Frank energy and $r_h$ the order parameter in the definition of the step-length tensor $L_h(n_h)$ introduced in \eqref{eq:steplength}.
Thus, the deviation of $L_h(n_h)^{-\frac12}$ from identity is given (up to a scaling with $h$) by the nonlinear map $B_h:\mathbb S^2\to\R^{3\times 3}$,
\begin{equation}\label{eq:Bh}
  B_h(n)\colonequals \frac1 h(L_h(n)^{-\frac12}-I)=\frac{(1+h\bar r)^{\frac16}-1}{h}I + (1+h\bar r)^\frac16\frac{(1+h\bar r)^{-\frac12}-1}{h}n\otimes n.
\end{equation}
Note that $B_h$ is only well-defined for $h>0$ with $|h\bar r|<1$.
We tacitly assume this condition from now on.

We suppose that the stored energy function $W$ satisfies the following properties:

\begin{assumption}[Stored energy function]\label{ass:W}
  Let $C_W>0$,  $q_W>4$, and let $r_W:[0,\infty)\to[0,\infty]$ denote a monotone function satisfying $\lim_{\delta\to0}r_W(\delta)=0$. We assume that $W:(-\frac12,\frac12)\times\R^{3\times 3}\to[0,\infty]$ is a Borel function such that for all $x_3\in(-\frac12,\frac12)$ we have:
\begin{enumerate}[label  = (W\arabic*)]
\item \label{item:ass:W:frame:ind}
	(frame indifference). $W(x_3,RF)=W(x_3,F)$ for all $R\in \SO 3$ and $F\in\R^{3\times 3}$.
\item \label{item:ass:W:non:degen:nat}
	(non-degeneracy and natural state).
 \begin{align*}
  W(x_3,F)&\geq \frac{1}{C_W} \dist^2(F,\SO 3)\qquad\mbox{for all $F\in\R^{3\times 3}$,}\\
  W(x_3,F)&\leq C_W \dist^2(F,\SO 3)\qquad\mbox{for all $F\in\R^{3\times 3}$ with $\dist^2(F,\SO 3)\leq \frac{1}{C_W}$.}
 \end{align*}
\item \label{item:ass:W:quadratic:exp}
	(quadratic expansion at $\SO3$). There exists a quadratic form $Q(x_3,\cdot):\R^{3\times 3}\to\R$ such that
  \begin{equation*}
    |W(x_3,I_{3\times 3}+G)-Q(x_3,G)|\leq |G|^2r_W(|G|)\qquad\mbox{for all $G\in\R^{3\times 3}$.}
  \end{equation*}
\item \label{item:ass:W:growth:cond}
	(growth condition).
	For all $F\in\R^{3\times 3}$,
  \begin{equation*}
    W(x_3,F)\geq
    \begin{cases}
      \frac{1}{C_W}\max\{|F|^{q_W},\det(F)^{-\frac{q_W}{2}}\}-C_W &\det F>0,\\
      +\infty&\text{else.}
    \end{cases}
  \end{equation*}
\end{enumerate}
\end{assumption}

\begin{remark}
\label{rem:W}
  Assumptions \ref{item:ass:W:frame:ind}, \ref{item:ass:W:non:degen:nat}, \ref{item:ass:W:quadratic:exp} are standard assumptions in the context of the derivation of plate theories from 3d-nonlinear elasticity. In particular, these assumptions allow us to linearize the material law at $\SO 3$. The quadratic form in assumption \ref{item:ass:W:quadratic:exp} is given by the Hessian of $W$ at identity, i.e.,
  \begin{equation*}
    Q(x_3,G)=\frac12 \frac{\partial^2W(x_3,I)}{\partial F\partial F}G\cdot G.
  \end{equation*}
  As a consequence of \ref{item:ass:W:frame:ind}, \ref{item:ass:W:non:degen:nat} and \ref{item:ass:W:quadratic:exp}, $Q$ satisfies
  \begin{equation}\label{eq:Q:1}
    \forall G\in\R^{3\times 3}\,:\qquad\frac{1}{C_W}|\sym G|^2\leq Q(x_3,G)\leq C_W|\sym G|^2.
  \end{equation}
  The growth condition \ref{item:ass:W:growth:cond} is less standard, however, in accordance with physical growth assumptions: it enforces a blow up of the elastic energy for compressions of the material to zero volume and for large changes of  area, and excludes (local) interpenetration. From the analytic perspective, we use \ref{item:ass:W:growth:cond} in Lemma~\ref{L:pf:apriori} to control the expression $F^{-1}\sqrt{\det F}$, which appears in the formulation of the Oseen-Frank energy, cf.~\eqref{eq:energy}.
\end{remark}

\subsection{The effective bending plate model and $\Gamma$-convergence}
\label{Sec:2d:assumptions:and:main:results}
As we shall see, in the limit $h\to 0$, we obtain as a $\Gamma$-limit an energy functional $\mathcal E(y,n)$ defined for a bending deformation $y:S\to\R^3$ and a director field $n:S\to\mathbb S^2$ in the domain
\begin{equation*}
  \mathcal A_2\colonequals \{(y,n)\,:\, y\in H^2_{\iso}(S;\R^3),\,n\in H^1(S;\mathbb S^2)\}.
\end{equation*}
The limiting functional $\mathcal E$ can be written as the sum of an elastic energy and surface Oseen-Frank energy,
\begin{equation}\label{eq:energy2d:total}
  \mathcal E(y,n)\colonequals \mathcal E_\mathrm{el}(y,n)+\bar\e^2\mathcal E_\mathrm{OF}(n),
\end{equation}
where
\begin{equation*}
\mathcal E_{OF}(n)\colonequals \frac12\int_S|\nabla'n|^2\,\mathrm{d} x'.
\end{equation*}
(Recall the notations $x'\colonequals (x_1,x_2)$ and $\nabla'\colonequals (\partial_1,\partial_2)$.)
The elastic energy is defined as
\begin{equation*}
  \begin{aligned}
    \mathcal E_\mathrm{el}(y,n)\colonequals  &\,\int_{S} Q_\mathrm{el}\Big(\II_y+\bar r\mathbb B\big(\tfrac13 I_{2\times 2}-\nabla'y^\top n\otimes \nabla'y^\top n\big)\Big)\\
    &\qquad +\bar r^2E_\mathrm{res}\big(\tfrac13 I_{2\times 2}-\nabla'y^\top n\otimes \nabla'y^\top n\big)\,\mathrm{d} x'.
  \end{aligned}
\end{equation*}

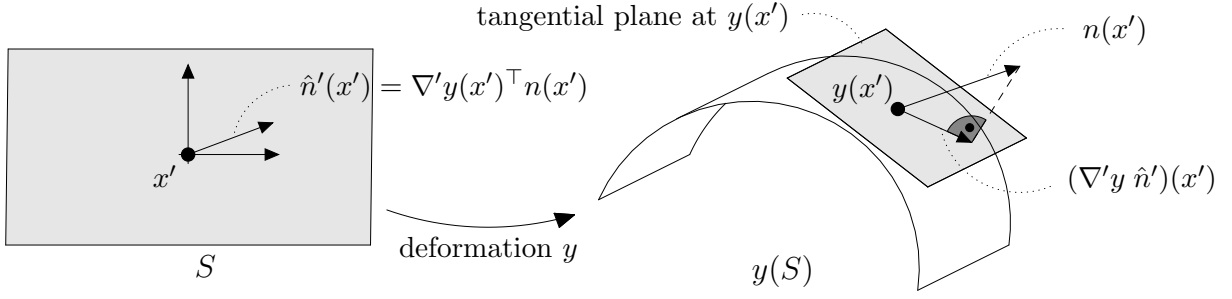
\begin{figure}
\begin{tikzpicture}[line cap=round,line join=round,>=triangle 45,x=.6cm,y=.6cm]
\clip(-6.3,1.0) rectangle (20.5,7.8);
\draw [shift={(15.184470297609847,4.264620292068637)},line width=.2pt,fill=black,fill opacity=0.42] (0,0) -- (57.99461679191625:0.5901046526559983) arc (57.99461679191625:155.4398927656544:0.5901046526559983) -- cycle;
\fill[line width=.2pt,fill=black,fill opacity=0.1] (11.137797943838862,5.689752729701015) -- (14.220314379355658,3.282892507625115) -- (16.383492799621237,4.364481717757904) -- (13.30097636410444,6.771341939833804) -- cycle;
\fill[line width=.2pt,fill=black,fill opacity=0.1] (-5.946537434084586,6.331717926057335) -- (2.0534625659154138,6.331717926057335) -- (2,2) -- (-6,2) -- cycle;
\fill[line width=.2pt,fill=black,fill opacity=1] (6.0259726284292165,2.701044965500647) -- (6.138630859067428,2.389068326810221) -- (6.473948526337566,2.6691530347436943) -- cycle;
\draw [->,line width=.2pt] (-2,4) -- (-2,6);
\draw [->,line width=.2pt] (-2,4) -- (0,4);
\draw [->,line width=.2pt] (-2,4) -- (-0.1288755144223903,4.7063236931265475);
\draw (0.22467416865026266,6.1) node[anchor=north west] {$\hat n'(x')=\nabla'y(x')^\top n(x')$};
\draw [line width=.2pt] (-2,4)-- (-2.2,4);
\draw [line width=.2pt] (-2,4)-- (-2,3.8);
\draw (-3,4) node[anchor=north west] {$x'$};
\draw [shift={(10.357142857142858,1.5)},line width=.2pt]  plot[domain=-0.13640260440094742:2.721395939980518,variable=\t]({1*3.6770107646382146*cos(\t r)+0*3.6770107646382146*sin(\t r)},{0*3.6770107646382146*cos(\t r)+1*3.6770107646382146*sin(\t r)});
\draw [shift={(12.357142857142831,2.5)},line width=.2pt]  plot[domain=2.3433960894852786:2.721395939980517,variable=\t]({1*3.677010764638187*cos(\t r)+0*3.677010764638187*sin(\t r)},{0*3.677010764638187*cos(\t r)+1*3.677010764638187*sin(\t r)});
\draw [shift={(12.357142857142858,2.5)},line width=.2pt]  plot[domain=-0.13640260440094742:2.0344439357957027,variable=\t]({1*3.6770107646382146*cos(\t r)+0*3.6770107646382146*sin(\t r)},{0*3.6770107646382146*cos(\t r)+1*3.6770107646382146*sin(\t r)});
\draw [line width=.2pt] (8.712733652396953,4.788818409491811)-- (10.712733652396953,5.788818409491811);
\draw [line width=.2pt] (7,3)-- (9,4);
\draw [line width=.2pt] (14,1)-- (16,2);
\draw [line width=.2pt] (11.137797943838862,5.689752729701015)-- (14.220314379355658,3.282892507625115);
\draw [line width=.2pt] (11.137797943838862,5.689752729701015)-- (13.30097636410444,6.771341939833804);
\draw [line width=.2pt] (16.383492799621237,4.364481717757904)-- (13.30097636410444,6.771341939833804);
\draw [line width=.2pt] (16.383492799621237,4.364481717757904)-- (14.220314379355658,3.282892507625115);
\draw (10.138432333271036,2) node[anchor=north west] {\large $y(S)$};
\draw (4.1,7.588413358146278) node[anchor=north west] {$\text{tangential plane at } y(x')$};
\draw [->,line width=.2pt] (13.563616887738325,5.005340264960481) -- (16.240123955115322,5.95366614407738);
\draw (17.396719560939815,7.3) node[anchor=north west] {$n(x')$};
\draw (11.9,6) node[anchor=north west] {$y(x')$};
\draw [line width=.2pt,dashed] (16.240123955115322,5.95366614407738)-- (15.184470297609847,4.264620292068637);
\draw [->,line width=.2pt] (13.563616887738325,5.005340264960481) -- (15.184470297609847,4.264620292068637);
\draw (17.042656769346216,4.0) node[anchor=north west] {$(\nabla'y \;\hat n') (x')$};
\draw [line width=.2pt] (11.137797943838862,5.689752729701015)-- (14.220314379355658,3.282892507625115);
\draw [line width=.2pt] (14.220314379355658,3.282892507625115)-- (16.383492799621237,4.364481717757904);
\draw [line width=.2pt] (16.383492799621237,4.364481717757904)-- (13.30097636410444,6.771341939833804);
\draw [line width=.2pt] (13.30097636410444,6.771341939833804)-- (11.137797943838862,5.689752729701015);
\draw [shift={(11.716785641188585,5.936022770565056)},line width=.2pt,dotted]  plot[domain=0.4630710472330943:1.8252666201029544,variable=\t]({1*1.1817294393547422*cos(\t r)+0*1.1817294393547422*sin(\t r)},{0*1.1817294393547422*cos(\t r)+1*1.1817294393547422*sin(\t r)});
\draw [shift={(17.103668966363088,5.0795089928336665)},line width=.2pt,dotted]  plot[domain=1.5199648360304359:2.770474767655488,variable=\t]({1*1.7032858853370845*cos(\t r)+0*1.7032858853370845*sin(\t r)},{0*1.7032858853370845*cos(\t r)+1*1.7032858853370845*sin(\t r)});
\draw [shift={(16.24003649400355,4.92926519179235)},line width=.2pt,dotted]  plot[domain=3.3443802452398206:5.154814206889009,variable=\t]({1*1.7673018482908822*cos(\t r)+0*1.7673018482908822*sin(\t r)},{0*1.7673018482908822*cos(\t r)+1*1.7673018482908822*sin(\t r)});
\draw [line width=.2pt] (-5.946537434084586,6.331717926057335)-- (2.0534625659154138,6.331717926057335);
\draw [line width=.2pt] (2.0534625659154138,6.331717926057335)-- (2,2);
\draw [line width=.2pt] (2,2)-- (-6,2);
\draw [line width=.2pt] (-6,2)-- (-5.946537434084586,6.331717926057335);
\draw (-2.076733976708131,2) node[anchor=north west] {\large $S$};
\draw [shift={(0.3900288583242757,4.117747904729774)},line width=.2pt,dotted]  plot[domain=1.6914836376757174:2.960237950799003,variable=\t]({1*1.4219958308805274*cos(\t r)+0*1.4219958308805274*sin(\t r)},{0*1.4219958308805274*cos(\t r)+1*1.4219958308805274*sin(\t r)});
\draw [shift={(4.579406560145415,9.027638230383404)},line width=.2pt]  plot[domain=4.370873083458915:5.001968502887352,variable=\t]({1*6.634728573561437*cos(\t r)+0*6.634728573561437*sin(\t r)},{0*6.634728573561437*cos(\t r)+1*6.634728573561437*sin(\t r)});
\draw [line width=.2pt] (6.473948526337566,2.6691530347436943)-- (6.0259726284292165,2.701044965500647);
\draw [line width=.2pt] (6.473948526337566,2.6691530347436943)-- (6.138630859067428,2.389068326810221);
\draw [line width=.2pt] (6.0259726284292165,2.701044965500647)-- (6.138630859067428,2.389068326810221);
\draw [line width=.2pt] (6.138630859067428,2.389068326810221)-- (6.473948526337566,2.6691530347436943);
\draw [line width=.2pt] (6.473948526337566,2.6691530347436943)-- (6.0259726284292165,2.701044965500647);
\draw (2.4,2.415162569862055) node[anchor=north west] {$\text{deformation } y$};
\begin{scriptsize}
\draw [fill=black] (13.563616887738325,5.005340264960481) circle (2.5pt); 
\draw [fill=black] (15.139621802462381,4.593846916107125) circle (1.5pt); 
\draw [fill=black] (-2,4) circle (2.5pt); 
\end{scriptsize}
\end{tikzpicture}
\caption{For each material point $x'\in S$, the director field's value can be expressed in global coordinates via $n(x')$ and in local coordinates via $\hat n(x')=R_y(x')^\top n(x')$.
The latter describes the director field with respect to the tangential plane and the surface normal in $y(x')$.
The value $\nabla'y(x')^\top n(x')=\hat n'(x')$ represents the tangential components.}
\label{Fig:2D}
\end{figure}

Note that $\nabla'y^\top n$ describes the in-plane projection of the director field $n$ in local coordinates, i.e., with respect to the tangential plane's basis $(\partial_1y,\partial_2y)$.
(See Figure \ref{Fig:2D}.)

Furthermore,
\begin{itemize}
\item
	$Q_\mathrm{el}:\R^{2\times 2}_{\mathrm{sym}}\to[0,\infty)$ denotes a quadratic form that describes the bending stiffness of the plate.
	It is positive definite on $\R^{2\times 2}_{\sym}$ and determined by the relaxation formula of Definition~\ref{def:QEB} from the quadratic form $Q$ in the expansion of $W$, cf.~\ref{item:ass:W:quadratic:exp}.
\item
	$E_\mathrm{res}$ is a positive definite quadratic form on $\R^{3\times 3}_{\sym}$ and describes the local residual energy coming from the LCE Eigenstrain.
	It is given in Definition~\ref{def:QEB} as a function of $Q$.
\item $\mathbb B:\R^{2\times 2}_{\sym}\to \R^{2\times 2}_{\sym}$ is a linear map that is related to the linearization of $B_h$.
	Its precise form, see Definition~\ref{def:QEB}, depends on $Q$.
\end{itemize}
The analytical main result of the paper is the following theorem, which states that $\mathcal E$ emerges as a rigorous $\Gamma$-limit of $\mathcal E_h$ for $h\to 0$:
\begin{theorem}[$\Gamma$-convergence]\label{T1}
  Let Assumption~\ref{ass:W} and \eqref{ass:domain} be satisfied.
  \begin{enumerate}[(a)]
  \item (Compactness). \label{item:T1:compactness}
    Let $(y_h,n_h)\subseteq H^1(\Omega;\R^3)\times H^1(\OmegaT;\mathbb S^2)$ satisfy
    \begin{equation}
      \label{EqFiniteEnergy}
      \limsup_{h\rightarrow 0}\mathcal E_h(y_h ,n_h )<\infty.
    \end{equation}
    Then there exist $(y,n)\in\mathcal A_2$, and $d\in L^2(\OmegaT;\R^3)$ such that (up to a subsequence),
    \begin{alignat}{2}
      \label{KompaktheitKonvergenzDef}
      (y_h-\fint_\Omega y_h,\nabla\!_h y_h)& \rightarrow(y,R_y) &&\qquad\text{strongly in }L^2(\Omega),\\
      \label{KompaktheitKonvergenzDir}
      n_h &\to n &&\qquad\text{strongly in }L^2(\Omega),\\
      \label{KompaktheitKonvergenzDir2}
      \nabla\!_hn_h &\rightharpoonup (\nabla'n,d) &&\qquad\text{weakly in }L^{p}(\OmegaT)\text{ with }p\colonequals \tfrac{2}{1+4/{q_W}}>1.
    \end{alignat}

  \item (Lower bound). \label{item:T1:lower:bound}
  Let $(y_h,n_h)\subseteq H^1(\Omega;\R^3)\times H^1(\OmegaT;\mathbb S^2)$ and $(y,n)\in\mathcal A_2$.
  Suppose $(y_h,n_h)\to (y,n)$ strongly in $L^2(\Omega)\times L^2(\OmegaT)$.
    Then,
    \begin{equation*} 
      \liminf_{h\rightarrow 0}\mathcal E_h(y_h ,n_h )\geq \mathcal E(y,n).
    \end{equation*}
  \item (Upper bound). \label{item:T1:upper:bound}
    For all $(y,n)\in\mathcal A_2$ there exists a sequence $(y_h,n_h)\subseteq H^1(\Omega;\R^3)\times H^1(\OmegaT;\mathbb S^2)$ with $(y_h,n_h)\to (y,n)$ strongly in $L^2(\Omega)\times L^2(\OmegaT)$ and
    \begin{equation}\label{eq:recov}
      \lim\limits_{h\to 0}\mathcal E_h(y_h,n_h)=\mathcal E(y,n).
    \end{equation}
  \end{enumerate}
\end{theorem}
(See Sections \ref{Sec:Proof:Compactness}, \ref{Sec:Proof:Lower:Bound}, \ref{Sec:Proof:Upper:Bound} for the proof.)\medskip

In the special case of a homogeneous, isotropic material (the situation that we shall investigate numerically in Section~\ref{Sec:Freiburg}) the energy densities $Q_\mathrm{el}$, $E_\mathrm{res}$ and the map $\mathbb B$ are given by explicit expressions:
\begin{lemma}[The homogeneous and homogeneous, isotropic case]
\label{L:representations:hom:iso}
  Let Assumption~\ref{ass:W} be sat\-is\-fied.
  Additionally assume that $W$ is independent of $x_3$. Then we have for all $A,U\in\R^{2\times 2}_{\mathrm{sym}}$,
  \begin{align}
  	\label{eq:repr:hom:iso:el}
    Q_\mathrm{el}(A)&=\frac1{12}\min_{d\in\R^3}Q\Big(\iota(A)+\sym(d\otimes e_3)\Big),\\
  	\label{eq:repr:hom:iso:Eres}
    E_\mathrm{res}(U)&=\frac{3}{16}Q_\mathrm{el}\big(U\big),\\
  	\label{eq:repr:hom:iso:B}
    \mathbb B(U)&=\frac{3}{4}U.
  \end{align}
  In addition, if the material is isotropic, that is,
  \begin{equation}\label{eq:Q:hom:iso}
    \forall A\in\R^{3\times 3}:\qquad Q(A)=\frac\lambda2(\operatorname{tr}A)^2+\mu|\sym A|^2,
  \end{equation}
  then \eqref{eq:repr:hom:iso:el} and \eqref{eq:repr:hom:iso:Eres} further simplify to
  \begin{align*}
    Q_\mathrm{el}(A)&=\frac{1}{12}\Big(\frac{\mu\lambda}{\lambda+2\mu}(\operatorname{tr}A)^2+\mu|A|^2\Big),\\
    E_\mathrm{res}(U)&=\frac{1}{64}\Big(\frac{\mu\lambda}{\lambda+2\mu}(\operatorname{tr} U)^2+\mu|U|^2\Big).
  \end{align*}
\end{lemma}
(See Section \ref{Section:Representation} for the proof.)\medskip

In the general case these quantities are defined as follows:
\begin{definition}[Definition of $Q_\mathrm{el}, E_\mathrm{res}$ and $\mathbb B$]\label{def:QEB}
  For $A,U\in\R^{2\times 2}_{\sym}$ and $x_3\in(-\frac12,\frac12)$ define
  \begin{align*}
    Q_2(x_3,A)&\colonequals \min_{d\in\R^{3}}Q(x_3,\iota(A)+\sym(d\otimes e_3)),\\
    Q_\mathrm{el}(A)&\colonequals \min_{M\in\R^{2\times 2}_{\sym}}\int_{-\frac12}^\frac12 Q_2(x_3,x_3 A+M)\,\mathrm{d} x_3,\\
    E_\mathrm{res}(U)&\colonequals \min_{A,M\in\R^{2\times 2}_{\sym}}\int_{-\frac12}^\frac12 Q_2\Big(x_3,x_3 A+M+\mathbf 1(x_3>0)\tfrac{1}{2} U\Big)\,\mathrm{d} x_3,
  \end{align*}
  where $\mathbf 1(x_3>0)$ denotes the indicator function of $(0,\tfrac12)$.
  For the definition of $\mathbb B$, we consider the matrices
  \begin{equation*}
    G_1=e_1\otimes e_1,\qquad G_2\colonequals e_2\otimes e_2,\qquad G_3\colonequals \frac1{\sqrt 2}(e_1\otimes e_2+e_2\otimes e_1),
  \end{equation*}
  which form an orthonormal basis of $\R^{2\times 2}_{\sym}$. For $i=1,2,3$, we denote by $A_i,M_i\in\R^{2\times 2}_{\sym}$ the unique minimizers of the strictly convex, quadratic function
  \begin{equation*}
    (\R^{2\times 2}_{\sym}\times\R^{2\times 2}_{\sym})\ni (A,M)\mapsto\int_{-\frac12}^\frac12 Q_2\Big(x_3,\mathbf 1(x_3>0)\tfrac{1}{2}G_i-(x_3 A+M)\Big)\,\mathrm{d} x_3.
  \end{equation*}
  We define $\mathbb B:\R^{2\times 2}_{\sym}\to\R^{2\times 2}_{\sym}$ via
  \begin{equation*}
    \mathbb B(U)\colonequals \sum_{i=1}^3\Big(U\cdot G_i\Big) A_i.
  \end{equation*}
  \qed
\end{definition}
The following structural properties are satisfied by $Q_{\mathrm{el}},E_{\mathrm{res}},\mathbb B$:
\begin{lemma}[Properties of $Q_{\mathrm{el}},E_{\mathrm{red}},\mathbb B$]
  \label{L:properties}
  Let $W$ satisfy Assumption~\ref{ass:W}, and let $Q_\mathrm{el}$, $E_\mathrm{res}$ and $\mathbb B$ be defined as in Definition~\ref{def:QEB}.
  Then,
  \begin{enumerate}[(a)]
  \item $Q_\mathrm{el}:\R^{2\times 2}_{\sym}\to[0,\infty)$ is quadratic and satisfies $\frac{1}{12\,C_W}|A|^2\leq Q_\mathrm{el}(A)\leq \frac{C_W}{12} |A|^2$.
  \item $E_\mathrm{res}:\R^{2\times 2}_{\sym}\to[0,\infty)$ is quadratic and satisfies $\frac{1}{64 C_W}|U|^2\leq E_\mathrm{res}(U)\leq \frac{C_W}{64}|U|^2$.
  \item $\mathbb B$ is linear and there is $C>0$ (depending only on $C_W$) such that $|B(U)|\leq C|U|$.
  \end{enumerate}
\end{lemma}
(See Section \ref{Sec:Proof:Properties:of:Limit} for the proof.)\medskip

As a consequence of Lemma \ref{L:properties}, the limiting energy features good compactness and continuity properties:
\begin{lemma}[Compactness and continuity of the limiting energy]\label{L:cont}
  Let Assumption~\ref{ass:W} be satisfied and let $S\subseteq\R^2$ be a Lipschitz domain. We extend $\mathcal E$ to a functional defined on $L^2(S;\R^3)\times L^2(S;\R^3)$ by setting $\mathcal E(y,n)\colonequals \infty$ for $(y,n)\not\in\mathcal A_2$. Then:
  \begin{enumerate}[(a)]
  \item (A priori estimates).   There exists $C=C(S,C_W, \overline{r})>0$ such that for all $(y,n)\in\mathcal A_2$ we have
    \begin{equation*}
      \int_S|\nabla'\nabla'y|^2\,\mathrm dx'\leq C(\mathcal E_\mathrm{el}(y,n)+1).
    \end{equation*}
  \item (Compactness).
  Sublevels of $\mathcal E$ of the form $\{(y,n)\in\mathcal A_2\,:\,\mathcal E(y,n)\leq C\text{ and }\|y\|_{L^2(S;\R^3)}\leq C\}$ with $C>0$ are compact in $\mathcal A_2$ considered with both, the strong topology of $L^2(S;\R^3)\times L^2(S;\R^3)$ and the weak topology of $H^2(S;\R^3)\times H^1(S;\R^2)$.
  \item (Lower semicontinuity). Let  $(y^k,n^k), (y,n)\in\mathcal A_2$ and suppose that $(y^k,n^k)\to (y,n)$ strongly in $L^2(\Omega)\times L^2(\OmegaT)$.
  Then
    \begin{align*}
      &\liminf\limits_{k\to\infty}\mathcal E_{OF}(n^k)\geq \mathcal E_{OF}(n),\\
      &\liminf\limits_{k\to\infty}\mathcal E_\mathrm{el}(y^k,n^k)\geq \mathcal E_\mathrm{el}(y,n),\text{ provided }\limsup\limits_{k\to\infty}\mathcal E_{OF}(n^k)<\infty.
    \end{align*}
    In particular,
    \begin{equation*}
      \liminf\limits_{k\to\infty}\mathcal E(y^k,n^k)\geq \mathcal E(y,n),
    \end{equation*}
  \item (Continuity). Let $(y^k,n^k),(y,n)\in\mathcal A_2$ and suppose that $y^k\to y$ strongly in $H^2(S;\mathbb R^3)$.
  Then
    \begin{alignat*}{2}
      &\mathcal E_\mathrm{el}(y^k,n^k)\to\mathcal E_\mathrm{el}(y,n)\qquad&&\text{if }n^k\wto n\text{ weakly in }H^1(\OmegaT),\\
      &\mathcal E_{OF}(n^k)\to\mathcal E_{OF}(n)\qquad&&\text{if }n^k\to n\text{ strongly in }H^1(\OmegaT).
    \end{alignat*}
  \end{enumerate}
\end{lemma}
(See Section \ref{Sec:Proof:Properties:of:Limit} for the proof.)\medskip

We note that in view of the previous lemma, the direct method can be applied to $\mathcal E$ to prove the existence of minimizers in $H^2_{\iso}(S;\R^3)\times H^1(S;\mathbb S^2)$.

\subsection{Boundary conditions for $y$}\label{S:BC}
We establish refined versions of Theorem \ref{T1}, where we take boundary conditions for $y$ and/or $n$ into account. In this section we focus on boundary conditions for $y$.
As it is well-known, cf.~\cite{FJM02}, not every type of boundary condition for the deformation is compatible with the bending scaling.
We therefore restrict our analysis to the natural situation of clamped, affine boundary conditions on ``straight'' parts of the lateral boundary $\partial S\times(-\frac12,\frac12)$ of $\Omega$.
For the precise statement we make the following assumption on the boundary data:
\begin{assumption}[Boundary conditions for the 2d-model]\label{ass:BC}
Consider $k_{BC}\in\N$ relatively open, non-empty line segments $\mathcal L_i\subseteq\partial S$ and a reference isometry $y_{BC}\in H^2_{\iso}(S;\R^3)$ such that for $i=1,\ldots,k_{BC}$  (the trace of) $\nabla' y_{BC}$ is constant on $\mathcal L_{i}$.
We introduce a space of bending deformations that satisfy the following boundary conditions on the line segments:
  \begin{equation*}
    \mathcal A_{BC}\colonequals \Big\{y\in H^2_{\iso}(S;\R^3)\,:\,y=y_{BC}\text{ and }\nabla' y=\nabla'y_{BC}\text{ on }\mathcal L_1\cup\ldots \cup \mathcal L_{k_{BC}}\}.
  \end{equation*}
\end{assumption}
As we shall see in Theorem \ref{T2}, the 2d-boundary conditions emerge from sequences of 3d-deformations $(y_h)\subseteq H^1(\Omega;\R^3)$ with finite bending energy and subject to the following clamped, affine 3d-boundary conditions,
\begin{equation}\label{eq:BC3d}
  \begin{aligned}
    &y_h(x',x_3)=(1-h \delta)y_{BC}(x')+hx_3b_{y_{BC}}(x')
  \qquad\text{on }(\mathcal L_1\cup\ldots\cup \mathcal L_{k_{BC}})\times(-\tfrac12,\tfrac12),
\end{aligned}
\end{equation}
in the sense of traces and for a compression factor $\delta\in\R$.
In the following we write
\begin{equation}\label{eq:BC3dA}
  \mathcal A_{BC,\delta,h}\colonequals \big\{y\in H^1(\Omega;\R^3)\,:\,y\text{ satisfies }\eqref{eq:BC3d}\,\big\}.
\end{equation}
\begin{theorem}[Clamped, affine boundary conditions for $y$]\label{T2}
  Let Assumptions~\ref{ass:W}, \ref{ass:BC} and \eqref{ass:domain} be satisfied.
  \begin{enumerate}[(a)]
  \item (Compactness in $\mathcal A_{BC,\delta,h}$).
  	\label{item:T2:Comp}
  	Let $\delta\in\mathbb R$.
  	Consider the situation of Theorem~\ref{T1} (a) and additionally suppose that $y_h\in \mathcal A_{BC,\delta,h}$ for all $h$.
  	Then there exist $y\in\mathcal A_{BC}$ and $n\in H^1(S;\mathbb S^2)$ such that (up to a subsequence) we have \eqref{KompaktheitKonvergenzDef}, $y_h\to y$ strongly in $L^2(\Omega;\mathbb R^3)$, \eqref{KompaktheitKonvergenzDir} and \eqref{KompaktheitKonvergenzDir2}.
  \item (Recovery sequence in $\mathcal A_{BC,\bar\delta,h}$).
  	\label{item:T2:Recov}
  	Consider the situation of Theorem~\ref{T1} (c) and additionally suppose that $y\in \mathcal A_{BC}$.
  	Then for any
        \begin{equation}
          \bar\delta\geq\bar{r}\,C_W\label{eq:bardelta},
        \end{equation}
        the statements of Theorem~\ref{T1} (c) hold for a sequence $(y_h)$ that additionally satisfies $y_h\in\mathcal A_{BC,\bar\delta,h}$.
  \end{enumerate}
\end{theorem}
(See Sections \ref{Sec:Proof:Compactness} and \ref{Sec:Proof:Upper:Bound} for the proof.)\medskip

\begin{remark}[Condition \eqref{eq:bardelta}]
  The construction of the recovery sequence on affine regions of  $y\in\mathcal A_{BC}$ is based on a wrinkling ansatz that allows for recovering order-$h$-stretches in in-plane directions.
  Compressions cannot be realized in this way.
  However, the nematic-elastic coupling tensor (c.f.~\eqref{eq:steplength}) requires to realize compressions of magnitude of the order $\bar r C_W$.
  The additional term $-h\bar\delta y_{BC}$ in the 3d-boundary condition (cf.~\eqref{eq:BC3d} and \eqref{eq:bardelta}) allows us to introduce a macroscopic in-plane compression, that, in combination with the wrinkling ansatz, yields a recovery sequence realizing both, order-$h$-in-plane stretches and the required compressions.
\end{remark}
A crucial ingredient in the construction of recovery sequences with prescribed boundary conditions is the following density argument, which is of independent interest.
It applies even to non-convex domains $S$:
\begin{proposition}[Density of smooth isometries in $\mathcal A_{BC}$]\label{P:approx}
  Let $S$ satisfy \eqref{ConditionStar} and let As\-sump\-tion~\ref{ass:BC} be satisfied.
  Then
  \begin{equation*}
    \mathcal A_{BC}\cap C^\infty(\overline S;\R^3)\text{ is dense in }(\mathcal A_{BC},\|\cdot\|_{H^2(S;\R^3)}).
  \end{equation*}
\end{proposition}
(See Section \ref{Sec:Proof:P:approx} for the proof.)\medskip

\begin{remark}[Extension to non-convex $S$]\label{R:nonconvex}
  In this paper, convexity of $S$ is only used for the construction of the recovery sequence (cf.~Theorem~\ref{T1} \ref{item:T1:upper:bound} and Theorem~\ref{T2} \ref{item:T2:Recov}).
With help of \cite{HornungApprox}, the results can be extended to non-convex domains satisfying the regularity condition \eqref{ConditionStar}.
  In fact, we prove Proposition~\ref{P:approx} already for non-convex $S$ by combining results from \cite{HornungApprox} with an extension argument.
However, for the construction of the recovery sequence we assume convexity of $S$ to keep the arguments less technical and self-contained:
Our construction uses a decomposition of the non-flat part $\{\II_y\neq 0\}$ (for smooth $y\in H^2_{\iso}(S;\R^3)$) into connected components.
Convexity of $S$ simplifies the geometry of such connected components, cf.~proof of Lemma~\ref{LemmaFinePropertiesOfSmoothisometries}.
\end{remark}

\subsection{Weak and strong anchoring for the director on bulk and boundary}\label{S:anchoring}
In this section we discuss a variety of anchorings for 3d-director fields $n_h:\OmegaT\rightarrow\mathbb S^2$ and their relation to anchorings for 2d-director fields $n:S\rightarrow\mathbb S^2$.
We shall distinguish \textit{weak} and \textit{strong anchoring}.
In the physics literature, strong anchoring refers to a point-wise condition for $n_h$.
Here, we impose them on two-dimensional surfaces of the body, in particular, on parts of the  lateral boundary, i.e., on
\begin{equation*}
\Gamma_{n}\times(0,\tfrac12)\quad\text{for some }\Gamma_n\subseteq\partial S.
\end{equation*}
On the other hand, weak anchoring is described by an additional term in the energy functional, which assigns positive, but finite energy to deviations of the director from a prescribed behavior.
We shall consider both surface and volumetric anchorings, where the anchoring is imposed on two-dimensional and three-dimensional subsets of the body, respectively.

From the modeling point of view, it is natural to phrase weak and strong anchorings for the director field in \textit{local coordinates}, i.e., in the coordinate system spanned by $\nabla\!_hy_h(x')$ where $x'\in\Omega$ denotes a material point and $y_h:\Omega\rightarrow\mathbb R^3$ the deformation.
Therefore, given a configuration $(y_h,n_h)\in H^1(\Omega;\R^3)\times H^1(\OmegaT;\mathbb S^2)$ we consider the transformed field
\begin{equation*}
  \hat n_h(x)\colonequals \frac{\operatorname{Cof}(\nabla\!_h y_h(x))^{-1} n_h(x)}{|\operatorname{Cof}(\nabla\!_h y_h(x))^{-1} n_h(x)|}=\frac{\nabla\!_h y_h(x)^\top n_h(x)}{|\nabla\!_h y_h(x)^\top n_h(x)|}.
\end{equation*}
We note that it is not possible to introduce surface anchorings for $\hat n_h$ directly, since, in general, the trace of $\hat n_h$ on two-dimensional surfaces is not well-defined. To overcome this obstacle we consider two types of workaround solutions:
\begin{itemize}
\item One possibility is to consider only parts of the boundary where clamped, affine boundary conditions for the deformation are assumed.
Indeed, if $\Gamma_{n}$ is contained in a flat part of $\partial S\times(0,\frac12)$, where $(y_h)$ satisfies \eqref{eq:BC3d} and has finite bending energy, then we have $\nabla\!_hy_h=R_\mathrm{BC}+O(h)$ close to $\Gamma_{n}$ in an $L^2$-sense for some $R_{BC}\in\SO 3$.
Thus, the ill-defined restriction of $\hat n_h$ to $\Gamma_{n}\times(0,\frac12)$ might be replaced by the trace of the field $R_\mathrm{BC}^\top n_h$, and we obtain a strong anchoring condition, which, for instance, takes the form
  \begin{equation*}
    R_\mathrm{BC}^\top n_h=\hat n_{BC}\qquad\text{on $\Gamma_{n}\times(0,\frac12)$ in the sense of traces},
  \end{equation*}
  for some $\hat n_{BC}\in H^1(S;\mathbb S^2)$.
  As we shall see, in the  limit $h\to 0$ the strong 2d-anchoring condition $\hat n=\hat n_{BC}$ emerges, where $\hat n=R_y^\top n\in H^1(S;\mathbb S^2)$
 denotes the director in local coordinates w.r.t.~the isometric immersion $y$.
\item Another possibility is to consider volumetric weak anchoring through penalization energies of the form
  \begin{equation*}
    \int_{\OmegaT}|\hat n_h-\hat n_{BC}|^2\rho\,\mathrm{d} x,
  \end{equation*}
  where $\rho\in L^1(\OmegaT)$ denotes a non-negative weight. In the limit $h\to 0$ we recover the weak anchoring $\int_{S}|\hat n-\hat n_{BC}|^2\bar\rho(x')\,\mathrm{d} x'$ where $\bar\rho(x')\colonequals \int_0^\frac12\rho(x',x_3)\,\mathrm{d} x_3$.
\end{itemize}
As we shall see, with help of volumetric weak anchoring, we can also justify
\begin{itemize}
\item weak and strong 2d-anchorings on parts of $\partial S$ where no conditions on the deformation are imposed (see Lemma \ref{L:weak:anchoring:with:boundary:slices}),
\item weak and strong 2d-anchorings on the mid-surface of the plate (see Lemma \ref{L:volweakanch}).
\end{itemize}
Finally, we remark that we do not only consider anchorings of Dirichlet type (i.e., where the value $\hat n$ is fully prescribed), but also conditions or penalizations that enforce tangentiality of the director in the 2d-model, i.e.,~$n\cdot b_y=0$, where $y$ is the 2d-deformation.

\paragraph{Strong anchoring on clamped boundaries.}
We first state a result that is a rather direct consequence of Theorem~\ref{T1} and treats general linear strong anchorings phrased for $n$ in global coordinates:
\begin{lemma}[Clamped, affine boundary conditions and strong anchoring]
\label{P:Gamma:Strong:Anchoring}
Let Assumption~\ref{ass:W} and \eqref{ass:domain} be satisfied, and let $\Gamma_{n}\subseteq\partial S$ be relatively open.
Consider $\Psi_{BC}:\Gamma_n\times\R^3\to\R^3$, $\Psi_{BC}(x',n)\colonequals m(x')+M(x')n$ with $m\in L^2(\Gamma_n;\mathbb R^3)$ and $M\in L^2(\Gamma_n;\R^{3\times 3})$.
  \begin{enumerate}[(a)]
  \item (Compactness). Consider a sequence $(y_h,n_h)$ with equibounded energy in the sense of \eqref{EqFiniteEnergy} and assume that
    \begin{equation}\label{BC:SD3}
      \Psi_{BC}(x',n_h(x))=0\qquad\text{for a.e.~}x=(x',x_3)\in\Gamma_n\times(-\frac12,\frac12).
    \end{equation}
    Then in addition to the statements of Theorem~\ref{T1} \ref{item:T1:compactness}, any limit $(y,n)$ of a subsequence of $(y_h,n_h)$ satisfies
    \begin{equation}\label{BC:SD2}
      \Psi_{BC}(x',n(x'))=0\qquad\text{for a.e.~}x'\in\Gamma_n.
    \end{equation}
  \item (Recovery sequence).
  Let $(y,n)\in\mathcal A_2$ satisfy \eqref{BC:SD2}.
  Then there exists a recovery sequence $(y_h,n_h)$ that satisfies the properties of Theorem~\ref{T1} \ref{item:T1:upper:bound}, and additionally \eqref{BC:SD3}.
  Furthermore, if Assumption~\ref{ass:BC} is satisfied and $y\in\mathcal A_{BC}$, then the recovery sequence $(y_{h}, n_{h})$ can be constructed such that it additionally satisfies $y_h\in\mathcal A_{BC,\bar\delta,h}$ with $\bar\delta$ satisfying \eqref{eq:bardelta}.
  \end{enumerate}
\end{lemma}
(See Section \ref{Sec:Proof:213:214} for the proof.)\medskip

With help of the previous lemma we may treat strong Dirichlet and tangential anchorings:
For an illustration, we denote by $\mathcal L\subseteq\partial S$ one of the line segments of Assumption~\ref{ass:BC} and consider some relatively open $\Gamma_{n}\subseteq\mathcal L$.
Note that for $y\in\mathcal A_{BC}$, we have $R_y=R_{BC}$ on $\mathcal L$ for a single rotation $R_{BC}\in\SO{3}$.
For $3d$-director fields $n_h\in H^1(\OmegaT;\mathbb S^2)$ and 2d-director fields $n\in H^1(S;\mathbb S^2)$ we may consider the following strong anchorings:
\begin{itemize}
\item \textit{(Strong Dirichlet anchoring).}
Let $\hat n_{BC}:\Gamma_n\to\mathbb S^2$ be given (and admissible in the sense that it is the trace of some director field in $H^1(S;\mathbb S^2)$).
Consider
\begin{equation}\label{eq:Psi:BC:Dirichlet}
\Psi_{BC}(x',n)\colonequals R^\top_{{BC}}(x')n-\hat n_{BC}(x').
\end{equation}
Then Lemma~\ref{P:Gamma:Strong:Anchoring} shows that the strong 3d-anchoring condition $R_{{BC}}^\top n_h\,=\,\hat n_{BC}\text{ on }\Gamma_{n}\times(0,\frac12)$ leads to the 2d-anchoring condition $R_{{BC}}^\top n\,=\,\hat n_{BC}\text{ on }\Gamma_{n}$.
\item \textit{(Strong tangential anchoring).}
Let $\nu_S$ denote the (constant) outer unit normal to $S$ on $\Gamma_n$.
Consider
\begin{equation}\label{eq:Psi:BC:Tan}
\Psi_{BC}(x',n)\colonequals \big(R^\top_{{BC}}(x')n\cdot(\nu_S,0)\big)e_1.
\end{equation}
Then Lemma~\ref{P:Gamma:Strong:Anchoring} shows  that the strong 3d-anchoring condition $R_{{BC}}^\top n_h\perp (\nu_S,0)$ on $\Gamma_{n}\times(0,\frac12)$ leads to the 2d-anchoring condition $R_{{BC}}^\top n\perp(\nu_S,0)$ on $\Gamma_{n}$.
\end{itemize}

\paragraph{Volumetric  anchoring.}
In this paragraph, we derive 2d-weak and strong anchorings imposed on the mid-surface starting with weak 3d-volumetric anchorings in the form of weighted integrals.
To that end let $\rho\in L^1(\OmegaT)$ be a non-negative weight and set $\bar \rho(x')\colonequals \int_0^\frac12\rho(x',x_3)\,\mathrm{d} x_3$.
Furthermore, we consider $\hat n_0,\hat\nu\in L^2(S;\mathbb S^2)$ for describing preferred orientations of the director on the mid-surface. We introduce functionals $\mathcal G_h,\mathcal H_h:L^2(\Omega;\R^3)\times L^2(\OmegaT;\R^3)\to[0,\infty]$ defined for $(y_h,n_h)\in H^1(\Omega;\mathbb R^3)\times H^1(\OmegaT;\mathbb S^2)$ with $\det\nabla\!_hy_h>0$ a.e.~in $\Omega$ by
  \begin{equation*}
    \mathcal G_h(y_h,n_h)\colonequals \int_{\OmegaT}\left|\frac{\nabla\!_h y^\top_h n_h}{|\nabla\!_h y^\top_h n_h|}-\hat n_0\right|^2\rho\,\mathrm{d} x,\qquad\mathcal H_h(y_h,n_h)\colonequals \int_{\OmegaT}\left|\frac{\nabla\!_h y^\top_h n_h}{|\nabla\!_h y^\top_h n_h|}\cdot\hat\nu\right|^2\rho\,\mathrm{d} x,
  \end{equation*}
  and $\mathcal G_h(y_h,n_h)\colonequals \mathcal H_h(y_h,n_h)\colonequals +\infty$ otherwise. To express the corresponding 2d-anchorings, we introduce the functionals
  $\mathcal G_{\mathrm{weak/strong}},\mathcal H_{\mathrm{weak/strong}}:L^2(S;\R^3)\times L^2(S;\R^3)\to[0,\infty]$ defined by
  \begin{alignat*}{3}
    \mathcal G_{\mathrm{weak}}(y,n)&\colonequals
    \begin{cases}
      \int_{S}|R^\top_y n-\hat n_0|^2\bar\rho\,\mathrm{d} x'&\text{if }(y,n)\in\mathcal A_2,\\
      +\infty&\text{else},
    \end{cases}\\
    \mathcal G_{\mathrm{strong}}(y,n)&\colonequals
    \begin{cases}
      0&\text{if }(y,n)\in\mathcal A_2\text{ and }(R^\top_y n-\hat n_0)\bar\rho=0\text{ a.e.~in }S,\\
      +\infty&\text{else,}
    \end{cases}\\
    \mathcal H_{\mathrm{weak}}(y,n)&\colonequals
    \begin{cases}
      \int_{S}|(R^\top_y n)\cdot\hat\nu|^2\bar\rho\,\mathrm{d} x'&\text{if }(y,n)\in\mathcal A_2,\\
      +\infty&\text{else},
    \end{cases}\\
    \mathcal H_{\mathrm{strong}}(y,n)&\colonequals
    \begin{cases}
      0&\text{if }(y,n)\in\mathcal A_2\text{ and }(R^\top_y n\cdot\hat\nu)\bar\rho=0\text{ a.e.~in }S,\\
      +\infty&\text{else.}
    \end{cases}
  \end{alignat*}

  \begin{lemma}[Derivation of anchorings on the mid-surface]\label{L:volweakanch}\mbox{}
    \begin{enumerate}[(a)]
    \item The statements of Theorem~\ref{T1} hold with $(\mathcal E_h,\mathcal E)$ being replaced by $(\mathcal E_h+\mathcal G_h,\mathcal E+\mathcal G_\mathrm{ weak})$.
    They also hold with $(\mathcal E_h,\mathcal E)$ being replaced by $(\mathcal E_h+\mathcal H_h,\mathcal E+\mathcal H_{\mathrm weak})$.
    \item
    	Let $0<\beta<1$.
    	Then the  statements of Theorem~\ref{T1} hold with $(\mathcal E_h,\mathcal E)$ being replaced by $(\mathcal E_h+h^{-\beta}\mathcal G_h,\mathcal E+\mathcal G_{\mathrm{strong}})$.
    They also hold with $(\mathcal E_h,\mathcal E)$ being replaced by $(\mathcal E_h+h^{-\beta}\mathcal H_h,\mathcal E+\mathcal H_{\mathrm{strong}})$.
    \end{enumerate}
\end{lemma}
(See Section \ref{Sec:Proof:213:214} for the proof.)\medskip

In the special case $\rho\colonequals \frac{2}{\sigma}$ and $\hat\nu\colonequals e_3$, we obtain by Lemma~\ref{L:volweakanch} (a) the 2d-energy
\begin{equation*}
  (y,n)\quad\mapsto\quad\mathcal E(y,n)+\frac1\sigma\int_{S}|n\cdot b_y|^2\,\mathrm{d} x',
\end{equation*}
which features a weak anchoring term that penalizes non-tangentiality of $n$.
Likewise Lemma~\ref{L:volweakanch} (b) leads to the point-wise, tangentiallity constraint
\begin{equation*} 
  n\cdot b_y=0\qquad\text{a.e.~in }S.
\end{equation*}
\paragraph{Anchoring on lateral boundaries $\Gamma_{n}\subseteq\partial S$.}
Next, we justify weak and strong 2d-anchorings imposed on boundary segments $\Gamma_{n}\subseteq\partial S$. In contrast to Lemma \ref{P:Gamma:Strong:Anchoring}, here, anchorings for the director in \emph{local} coordinates are treated.
To simplify the analysis, we make the following assumption:
\begin{equation}
  \label{eq:gamma_n}
\begin{aligned}
  &\text{
  	$\Gamma_n\subseteq \partial S$ is relatively open with finitely many connected components, and
  }\\
  &\text{
   the closure $\overline{\Gamma_n}$ is contained in a one-dimensional $C^2$-submanifold of $\R^2$.
  }
\end{aligned}
\end{equation}
As mentioned before, it is not possible to directly consider anchorings for the director in local coordinates on a 2d-surface $\Gamma_{n}\times(0,\frac12)$. We therefore introduce the set
\begin{equation*}
  \Gamma_{n}^\sigma\colonequals \{\xi'-s\,\nu_{S}(\xi')\,:\,\xi'\in\Gamma_{n},\,0<s<\sigma\,\},
\end{equation*}
where $\nu_S$ is the outer unit normal to $S$.
Note that $\Gamma_n^\sigma$ is a thin ``boundary strip'' that intersects $S$, and that for $\sigma\to 0$ ``approximates'' the boundary segment $\Gamma_{n}$.
We shall consider weak volumetric anchoring energies on $(S\cap\Gamma_n^\sigma)\times(0,\frac12)$.
To formulate these, we need to extend the given boundary director field $\hat n_{BC}:\Gamma_n\to\mathbb S^2$ to the strip $\Gamma_n^\sigma$.
We consider the following canonical extension that extends to a field that is constant along fibers orthogonal to $\Gamma_n$:
\begin{lemma}[Canonical extension]\label{L:ext}
  Let $S\subseteq \mathbb R^2$ be a bounded Lipschitz domain.
  Assume \eqref{eq:gamma_n}.
  Then there exists $\bar\sigma>0$ (only depending on $\Gamma_n$ and $S$) such that for all $g\in L^2(\Gamma_n;\mathbb S^2)$, the identity
  \begin{equation*}
    (\mathrm{E}g)(x'-s\nu_S(x'))\colonequals g(x')\qquad\text{for a.e.~}x'\in\Gamma_n\text{ and }s\in(0,\bar\sigma)
  \end{equation*}
  defines a unique function $\mathrm{E}g\in L^2(\Gamma_n^{\bar\sigma};\mathbb S^2)$. We call $\mathrm{E}g$ the \emph{canonical extension} of $g$.
\end{lemma}
(See Section \ref{Sec:Proof:215:217} for the proof.)\medskip

We are now in the position to introduce the functionals that model the volumetric weak anchorings:
For given data $\hat n_{BC}$, $\hat\nu\in L^2(\Gamma_n;\mathbb S^2)$ and $\sigma>0$  we consider $\mathcal G_h^\sigma,\mathcal H_h^\sigma:L^2(\Omega;\R^3)\times L^2(\OmegaT;\R^3)\to[0,\infty]$ defined for $(y_h,n_h)\in H^1(\Omega;\mathbb R^3)\times H^1(\OmegaT;\mathbb S^2)$ with $\det \nabla\!_hy_h>0$ a.e.~in $\Omega$ by
  \begin{align*}
    \mathcal G_h^\sigma(y_h,n_h)&\colonequals \frac2\sigma\int_{(S\cap\Gamma_n^\sigma)\times(0,\frac12)}\left|\frac{\nabla\!_h y^\top_h n_h}{|\nabla\!_h y^\top_h n_h|}-\mathrm{E}\hat n_{BC}\right|^2\,\mathrm{d} x,\\
    \mathcal H_h^\sigma(y_h,n_h)&\colonequals \frac2\sigma\int_{(S\cap\Gamma_n^\sigma)\times(0,\frac12)}\left|\frac{\nabla\!_h y^\top_h n_h}{|\nabla\!_h y^\top_h n_h|}\cdot\mathrm{E}\hat\nu\right|^2\,\mathrm{d} x,
  \end{align*}
  and  $\mathcal G_h^\sigma(y_h,n_h)\colonequals \mathcal H_h^\sigma(y_h,n_h)\colonequals +\infty$ otherwise. Above,   $\mathrm{E}\hat n_{BC}$ and $\mathrm{E}\hat\nu$ denote the canonical extension of $\hat n_{BC}$ and $\hat\nu$.

  By appealing to Lemma~\ref{L:volweakanch} (a), for $\sigma>0$ fixed, we shall obtain weak 2d-anchoring energies on the boundary strip $(S\cap\Gamma_n^\sigma)$.
   The limiting anchorings are given by the functionals $\mathcal G^\sigma,\mathcal H^\sigma\,:\,L^2(S;\R^3)\times L^2(S;\R^3)\rightarrow [0,\infty]$ defined for $(y,n)\in \mathcal A_2$ by
  \begin{align*}
    \mathcal G^\sigma(y,n)\colonequals \frac1\sigma\int_{(S\cap\Gamma_n^\sigma)}|R^\top_y n-\mathrm{E}\hat n_{BC}|^2\,\mathrm{d} x',\qquad
    \mathcal H^\sigma(y,n)\colonequals \frac1\sigma\int_{(S\cap\Gamma_n^\sigma)}|R^\top_y n\cdot\mathrm{E}\hat\nu|^2\,\mathrm{d} x',
  \end{align*}
  and $+\infty$ otherwise.
  \begin{corollary}\label{L:preliminary:weak:anchoring:with:boundary:slices}
    Assume \eqref{eq:gamma_n}.
    Then for all scaling parameters $\beta\geq 0$ and $\sigma\in(0,\bar\sigma)$ (where $\bar\sigma>0$ only depends on $\Gamma_n$ and $S$), the statements of Theorem~\ref{T1} hold with $(\mathcal E_h, \mathcal E)$ being replaced by $(\mathcal E_h+\sigma^{-\beta}\mathcal G_h^\sigma, \mathcal E+\sigma^{-\beta}\mathcal G^\sigma)$ and also with $(\mathcal E_h, \mathcal E)$ being replaced by $(\mathcal E_h+\sigma^{-\beta}\mathcal H_h^\sigma, \mathcal E+\sigma^{-\beta}\mathcal H^\sigma)$.
  \end{corollary}
  The proof directly follows from Lemma~\ref{L:ext} and Lemma~\ref{L:volweakanch} (a) and is thus left to the reader.

  In a second step, we argue that a subsequent $\Gamma$-limit for $\sigma\to 0$ yields weak and strong anchorings on $\Gamma_{n}$.
  To describe the latter, we introduce the functionals
  $\mathcal G_{\mathrm{weak/strong},\Gamma_n},\mathcal H_{\mathrm{weak/strong},\Gamma_n}:L^2(S;\R^3)\times L^2(S;\R^3)\to[0,\infty]$,
  \begin{alignat*}{3}
    \mathcal G_{\mathrm{weak},\Gamma_n}(y,n)&\colonequals
    \begin{cases}
      \int_{\Gamma_n}|R^\top_y n-\hat n_{BC}|^2\,\mathrm d\mathcal H^1&\text{if }(y,n)\in\mathcal A_2,\\
      +\infty&\text{else},
    \end{cases}\\
    \mathcal G_{\mathrm{strong},\Gamma_n}(y,n)&\colonequals
    \begin{cases}
      0&\text{if }(y,n)\in\mathcal A_2\text{ and }R^\top_y n=\hat n_{BC}\text{ a.e.~on }\Gamma_n,\\
      +\infty&\text{else,}
    \end{cases}\\
    \mathcal H_{\mathrm{weak},\Gamma_n}(y,n)&\colonequals
    \begin{cases}
      \int_{S}|(R^\top_y n)\cdot\hat\nu|^2\,\mathrm{d}\mathcal H^1&\text{if }(y,n)\in\mathcal A_2,\\
      +\infty&\text{else},
    \end{cases}\\
    \mathcal H_{\mathrm{strong},\Gamma_n}(y,n)&\colonequals
    \begin{cases}
      0&\text{if }(y,n)\in\mathcal A_2\text{ and }R^\top_y n\cdot\hat\nu=0\text{ a.e.~on }\Gamma_n,\\
      +\infty&\text{else.}
    \end{cases}
  \end{alignat*}
\begin{lemma}
 \label{L:weak:anchoring:with:boundary:slices}
 Let $S\subseteq \mathbb R^2$ be a bounded Lipschitz domain.
  Assume \eqref{eq:gamma_n}.
  Then the following holds:
  \begin{enumerate}[(a)]
  \item Let $\beta=0$. As $\sigma\to 0$ the functionals $\mathcal E+\mathcal G^\sigma$ and $\mathcal E+\mathcal H^\sigma$ $\Gamma$-converge (w.r.t.~strong convergence in $L^2(S)$) to the functionals $\mathcal E+\mathcal G_{\mathrm{weak},\Gamma_n}$ and $\mathcal E+\mathcal H_{\mathrm{weak},\Gamma_n}$, respectively.
  \item
  Let $0<\beta<\frac12$.
  As $\sigma\to 0$ the functionals $\mathcal E+\sigma^{-\beta}\mathcal G^\sigma$ and $\mathcal E+\sigma^{-\beta}\mathcal H^\sigma$ $\Gamma$-converge (w.r.t.~strong convergence in $L^2(S)$) to the functionals $\mathcal E+\mathcal G_{\mathrm{strong},\Gamma_n}$ and $\mathcal E+\mathcal H_{\mathrm{strong},\Gamma_n}$, respectively.
  \end{enumerate}
\end{lemma}
(See Section \ref{Sec:Proof:215:217} for the proof.)\medskip

\begin{remark}[Relation to simulations of Section \ref{Sec:Freiburg}]
  \label{R:Bridge:Analysis:Numerics}

  The results of this and the previous section explain how 2d-boundary conditions for the deformation and the director field emerge from 3d-models.
  Let us anticipate that in Section \ref{Sec:Freiburg} we present numerical experiments for 2d-models with different types of boundary conditions.
  In this remark we comment on their derivation from 3d-models.
  \begin{enumerate}[(a)]
  \item In Example~\ref{ex:compress} a 2d-model is considered with $S\colonequals (-5,5)\times (-1,1)$, subject to boundary conditions on $\mathcal L_1=\{-5\} \times [-1,1]$ and $\mathcal L_2=\{5\} \times [-1,1]$ for the 2d-deformation, its gradient and the director.
  Furthermore, tangentiality of $n$ is enforced.
    The situation can be described by the 2d-energy $(\mathcal E+\mathcal H_\mathrm{strong})$, where $\mathcal H_{\mathrm{strong}}$ enforces the tangentiality of $n$ and is defined with $\bar\rho\colonequals 1$ and $\hat n\colonequals e_3$.
    The clamped, affine boundary conditions for $y$ on $\mathcal L_1\cup\mathcal L_2$ can be encoded by a suitable choice of $\mathcal A_{BC}$, while the strong anchoring conditions of $n$ take the form of the pointwise constraint \eqref{BC:SD2} with $\Psi_{BC}$ as in \eqref{eq:Psi:BC:Dirichlet}.
    This situation can be derived as a $\Gamma$-limit from a 3d-model of the form $(\mathcal E_h+h^{-\beta}\mathcal H_h)(y_h,n_h)$ (c.f.~Lemma \ref{L:volweakanch} (b)) subject to the clamped, affine boundary condition $y_h\in\mathcal A_{BC,\bar\delta,h}$ and the pointwise constraint $\Psi_{BC}(x',n_h(x))=0$ on $(\mathcal L_1\cup\mathcal L_2)\times(-\tfrac12,\tfrac12)$ for the director field.
    The passage to 2d then follows by combining Theorem \ref{T2}, Lemma \ref{P:Gamma:Strong:Anchoring} and a variant of Lemma \ref{L:volweakanch}.
    Indeed, one can check that Lemma \ref{L:volweakanch} (b) extends to the case where additional boundary conditions for the deformation and strong anchorings for the director field are considered.

  \item
  In Example~\ref{ex:clampedplate}, the deformation and its gradient, as well as the director, are fixed on a common boundary segment.
This situation arises as a $\Gamma$-limit via directly combining Theorem \ref{T1} (b), Theorem \ref{T2} and Lemma  \ref{P:Gamma:Strong:Anchoring} with $\Psi_{BC}$ as in \eqref{eq:Psi:BC:Dirichlet}.

\item
In Example \ref{ex:circular} boundary conditions are only prescribed for the 2d-director (in local coordinates).
This setting is motivated by the $\Gamma$-limits of Corollary \ref{L:preliminary:weak:anchoring:with:boundary:slices} and Lemma \ref{L:weak:anchoring:with:boundary:slices}~(b).
  \end{enumerate}
\end{remark}

\subsection{Modeling and non-dimensionalization}\label{S:modeling}
In this section, we explain the derivation of the scaled and non-dimensionalized energy $\mathcal E_h$ defined in \eqref{eq:energy}. In the following,  variables with a physical dimension are underlined, e.g, we write $\underline x$ for the position vector with physical dimensions and $x$ for the position vector in non-dimensionalized form.
We consider a plate with reference domain $\underline\Omega\colonequals \underline S\times (-\underline h/2,\underline h/2)$, mid-surface $\underline S\subseteq\R^2$, thickness $\underline h$, top layer $\underline{\Omega}_{\text{top}}\colonequals \underline S\times (0,\underline h/2)$, bottom layer $\underline{\Omega}_{\text{bot}}\colonequals \underline S\times (-\underline h/2,0)$, and the energy functional
\begin{equation}\label{eq:energy:dimensions}
  \begin{aligned}
    &\underline{\mathcal E}(\underline y,\underline n)\colonequals \,\int_{\underline\Omega_{\text{bot}}} \underline W\big(\underline{x}_3,\nabla\, \underline y\big)\mathrm{d}\underline{x} +\int_{\underline\Omega_{\text{top}}}\underline W\left(\underline{x}_3,L(\underline n)^{-\frac12}\underline y\right)\mathrm{d}\underline{x}+{\boldsymbol{\kappa}}\int_{\underline\Omega_{\text{top}}}|\nabla\underline n(\nabla\underline y)^{-1}|^2\det(\nabla\underline y)\,\mathrm{d}\underline{x},
  \end{aligned}
\end{equation}
where  $\underline y:\underline\Omega\to\R^3$ denotes the deformation of the plate and  $\underline n:\underline \Omega_{\text{top}}\to\mathbb S^2$ the LC orientation in the top layer.
The first two integrals are the elastic energy stored in the deformed body with a material law described by the stored energy function $\underline W$. The third integral is the one-constant approximation of the Oseen-Frank energy (cf.~\cite{calderer2005model, lin2001static, B_Warner07}), which penalizes spatial variations of $\underline n$.
It invokes as a material parameter the Frank elasticity constant $\boldsymbol{\kappa}$ and is originally formulated as $\boldsymbol{\kappa}\int_{\underline y(\underline\Omega_{\text{top}})}|\nabla\tilde{\underline{n}}|^2\,\mathrm{d}\underline{x}$ for the director field $\tilde{\underline n}\colonequals\underline n\circ \underline y^{-1}$ defined on the deformed configuration $\underline y(\underline\Omega_{\text{top}})$.
The third integral in \eqref{eq:energy:dimensions} is then obtained by (formally) pulling back the director field to the reference configuration.

The second integral in \eqref{eq:energy:dimensions} invokes the step-length tensor $L(\underline n)\colonequals r^{-\frac13}(I+(r-1)\underline n\otimes \underline n)$, which has been introduced in \cite{bladon1993transitions, B_Warner07} to model the nematic-elastic coupling of LCEs.
Here,  $r>0$ is an order-parameter.
The product $L(\underline n)^{-\frac12}\nabla \underline y$ can be understood in view of the multiplicative decomposition $\nabla \underline y = \sqrt{L(\underline n)}F_\mathrm{elast}$ of the deformation gradient into a ``poststrain'' $\sqrt{L(\underline n)}$ and an elastic part $F_\mathrm{elast}$. The elastic stress tensor only depends on the latter. We note that for $r>1$, $\sqrt{L(\underline n)}$ describes stretching by a factor $r^{\frac13}$ in directions parallel to $\underline n$ and a contraction by a factor $r^{-\frac16}$ in directions perpendicular to $\underline n$. The case $r= 1$ refers to the isotropic case, i.e., when the nematic-elastic coupling is absent.
\medskip

In order to non-dimensionalize the energy, we introduce a reference length scale $\boldsymbol{\ell}$ and denote by $\boldsymbol{p}$ a pressure scale.
We then introduce the non-dimensional and scaled quantities
\begin{gather*}
  S\colonequals \boldsymbol\ell^{-1}\underline S,\quad h\colonequals \boldsymbol\ell^{-1}\underline h,\quad x\colonequals \boldsymbol\ell^{-1} (\underline x',h^{-1}\underline x_3),\\
  y_h(x)\colonequals \boldsymbol\ell^{-1}\underline y(\boldsymbol\ell x',\boldsymbol\ell h x_3),\quad n_h(x)\colonequals \underline n(\boldsymbol\ell x',\boldsymbol\ell h x_3),\qquad W(x_3,F)\colonequals \boldsymbol{p}^{-1}\underline W(\boldsymbol\ell hx_3,F).
\end{gather*}
We observe that
\begin{equation*}
  \begin{aligned}
    &\frac{1}{\boldsymbol{p}\boldsymbol\ell^3h^3}\underline{\mathcal E}(\underline y,\underline n)=\,\mathcal E_h(y_h,n_h),
  \end{aligned}
\end{equation*}
where $\mathcal E_h$ (as defined in \eqref{eq:energy}) invokes the non-dimensional and scaled parameters
\begin{equation*}
  \bar\e^2\colonequals \frac{{\boldsymbol{\kappa}}}{\boldsymbol{p}\boldsymbol\ell^2h^2}\quad\text{ and }\quad\bar r\colonequals \frac{r-1}{h}.
\end{equation*}

For homogeneous and isotropic materials, a typical choice for $\boldsymbol p$ is the shear modulus $\boldsymbol \mu$ of $\underline W$. For this choice,  the quadratic term $Q$ in the expansion of $W$ at $I_{3\times 3}$ (cf.~\ref{item:ass:W:quadratic:exp}) satisfies \eqref{eq:Q:hom:iso} with parameters $\mu =1$ and $\lambda = \frac{\boldsymbol\lambda}{\boldsymbol \mu}$, where $\boldsymbol \mu$ and $\boldsymbol \lambda$ are the physical Lam\'e-parameters.
Furthermore, $\bar\e$ is precisely the ratio between the length scale $\sqrt{\frac{\boldsymbol{\kappa}}{\boldsymbol{p}}}$ and the thickness of the plate $\underline h=h\boldsymbol l$.
With the choice $\boldsymbol p:=\boldsymbol\mu$, the length scale $\sqrt{\frac{\boldsymbol{\kappa}}{\boldsymbol{p}}}$ is typically in the range $10^{-9}$--$10^{-7}\mathrm{m}$ (e.g.,~see \cite{verwey1996elastic,B_Warner07}), and thus $\bar\e$ is typically rather small -- yet values, of the order $10^{-3}$--$10^{-1}$ might be realistic, e.g.,~see \cite{haseloh2011nanosized} for composites with nano-sized nematic LCE grains, and see \cite{C9MH00951E} for elastomers with $\boldsymbol\mu$ being in the range of $10^3$--$10^5\mathrm{Pa}$.
Let us anticipate that in some of the numerical simulations in Section~\ref{Sec:Freiburg} we set $\bar\e=1$ (or even larger) to pronounce the effect of the coupling with the surface Oseen-Frank energy. In particular, with $\bar r=5$, $\bar\varepsilon =1$, $S=(-1,1)^2$, $\sqrt{\frac{\boldsymbol{\kappa}}{\boldsymbol{p}}}\sim 10^{-7}\mathrm{m}$, and the choice $\boldsymbol\ell = 10^{-4}\mathrm{m}$, we obtain $h=10^{-4}$ and $r=1+h\bar r=1.0005$.
In comparison, if we consider a plate with physical thickness $\underline h\sim 10^{-5}\mathrm{m}$ (as in \cite{SUTDT10}) and diameter $\boldsymbol\ell\sim 10^{-4}\mathrm{m}$, we arrive at $\overline\e=10^{-2}$, $h=10^{-1}$ and $r=1.5$.


\section{Simulation and two-dimensional model exploration} \label{Sec:Freiburg}
To simplify notation, we drop in this section the convention of denoting coordinates and operators in the two-dimensional setting with a prime and simply write $x$ and $\nabla$ instead of $x'$ and~$\nabla'$.
We present in the following a numerical scheme for approximating stationary points of the dimensionally reduced energy functional in the case of a homogeneous and isotropic material.
In this case the elastic and local residual energy contributions are given as in Lemma~\ref{L:representations:hom:iso} and the resulting energy functional is
\begin{equation*}
\begin{aligned}
  \mathcal{E}(y,n)
  = &\frac{1}{12} \int_S
      \mu \left| \II_y + \frac{3\bar{r}}{4}\left(
        \frac{1}{3} \id -  \nabla y^\top n \otimes \nabla y^\top n
      \right)\right|^2
    \mathrm{d}x \\
    &+ \frac{1}{12} \int_S
      \bar{\lambda} \operatorname{tr}\left( \II_y + \frac{3\bar{r}}{4}\left(
        \frac{1}{3} \id -  \nabla y^\top n \otimes \nabla y^\top n
      \right)\right)^2
    \mathrm{d} x \\
    &+ \frac{1}{64} \int_S \mu \bar{r}^2 \left| \frac{1}{3} \id -  \nabla y^\top n \otimes \nabla y^\top n \right|^2
      + \bar{\lambda} \bar{r}^2 \operatorname{tr} \left( \frac{1}{3} \id -  \nabla y^\top n \otimes \nabla y^\top n \right)^2 \mathrm{d} x\\
    &+ \frac{\bar{\e}^2}{2} \int_S \left|\nabla n \right|^2 \mathrm{d}x,
\end{aligned}
\end{equation*}
where $\bar{\lambda} = \mu\lambda(2\mu + \lambda)^{-1}$.

\subsection{Finite element discretization}\label{subsec:discret}
We reformulate the energy functional before discretizing as follows.
For any $s \in \R^2$ it holds that $|s \otimes s|^2 = \operatorname{tr}(s \otimes s)^2$.
Furthermore, for isometries $y \in H^2_\mathrm{iso}(S;\R^3)$ we have that $\II_y = -(\partial_1 y \wedge \partial_2 y)^\top \nabla \nabla y$ as well as  $|\II_y| = |\operatorname{tr} \II_y| = |\nabla \nabla y|$. After expanding the quadratic terms and reordering, the energy can be written as
\begin{equation*}
\begin{aligned}
  \mathcal{E}(y,n) = &\frac{1}{12} (\mu + \bar{\lambda}) \int_S \left|\nabla \nabla y \right|^2 \mathrm{d}x
  + \frac{\bar{r}}{8} \int_S
    \mu \II_y\cdot\left( \frac{1}{3} \id -  \nabla y^\top n \otimes \nabla y^\top n \right)
  \mathrm{d}x \\
  & + \frac{\bar{r}}{8} \int_S
    \bar{\lambda} \left(
      \Delta y \cdot \left( \partial_1 y \wedge \partial_2 y \right)
    \right)
    \operatorname{tr} \left(
      \frac{1}{3} \id -  \nabla y^\top n \otimes \nabla y^\top n
    \right)
  \mathrm{d}x \\
  & + \frac{\bar{r}^2}{16} \int_S
    \left( \mu + \bar{\lambda} \right) \operatorname{tr} \left(
      \nabla y^\top n \otimes \nabla y^\top n
    \right)^2
    - \frac{2\mu + 4 \bar{\lambda}}{3} \operatorname{tr} \left(
      \nabla y^\top n \otimes \nabla y^\top n
    \right)
  \mathrm{d}x \\
  & + \frac{\bar{r}^2}{72} (\mu + 2\bar{\lambda}) |S|
  + \frac{\bar{\e}^2}{2} \int_S |\nabla n|^2 \mathrm{d}x.
\end{aligned}
\end{equation*}
We assume in the following that the two-dimensional domain $S$ is polygonal and divided into triangulations $\mathcal{T}_\mathfrak{h}$, $\mathfrak{h}>0$, consisting of triangles with maximum diameter less than or equal to $\mathfrak{h}$, such that the Dirichlet boundaries $\Gamma_y = \bigcup_{i=1}^{k_\mathrm{BC}} \mathcal{L}_i$ and $\Gamma_n$ are matched exactly by a union of element sides, cf.~Assumption~\ref{ass:BC} and Lemma~\ref{P:Gamma:Strong:Anchoring} with $\Psi_{BC}$ as in~\eqref{eq:Psi:BC:Dirichlet}.
We denote with $\mathcal{S}_\mathfrak{h}$ the set of all sides of elements and with $\mathcal{N}_\mathfrak{h}$ the set of all vertices in the triangulation $\mathcal{T}_\mathfrak{h}$.
For $k>0$ let $P_k(T)$ be the set of all polynomials of degree $k$ or less on $T\in\mathcal{T}_\mathfrak{h}$, and
\[
P_3^\mathrm{red}(T) \colonequals  \bigg\{ p \in P_3(T) : p(x_T) = \frac{1}{3}\sum_{z \in \mathcal{N}_\mathfrak{h} \cap T} \Big(p(z) + \frac{1}{2}\nabla p(z) \cdot (x_T - z)\Big)\bigg\}.
\]
the reduced set of cubic polynomials resulting from $P_3(T)$ after prescribing the function value at the center of mass $x_T = (1/3) \sum_{z \in \mathcal{N}_\mathfrak{h} \cap T} z$.
Our approximation employs the finite element spaces
\begin{align*}
  \mathcal{S}^1(\mathcal{T}_\mathfrak{h}) &\colonequals  \left\{
    n_\mathfrak{h} \in C(S) : n_\mathfrak{h}|_T \in P_1(T) \text{ for all } T \in \mathcal{T}_\mathfrak{h}
  \right\},\\
  \mathcal{S}^\mathrm{DKT}(\mathcal{T}_\mathfrak{h}) &\colonequals  \left\{
    y_\mathfrak{h} \in C(S) : y_\mathfrak{h}|_T \in P_3^\mathrm{red}(T) \text{ for all } T \in \mathcal{T}_\mathfrak{h}, \; \nabla y_\mathfrak{h} \text{ continuous in all } z \in \mathcal{N}_\mathfrak{h}
  \right\},\\
  \mathcal{S}^2(\mathcal{T}_\mathfrak{h}) &\colonequals  \left\{
    \theta_\mathfrak{h} \in C(S) : \theta_\mathfrak{h} \in P_2(T) \text{ for all } T \in \mathcal{T}_\mathfrak{h}
  \right\}.
\end{align*}
In the spatial discretization of the energy, we employ the space of piecewise linear finite elements $\mathcal{S}^1(\mathcal{T}_\mathfrak{h};\R^3) \subseteq H^{2}(S; \R^3)$ for the director field and nonconforming discrete Kirchhoff triangular (DKT) elements $\mathcal{S}^\mathrm{DKT}(\mathcal{T}_\mathfrak{h}; \R^3)$ for the deformation.
The space $\mathcal{S}^2(\mathcal{T}_\mathfrak{h})$ serves as an $H^{2}$-conforming approximation space for gradients of discrete deformations from $\mathcal{S}^\mathrm{DKT}(\mathcal{T}_\mathfrak{h}; \R^3)$ via the construction of a discrete gradient operator
\[
\Theta_\mathfrak{h} \colon \mathcal{S}^\mathrm{DKT}(\mathcal{T}_\mathfrak{h}) \to \mathcal{S}^2(\mathcal{T}_\mathfrak{h}; \R^2)
\]
as follows.
First, we note that the degrees of freedom in $\mathcal{S}^\mathrm{DKT}$ are the function values and gradients at the nodes of the triangulation whereas the degrees of freedom in $\mathcal{S}^2(\mathcal{T}_\mathfrak{h})$ are the function values at the nodes and midpoints of element sides.
Denoting with $z_E^1,z_E^2 \in \mathcal{N}_\mathfrak{h}$ the endpoints and with $z_E = (1/2)(z_E^1 + z_E^2)$ the midpoint of a given side $E \in \mathcal{S}_\mathfrak{h}$, and choosing normal and tangent vectors $\tau_E$ and $\nu_E$ for every side $E \in \mathcal{S}_\mathfrak{h}$ such that $|\tau_E| = |\nu_E| = 1$, we define the unique discrete gradient field $\Theta_\mathfrak{h} y_\mathfrak{h} \in \mathcal{S}^2(\mathcal{T}_\mathfrak{h}; \R^2)$ of a function $y_\mathfrak{h} \in \mathcal{S}^\mathrm{DKT}(\mathcal{T}_\mathfrak{h})$ via conditions
\begin{align*}
  ( \Theta_\mathfrak{h} y_\mathfrak{h} )(z) &= \nabla y_\mathfrak{h} (z), \\
  ( \Theta_\mathfrak{h} y_\mathfrak{h} )(z_E) \cdot \tau_E &= \nabla y_\mathfrak{h} (z_E) \cdot \tau_E,\\
  ( \Theta_\mathfrak{h} y_\mathfrak{h} )(z_E) \cdot \nu_E &= \frac{1}{2} \left(\nabla y_\mathfrak{h} (z_E^1)  + \nabla y_\mathfrak{h} (z_E^2)\right)\cdot \nu_E,
\end{align*}
for all $z \in \mathcal{N}_\mathfrak{h}$ and all $E \in \mathcal{S}_\mathfrak{h}$.
The definition of the discrete gradient is extended to vector fields $y_\mathfrak{h} \in \mathcal{S}^\mathrm{DKT}(\mathcal{T}_\mathfrak{h}; \R^3)$ by applying $\Theta_\mathfrak{h}$ to each component of $y_\mathfrak{h}$.
Approximation properties of $\Theta_\mathfrak{h}$ show that the discrete gradient defines an interpolation operator on the space of derivatives of functions in $H^3(S)$, see Chapter~8.2 of~\cite{bartels2015numerical}.
In the continuous setting the director field $n$ takes values in the sphere $\mathbb{S}^2$ and the deformation $y$ is an isometry in $H^2_\mathrm{iso}(S; \R^3)$. In the discrete setting we impose the unit-length constraint as well as the isometry constraint up to a tolerance $\zeta_\mathfrak{h} > 0$ in the nodes of the triangulation, i.\,e.~we consider the discrete admissible sets
\begin{align*}
  \mathcal{A}_{n,\mathfrak{h}} \colonequals  \left\{
    n_\mathfrak{h} \in \mathcal{S}^1(\mathcal{T}_\mathfrak{h}; \R^3) : n_\mathfrak{h}(z) = n_\mathrm{BC}(z) \text{ for all } z \in \mathcal{N}_\mathfrak{h} \cap \Gamma_n,
    \; \left|n_\mathfrak{h}(z) - 1 \right| \le \zeta_\mathfrak{h} \text{ for all } z \in \mathcal{N}_\mathfrak{h}
  \right\}
\end{align*}
and, assuming that $\Phi_\mathrm{BC} = \nabla y_\mathrm{BC} \in C(\Gamma_y; \R^{3\times2})$ with $y_\mathrm{BC}$ as in Assumption~\ref{ass:BC},
\begin{align*}
  \mathcal{A}_{y,\mathfrak{h}} \colonequals  \Big\{
    y_\mathfrak{h} \in \mathcal{S}^\mathrm{DKT}(\mathcal{T}_\mathfrak{h}; \R^3) : y_\mathfrak{h}(z) = y_\mathrm{BC}(z),\; \nabla y_\mathfrak{h}(z) = \Phi_\mathrm{BC}(z)
    \text{ for all } z \in \mathcal{N}_\mathfrak{h} \cap \Gamma_y,\\
    \left|\left(\nabla y_\mathfrak{h}(z)\right)^\top \nabla y_\mathfrak{h}(z) - \Id \right| \le \zeta_\mathfrak{h}  \text{ for all } z \in \mathcal{N}_\mathfrak{h}
  \Big\}.
\end{align*}
We use the shorthand notations $\Delta_\mathfrak{h} \colonequals  \operatorname{div} \Theta_\mathfrak{h}$ for the divergence of the discrete gradient and $b_\mathfrak{h}(z) \colonequals  \partial_1 y_\mathfrak{h}(z) \wedge \partial_2 y_\mathfrak{h}(z)$ for the nodal (almost-)unit normal and denote with $\mathcal{I}_\mathfrak{h}$ the interpolant operator into element-wise linear functions.
For $(y_\mathfrak{h},n_\mathfrak{h}) \in \mathcal{A}_{y,\mathfrak{h}} \times \mathcal{A}_{n,\mathfrak{h}}$ the energy is now discretized via
\begin{equation}\label{eq:energy-discrete1}
  \mathcal{E}_\mathfrak{h}(y_\mathfrak{h},n_\mathfrak{h}) \colonequals  \frac{1}{12} (\mu + \bar{\lambda}) \int_S
    \left |\nabla \Theta_\mathfrak{h} y_\mathfrak{h} \right|^2
  \mathrm{d}x
  + N_{\bar{r},\mathfrak{h}}(y_\mathfrak{h},n_\mathfrak{h})
  + \frac{\bar{\e}^2}{2} \int_S |\nabla n_\mathfrak{h}|^2 \mathrm{d}x,
\end{equation}
where $N_{\bar{r},\mathfrak{h}}$ is a quadrature operator for the non-convex terms given by
\begin{equation*}
\begin{aligned}
  N_{\bar{r},\mathfrak{h}}(y_\mathfrak{h},n_\mathfrak{h}) \colonequals  &\frac{\bar{r}}{8} \int_S \mathcal{I}_\mathfrak{h} \left[
    \mu \II_{y_\mathfrak{h}} \cdot \left( \frac{1}{3} \id -  \nabla y_\mathfrak{h}^\top n_\mathfrak{h} \otimes \nabla y_\mathfrak{h}^\top n_\mathfrak{h} \right)
  \right] \mathrm{d}x \\
  & + \frac{\bar{r}}{8} \int_S \mathcal{I}_\mathfrak{h} \left[
    \bar{\lambda} (\Delta_\mathfrak{h} y_\mathfrak{h} \cdot b_\mathfrak{h}) \operatorname{tr} \left( \frac{1}{3} \id -  \nabla y_\mathfrak{h}^\top n_\mathfrak{h} \otimes \nabla y_\mathfrak{h}^\top n_\mathfrak{h} \right)
  \right] \mathrm{d}x \\
  & + \frac{\bar{r}^2}{16} \int_S
    \mathcal{I}_\mathfrak{h} \left[
      \left( \mu + \bar{\lambda} \right) \operatorname{tr} \left(
        \nabla y_\mathfrak{h}^\top n_\mathfrak{h} \otimes \nabla y_\mathfrak{h}^\top n_\mathfrak{h}
      \right)^2
      - \frac{2\mu + 4 \bar{\lambda}}{3} \operatorname{tr} \left(
        \nabla y_\mathfrak{h}^\top n_\mathfrak{h} \otimes \nabla y_\mathfrak{h}^\top n_\mathfrak{h}
      \right)
    \right]
  \mathrm{d}x \\
  & + \frac{\bar{r}^2}{72} (\mu + 2\bar{\lambda}) |S|.
\end{aligned}
\end{equation*}
Derivatives of $N_{\bar{r},\mathfrak{h}}$ at $(y_\mathfrak{h},n_\mathfrak{h})$ will be denoted with $\partial_y N_{\bar{r},\mathfrak{h}}(y_\mathfrak{h},n_\mathfrak{h})$ and $\partial_n N_{\bar{r},\mathfrak{h}}(y_\mathfrak{h},n_\mathfrak{h})$.

\subsection{Numerical minimization via discrete gradient flow}\label{subsec:gradflow}
We employ a semi-implicit discrete gradient flow scheme for the approximation of stationary points of the discretized energy. During the iteration each (pseudo-)time step is restricted to the corresponding tangent space of the discrete admissible set at the current iterate. The tangent spaces are obtained from linearizing the unit-length and isometry constraints and are defined at $n_\mathfrak{h} \in \mathcal{A}_{n,\mathfrak{h}}$, $y_\mathfrak{h} \in \mathcal{A}_{y,\mathfrak{h}}$ via
\begin{align*}
  \mathcal{F}_{n,\mathfrak{h}}(n_\mathfrak{h}) \colonequals  \left\{
    s_\mathfrak{h} \in \mathcal{S}^1(\mathcal{T}_\mathfrak{h}; \R^3) : s_\mathfrak{h}(z) = 0 \text{ for all } z \in \mathcal{N}_\mathfrak{h} \cap \Gamma_n,
    \; s_\mathfrak{h}(z) \cdot n_\mathfrak{h}(z) = 0 \text{ for all } z \in \mathcal{N}_\mathfrak{h}
  \right\}
\end{align*}
and
\begin{align*}
  \mathcal{F}_{y,\mathfrak{h}}(y_\mathfrak{h}) \colonequals  \Big\{
    w_\mathfrak{h} \in \mathcal{S}^\mathrm{DKT}(\mathcal{T}_\mathfrak{h}; \R^3) : w_\mathfrak{h}(z) = 0,\; \nabla w_\mathfrak{h}(z) = 0
    \text{ for all } z \in \mathcal{N}_\mathfrak{h} \cap \Gamma_y, \\
    \left(\nabla w_\mathfrak{h}(z)\right)^\top \nabla y_\mathfrak{h}(z) + (\nabla y_\mathfrak{h}(z))^\top \nabla w_\mathfrak{h}(z) = 0  \text{ for all } z \in \mathcal{N}_\mathfrak{h}
  \Big\}.
\end{align*}

We now define the discrete gradient flow via semi-implicitly approximating the variation of the energy with respect to the deformation and the director field, relying on an implicit treatment of the convex quadratic terms while making use of an explicit treatment of nonlinear parts.
We write
\[
H_{\Gamma}^1(S; \R^3) \colonequals  \left\{ v \in H^1(S; \R^3) : v|_{\Gamma} = 0 \right\}
\]
for the set of functions in $H^1(S; \R^3)$ that vanish on a closed, non-zero measured subset $\Gamma \subseteq \partial \Omega$ in the sense of traces.
The scalar products on $\mathcal{S}^1_{\Gamma_n}(\mathcal{T}_\mathfrak{h}; \R^3) \colonequals  \mathcal{S}^1(\mathcal{T}_\mathfrak{h}; \R^3) \cap H_{\Gamma_n}^1(S; \R^3)$ and $\mathcal{S}^\mathrm{DKT}_{\Gamma_y}(\mathcal{T}_\mathfrak{h}; \R^3) \colonequals  \mathcal{S}^\mathrm{DKT}(\mathcal{T}_\mathfrak{h}; \R^3) \cap H_{\Gamma_y}^1(S;\R^3)$ are denoted with $(\,\cdot\,,\cdot\,)_*$ and $(\,\cdot\,,\cdot\,)_{**}$, respectively.
\begin{algorithm}\label{alg:gradflow}
  INPUT: initial values $y^0 \in \mathcal{A}_{y,\mathfrak{h}}$, $n^0 \in \mathcal{A}_{n,\mathfrak{h}}$, step size $\tau > 0$, stopping criterion $\e_\mathrm{stop}>0$. \\
  (1) Initialize $k \colonequals  1$.\\
  (2) Compute $y_\mathfrak{h}^k \colonequals  y_\mathfrak{h}^{k-1} + \tau d_t y_\mathfrak{h}^k$ with $d_t y_\mathfrak{h}^k \in \mathcal{F}_{y,\mathfrak{h}}(y_\mathfrak{h}^{k-1})$ such that
  \[
    (d_t y_\mathfrak{h}^k, w_\mathfrak{h})_{**} =
    -\frac{\mu + \bar{\lambda}}{6}
    \left(
      \nabla \Theta_\mathfrak{h} y_\mathfrak{h}^{k}, \nabla \Theta_\mathfrak{h} w_\mathfrak{h}
    \right)_{L^2(S)}
    -\partial_{y} N_{\bar{r},\mathfrak{h}} \left( y_\mathfrak{h}^{k-1},n_\mathfrak{h}^{k-1} \right)[w_\mathfrak{h}]
  \]
  for all $w_\mathfrak{h} \in \mathcal{F}_{y,\mathfrak{h}}(y_\mathfrak{h}^{k-1})$. \\
  (3) Compute $n_\mathfrak{h}^k \colonequals  n_\mathfrak{h}^{k-1} + \tau d_t n_\mathfrak{h}^k$ with $d_t n_\mathfrak{h}^k \in \mathcal{F}_{n,\mathfrak{h}}(n_\mathfrak{h}^{k-1})$ such that
  \[
    (d_t n_\mathfrak{h}^k, s_\mathfrak{h})_{*} =
    -\bar{\e}^2 \left( \nabla n_\mathfrak{h}^k, \nabla s_\mathfrak{h}^k \right)_{L^2(S)}
    -\partial_{n} N_{\bar{r},\mathfrak{h}}\left( y_\mathfrak{h}^{k}, n_\mathfrak{h}^{k-1}\right)[s_\mathfrak{h}]
  \]
  for all $s_\mathfrak{h} \in \mathcal{F}_{n,\mathfrak{h}}(n_\mathfrak{h}^{k-1})$. \\
  (4) If $\left( (d_t y_\mathfrak{h}^k,d_t y_\mathfrak{h}^k)_{**} + (d_t n_\mathfrak{h}^k,d_t n_\mathfrak{h}^k)_* \right)^{1/2} < \e_\mathrm{stop}$, stop the iteration at the approximate equilibrium $(y_\mathfrak{h}^\infty, n_\mathfrak{h}^\infty) \colonequals  (y_\mathfrak{h}^k, n_\mathfrak{h}^k)$. Else, increase $k$ via $k \mapsto k+1$ and continue with~(2).
\end{algorithm}

\subsection{Numerical experiments}\label{subsec:numexp}
The experiments presented here serve to illustrate the practical properties of the proposed scheme as well as characteristic features of the LCE model and its energy landscape.
Its stability as well as the convergence of the discretization are currently a work in progress and will be presented in a follow-up paper.
In the following we denote with $(y_\mathfrak{h}^\infty, n_\mathfrak{h}^\infty)$ the numerical equilibrium state obtained from the discrete gradient flow.
Since in general, we do not know analytical minimizers nor their minimal energies we use the experimental order of convergence,
\[
\mathrm{EOC}^{\mathcal{E}}_\mathfrak{h} \colonequals  \log_2 \left(
  \frac{
    | \mathcal{E}_{4\mathfrak{h}}(y_{4\mathfrak{h}}^\infty, n_{4\mathfrak{h}}^\infty)
    - \mathcal{E}_{2\mathfrak{h}}(y_{2\mathfrak{h}}^\infty, n_{2\mathfrak{h}}^\infty) |
  } {
    | \mathcal{E}_{2\mathfrak{h}}(y_{2\mathfrak{h}}^\infty, n_{2\mathfrak{h}}^\infty)
    - \mathcal{E}_\mathfrak{h}(y_\mathfrak{h}^\infty, n_\mathfrak{h}^\infty) |
  }
\right),
\]
as an approximate measure of the convergence rate for the discrete energies.
Furthermore, we let
\[
\mathrm{Err}_\mathrm{1}(n_\mathfrak{h}) \colonequals  \max_{z \in \mathcal{N}_\mathfrak{h}} \big| |n_\mathfrak{h}| - 1 \big|
\]
and
\[
\mathrm{Err}_\mathrm{iso}(y_\mathfrak{h}) \colonequals  \max_{z \in \mathcal{N}_\mathfrak{h}} \big| \mathcal{I}_\mathfrak{h}[\nabla y_\mathfrak{h}^\top \nabla y_\mathfrak{h}] - \Id \big|
\]
measure violations of the unit-length constraint and isometry constraint that occur in the discrete gradient flow scheme.
Finally, we denote with
\[
\mathcal{E}_{\mathfrak{h},\mathrm{OF}} \colonequals  \frac{\bar{\e}^2}{2} \int_S |\nabla n_\mathfrak{h}^\infty|^2 \mathrm{d}x
\]
the Oseen-Frank energy of the numerical equilibrium configurations and, occasionally, use the shorthand notation $\mathcal{E}_\mathfrak{h}\colonequals  \mathcal{E}_\mathfrak{h}(y_\mathfrak{h}^\infty, n_\mathfrak{h}^\infty)$ for the corresponding total energy.
The surface coloring in the following plots always corresponds to the potential of the Oseen-Frank energy.

Note that the relaxation of the constraints that results from the linearization in the discrete gradient flow effectively enlarges the admissible sets.
This may lead to discrete minimizers with lower energies compared to minimizers that one would obtain if $\zeta_\mathfrak{h} = 0$, i.\,e. if the nodal constraints were satisfied exactly.
The experimental data supports the hypothesis that constraint violation is controlled by the time step size which we usually choose proportional to the mesh size.
As a consequence in practice, the admissible sets $\mathcal{A}_{n,\mathfrak{h}/2}$, $\mathcal{A}_{y,\mathfrak{h}/2}$ that correspond to a triangulation $\mathcal{T}_{\mathfrak{h}/2}$ after one red-refinement of a mesh are, in general,
no subsets of the admissible sets $\mathcal{A}_{n,\mathfrak{h}}$, $\mathcal{A}_{y,\mathfrak{h}}$ that correspond to the unrefined triangulation $\mathcal{T}_{\mathfrak{h}}$.
This contributes to the reasons why the observed experimental convergence orders that we observe in the following are not always consistent.
Another reason lies in the fact that asymptotic ranges might not be reached due to our limited computational resources: all computations were performed on a standard Intel\textsuperscript{\textregistered} Core\textsuperscript{TM} i3-3240 CPU @\,3.40GHz\,$\times$\,4 desktop computer.

As a final remark before presenting the experiments, we note that in the examples the choice of required parameters is motivated by realistic materials, but adjusted to emphasize the interaction mechanisms of our model, cf.~Section~\ref{S:modeling}.
An extensive research on what materials are technically possible lies beyond the scope of this work.

\begin{example}[Horizontally clamped plate]\label{ex:clampedplate}
  In this first example we let $\mu = 1$, $\lambda = 1000$ and consider model parameters $\bar{r}, \bar{\e} \in \{1,5\}$. These choices pronounce characteristic model traits. Taking into account the considerations in Section~\ref{S:modeling}, we also perform the calculations for $\bar{r} = 1$, $\bar{\e} = 0.01$.

  We consider the two-dimensional reference configuration $S = (-1,1) \times (-1,1)$ with Dirichlet boundary $\Gamma = \Gamma_y = \Gamma_n = \{-1\} \times [-1,1]$. We prescribe clamped, affine boundary conditions,
  \[
  y_\mathrm{BC}(x) = (x,0)^\top, \; \Phi_\mathrm{BC} = (\Id, 0)^\top,
  \]
  for the deformation $y$ and one of two different strong Dirichlet anchorings,
  \[
  n^\mathrm{x}_\mathrm{BC} = (1,0,0)^\top \text{, \quad or \quad} n^\mathrm{y}_\mathrm{BC} = (0,1,0)^\top,
  \]
  for the director $n$.
  For triangulations $\mathcal{T}_\mathfrak{h} = \mathcal{T}_k$, $k=3,4,5,6$, of $S$ into halved squares with side lengths $\hat{\mathfrak{h}} = 2^{-k}$ we compute approximations of stationary points of the energy~\eqref{eq:energy-discrete1} using the alternating discrete gradient flow Algorithm~\ref{alg:gradflow} with stopping criterion $\e_\mathrm{stop} = 2 \times 10^{-3}$.
  Initial values are chosen as the flat deformation $y_\mathfrak{h}^0(x) = (x_1, x_2, 0)$ and, depending on the anchoring, $n_\mathfrak{h}^0 \equiv n^\mathrm{x}_\mathrm{BC}$ or $n_\mathfrak{h}^0 \equiv n^\mathrm{y}_\mathrm{BC}$.
  The finest triangulation $\mathcal{T}_6$ consists of 32'768 triangles and corresponding discrete deformations and director fields are determined by 149'769 and 49'923 degrees of freedom, respectively.
  For each triangulation we choose the (pseudo-) time step size $\tau = \hat{\mathfrak{h}}/10$.
  An exception to this choice of step size are the settings $\bar{r} = 1$, $\bar{\e} = 0.01$, where lower energies and constraint violations justify large step sizes $\tilde{\tau} = 4\hat{\mathfrak{h}}$ to avoid unnecessary high iteration numbers.
  The resulting deformations for $k=6$ are depicted in Figure~\ref{fig:ex1}.
  Iteration numbers, final energies, constraint violations and experimental convergence rates are listed in Tables~\ref{tab:clampedplatex} and~\ref{tab:clampedplatey} for strong anchorings $n^\mathrm{x}_\mathrm{BC}$ and $n^\mathrm{y}_\mathrm{BC}$, respectively.
  Notably, with the anchoring $n^\mathrm{y}_\mathrm{BC}$, the experimental order of convergence of the discrete energies is negative in many cases.
  This might be a consequence of relaxing of the node-wise unit-length and isometry constraints, as well as it might indicate that the asymptotic range has not been reached in terms of grid size.
  Violations of the point-wise unit-length and isometry constraints appear to decay linearly with respect to $\tau \sim \mathfrak{h}$.

  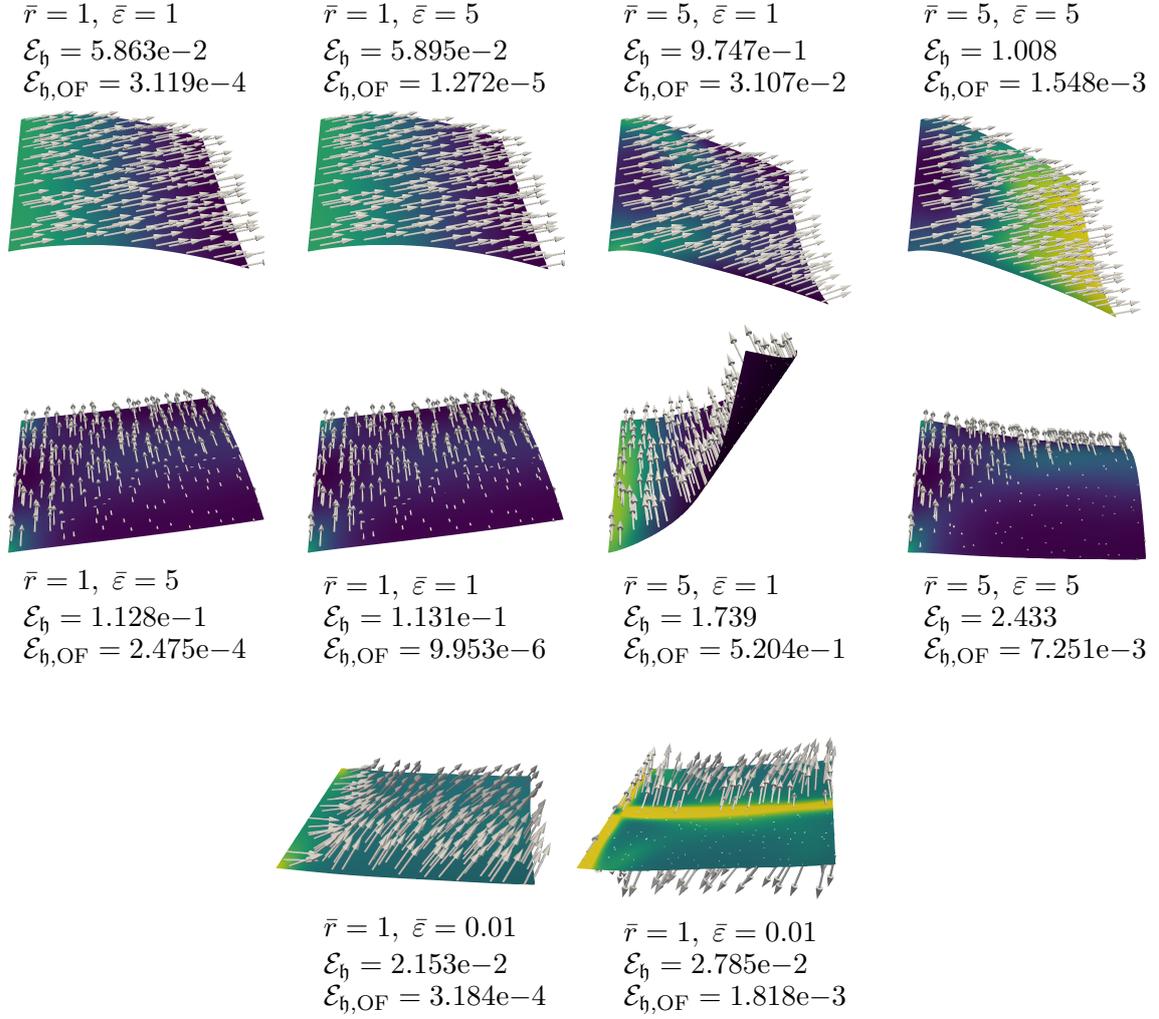
\begin{figure}
  \centering
  \begin{subfigure}{0.24\textwidth}
    \centering
    \includegraphics[trim={160pt 20pt 150pt 20pt},clip,width=\textwidth]{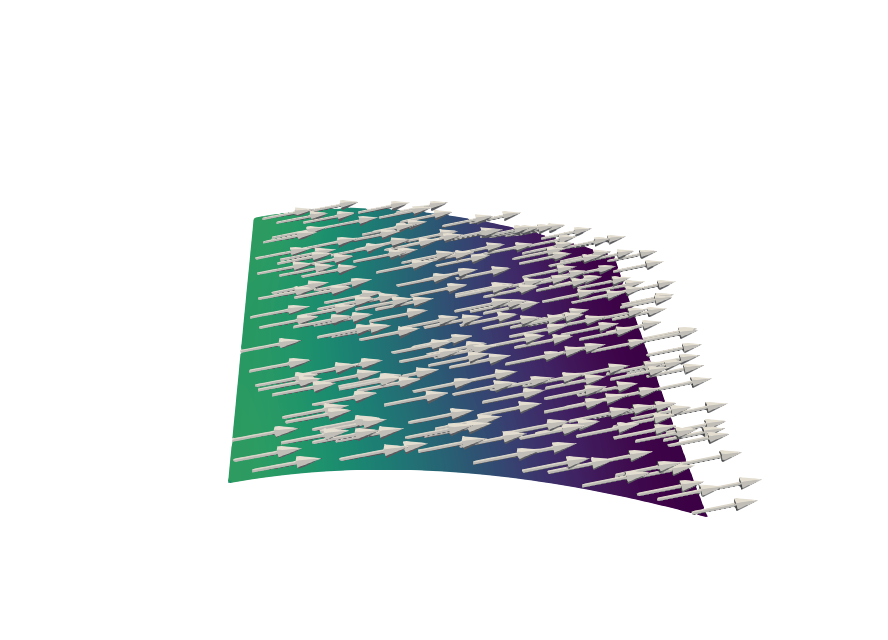}
    \put(-90,112){$\bar{r}=1,\;\bar{\varepsilon}=1$}
    \put(-90,98){$\mathcal{E}_\mathfrak{h}=\num{0.058628}$}
    \put(-90,86){$\mathcal{E}_{\mathfrak{h},\mathrm{OF}}=\num{0.000311899}$}
  \end{subfigure}
  \begin{subfigure}{0.24\textwidth}
    \centering
    \includegraphics[trim={160pt 20pt 150pt 20pt},clip,width=\textwidth]{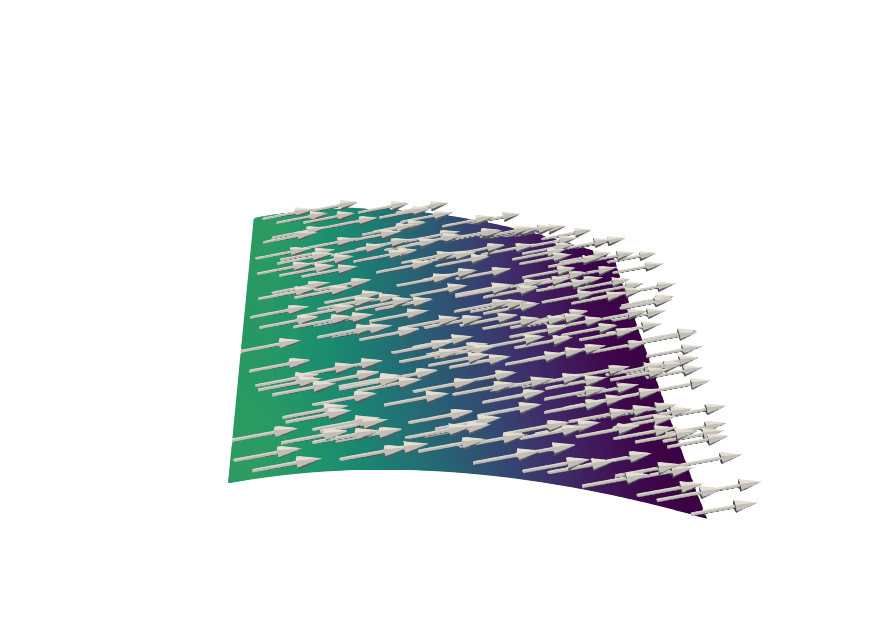}
    \put(-90,112){$\bar{r}=1,\;\bar{\varepsilon}=5$}
    \put(-90,98){$\mathcal{E}_\mathfrak{h}=\num{0.0589452}$}
    \put(-90,86){$\mathcal{E}_{\mathfrak{h},\mathrm{OF}}=\num{1.2724e-05}$}
  \end{subfigure}
  \begin{subfigure}{0.24\textwidth}
    \centering
    \includegraphics[trim={160pt 20pt 150pt 20pt},clip,width=\textwidth]{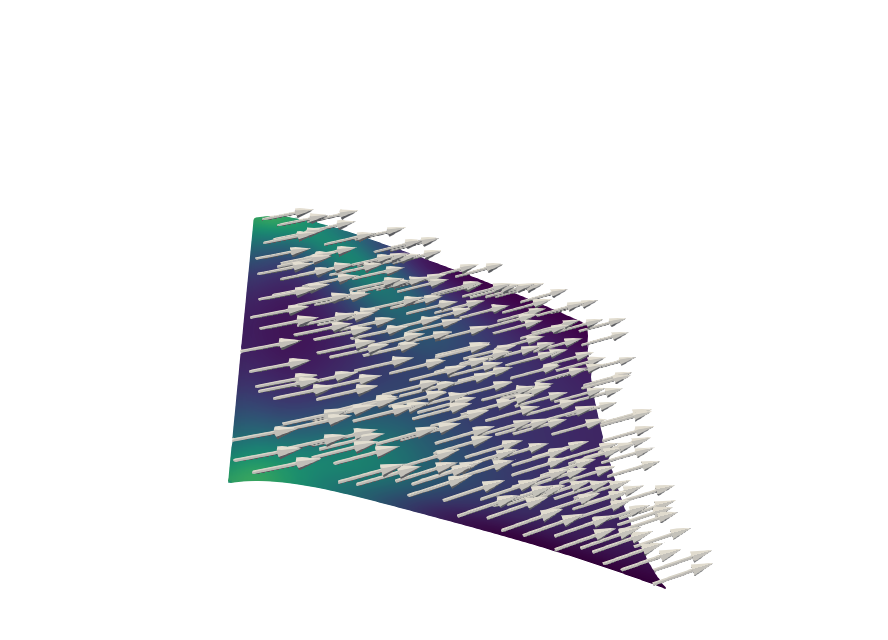}
    \put(-90,112){$\bar{r}=5,\;\bar{\varepsilon}=1$}
    \put(-90,98){$\mathcal{E}_\mathfrak{h}=\num{0.974691}$}
    \put(-90,86){$\mathcal{E}_{\mathfrak{h},\mathrm{OF}}=\num{0.0310667}$}
  \end{subfigure}
  \begin{subfigure}{0.24\textwidth}
    \centering
    \includegraphics[trim={160pt 20pt 150pt 20pt},clip,width=\textwidth]{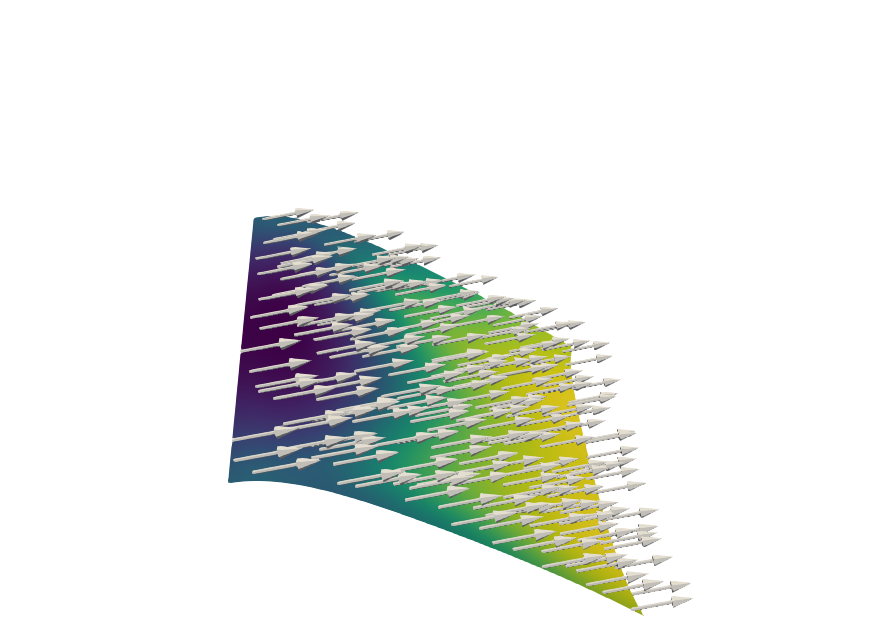}
    \put(-90,112){$\bar{r}=5,\;\bar{\varepsilon}=5$}
    \put(-90,98){$\mathcal{E}_\mathfrak{h}=\num{1.00828}$}
    \put(-90,86){$\mathcal{E}_{\mathfrak{h},\mathrm{OF}}=\num{0.00154843}$}
  \end{subfigure}
  \\
  \begin{subfigure}{0.24\textwidth}
    \centering
    \includegraphics[trim={160pt 20pt 150pt 20pt},clip,width=\textwidth]{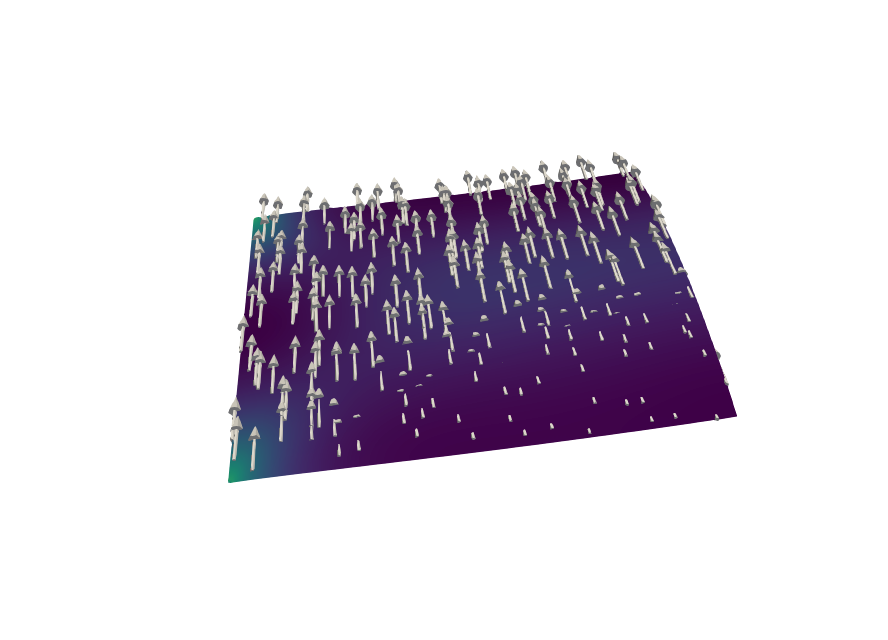}
    \put(-90,12){$\bar{r}=1,\;\bar{\varepsilon}=5$}
    \put(-90,-2){$\mathcal{E}_\mathfrak{h}=\num{0.112841}$}
    \put(-90,-14){$\mathcal{E}_{\mathfrak{h},\mathrm{OF}}=\num{0.000247541}$}
  \end{subfigure}
  \begin{subfigure}{0.24\textwidth}
    \centering
    \includegraphics[trim={160pt 20pt 150pt 20pt},clip,width=\textwidth]{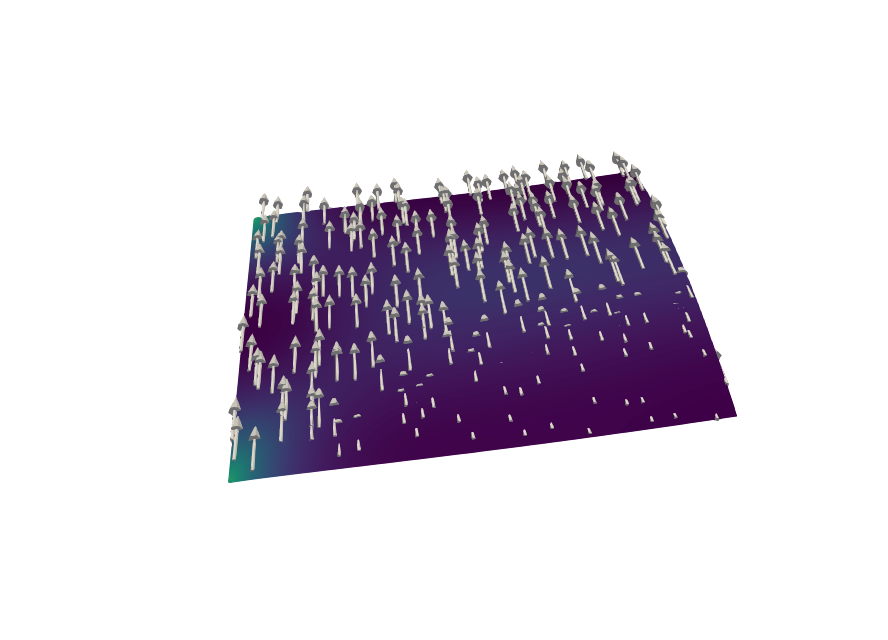}
    \put(-90,10){$\bar{r}=1,\;\bar{\varepsilon}=1$}
    \put(-90,-2){$\mathcal{E}_\mathfrak{h}=\num{0.113083}$}
    \put(-90,-14){$\mathcal{E}_{\mathfrak{h},\mathrm{OF}}=\num{9.95273e-06}$}
  \end{subfigure}
  \begin{subfigure}{0.24\textwidth}
    \centering
    \includegraphics[trim={160pt 20pt 150pt 20pt},clip,width=\textwidth]{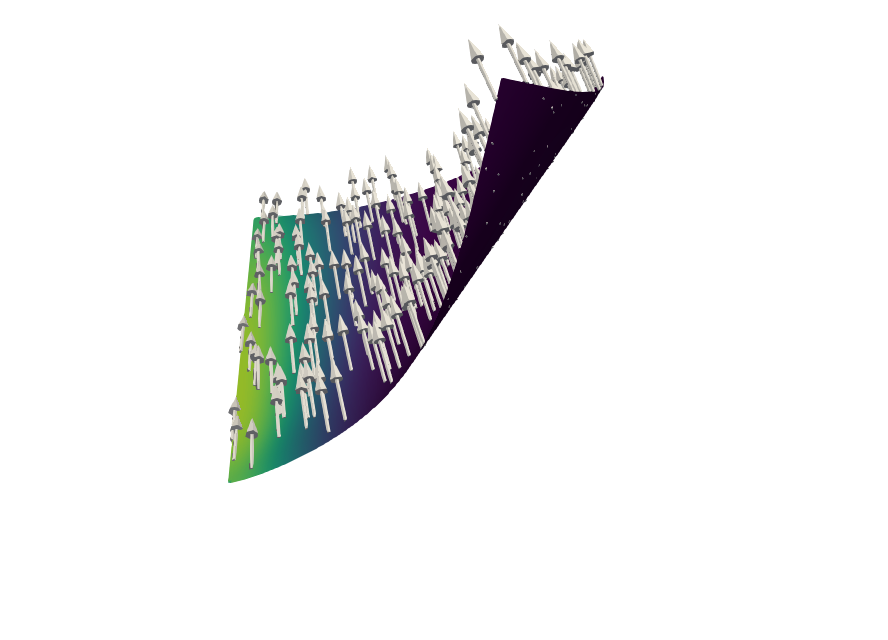}
    \put(-90,10){$\bar{r}=5,\;\bar{\varepsilon}=1$}
    \put(-90,-2){$\mathcal{E}_\mathfrak{h}=\num{1.73903}$}
    \put(-90,-14){$\mathcal{E}_{\mathfrak{h},\mathrm{OF}}=\num{0.520428}$}
  \end{subfigure}
  \begin{subfigure}{0.24\textwidth}
    \centering
    \includegraphics[trim={160pt 20pt 150pt 20pt},clip,width=\textwidth]{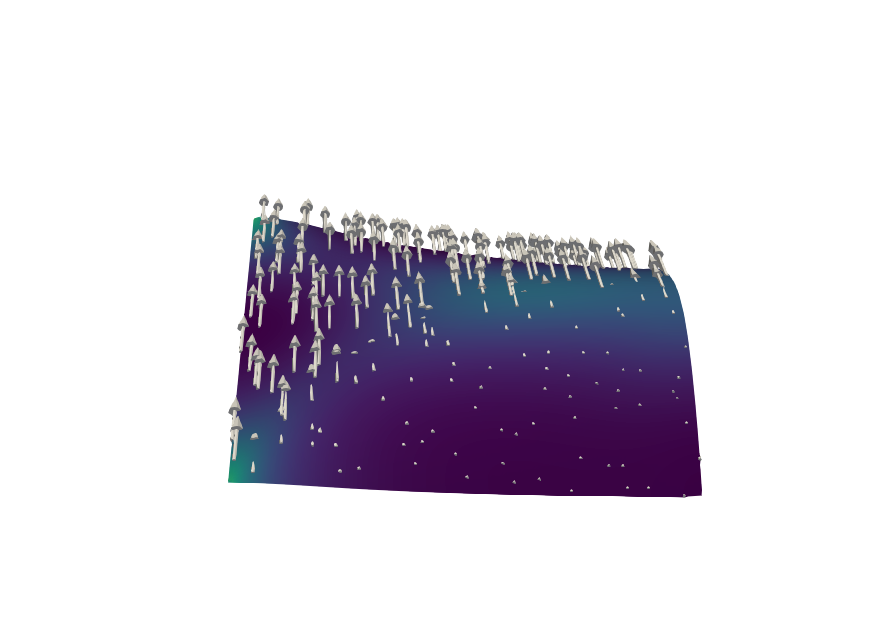}
    \put(-90,10){$\bar{r}=5,\;\bar{\varepsilon}=5$}
    \put(-90,-2){$\mathcal{E}_\mathfrak{h}=\num{2.43331}$}
    \put(-90,-14){$\mathcal{E}_{\mathfrak{h},\mathrm{OF}}=\num{0.00725076}$}
  \end{subfigure}
  \\
  \begin{subfigure}{0.24\textwidth}
    \centering
    \includegraphics[trim={160pt 20pt 150pt 20pt},clip,width=\textwidth]{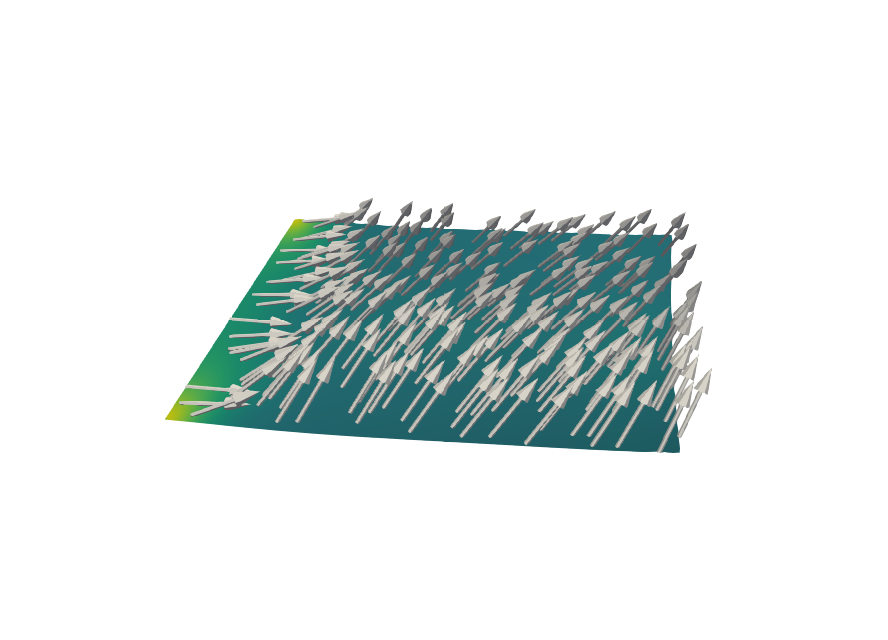}
    \put(-90,12){$\bar{r}=1,\;\bar{\varepsilon}=0.01$}
    \put(-90,-2){$\mathcal{E}_\mathfrak{h}=\num{0.0215294}$}
    \put(-90,-14){$\mathcal{E}_{\mathfrak{h},\mathrm{OF}}=\num{0.000318373}$}
  \end{subfigure}
  \begin{subfigure}{0.24\textwidth}
    \centering
    \includegraphics[trim={160pt 20pt 150pt 20pt},clip,width=\textwidth]{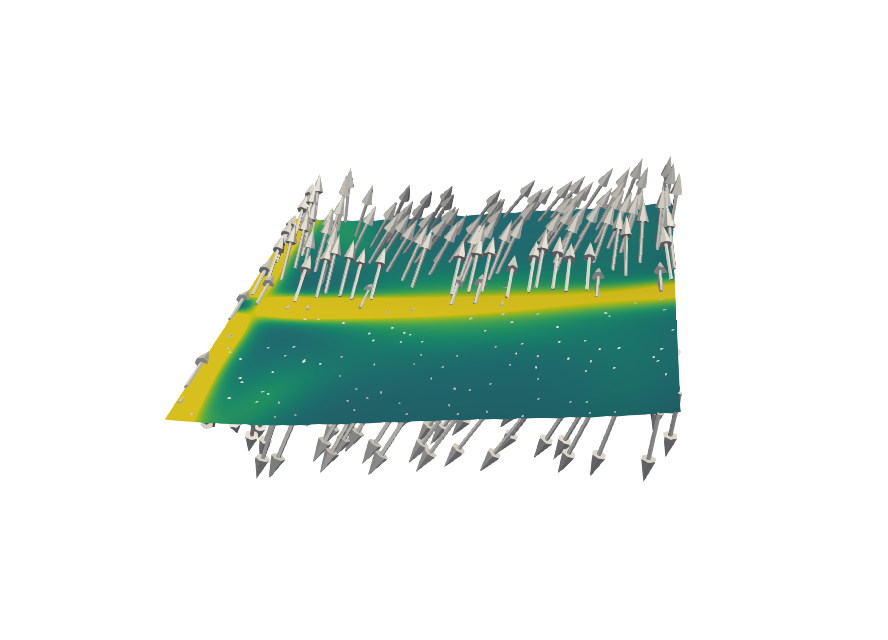}
    \put(-90,10){$\bar{r}=1,\;\bar{\varepsilon}=0.01$}
    \put(-90,-2){$\mathcal{E}_\mathfrak{h}=\num{0.0278454}$}
    \put(-90,-14){$\mathcal{E}_{\mathfrak{h},\mathrm{OF}}=\num{0.00181775}$}
  \end{subfigure}
  \caption{Deformations and director fields of a horizontally clamped LCE plate for different strong Dirichlet anchorings on the clamped boundary with varying model parameters in Example~\ref{ex:clampedplate};
  top row: perpendicular anchoring,
  middle row: parallel anchoring (with respect to the clamped boundary),
  bottom row: perpendicular (left) and parallel (right) anchoring.}
  \label{fig:ex1}
  \end{figure}

  \begin{table}
    \centering
    \begin{subtable}{.8\textwidth}
      \centering
      \begin{tabular}{ l c c r c c c c }
        $k$ & $\bar{r}$ & $\bar{\e}$ & \#iter & $\mathcal{E}_\mathfrak{h}$
        & $\mathrm{Err}_\mathrm{1}(n_\mathfrak{h}^\infty)$
        & $\mathrm{Err}_\mathrm{iso}(y_\mathfrak{h}^\infty)$
        & $\mathrm{EOC}_\mathfrak{h}^\mathcal{E}$ \\\hline\hline
        3 & 1 & 1 & 1103  & \num{0.0575931} & \num{0.000625947}
          & \num{0.00123226} & --- \\
        4 &   &   & 2159  & \num{0.0579885} & \num{0.00031294}
          & \num{0.000615916} & --- \\
        5 &   &   & 4241  & \num{0.0583321} & \num{0.000156418}
          & \num{0.000307807} & \num{0.202583} \\
        6 &   &   & 8333 & \num{0.058628}  & \num{7.8155e-05}
          & \num{0.000153781} & \num{0.21562} \\\hline
        3 & 1 & 5 & 1114  & \num{0.0580031} & \num{0.000629494}
          & \num{0.00125931} & --- \\
        4 &   &   & 2188  & \num{0.0583658} & \num{0.000314702}
          & \num{0.000629463} & --- \\
        5 &   &   & 4310  & \num{0.0586772} & \num{0.000157283}
          & \num{0.00031457} & \num{0.220008} \\
        6 &   &   & 8494 & \num{0.0589452} & \num{7.85819e-05}
          & \num{0.00015716} & \num{0.216536} \\\hline
        3 & 5 & 1 & 1487  & \num{0.91406}   & \num{0.012742}
          & \num{0.0181524} & --- \\
        4 &   &   & 2845  & \num{0.941151}  & \num{0.00625554}
          & \num{0.00891069} & --- \\
        5 &   &   & 5871 & \num{0.959666}  & \num{0.0030998}
          & \num{0.00441443} & \num{0.549119} \\
        6 &   &   & 12858 & \num{0.974691}  & \num{0.00154266}
          & \num{0.00219623} & \num{0.30133} \\\hline
        3 & 5 & 5 & 1446  & \num{0.946309}  & \num{0.0130965}
          & \num{0.0252902} & --- \\
        4 &   &   & 2690  & \num{0.97363}   & \num{0.00645727}
          & \num{0.0124195} & --- \\
        5 &   &   & 5266 & \num{0.992333}  & \num{0.00320551}
          & \num{0.00615247} & \num{0.546741} \\
        6 &   &   & 10621 & \num{1.00828} & \num{0.00159636}
          & \num{0.00306065} & \num{0.229985} \\\hline
        4\textsuperscript{*} & 1 & \num{0.01} & 2996  & \num{0.0167833}   & \num{0.0143042}
          & \num{0.023743} & --- \\
        5\textsuperscript{*} &   &   & 5278 & \num{0.0199209}  & \num{0.00703343}
          & \num{0.0118994} & \num{1.6720156034654983} \\
        6\textsuperscript{*} &   &   & 10058 & \num{0.0215294} & \num{0.00346909}
          & \num{0.00594608} & \num{0.9639455046520449} \\
      \end{tabular}
      \caption{Strong Dirichlet anchoring $n_\mathrm{BC}^\mathrm{x}$ (perpendicular to the clamped boundary).}
      \label{tab:clampedplatex}
    \end{subtable}

    \vspace{10pt}
    \begin{subtable}{.8\textwidth}
      \centering
      \begin{tabular}{ l c c c c c c c }
        $k$ & $\bar{r}$ & $\bar{\e}$ & \#iter & $\mathcal{E}_\mathfrak{h}$
        & $\mathrm{Err}_\mathrm{1}(n_\mathfrak{h}^\infty)$
        & $\mathrm{Err}_\mathrm{iso}(y_\mathfrak{h}^\infty)$
        & $\mathrm{EOC}_\mathfrak{h}^\mathcal{E}$ \\\hline\hline
        3 & 1 & 1 & 1045   & \num{0.110921} & \num{9.70513e-05}
          & \num{0.000193903} & --- \\
        4 &   &   & 2085   & \num{0.111095} & \num{4.8009e-05}
          & \num{9.60096e-05} & --- \\
        5 &   &   & 4164   & \num{0.111729} & \num{2.35571e-05}
          & \num{4.72638e-05} & \num{-1.8654} \\
        6 &   &   & 8332  & \num{0.112841} & \num{1.14114e-05}
          & \num{2.30742e-05} & \num{-0.810602} \\\hline
        3 & 1 & 5 & 1051   & \num{0.111197} & \num{9.64072e-05}
          & \num{0.000192905} & --- \\
        4 &   &   & 2079   & \num{0.111361} & \num{4.76998e-05}
          & \num{9.55449e-05} & --- \\
        5 &   &   & 4152   & \num{0.111984} & \num{2.34002e-05}
          & \num{4.70398e-05} & \num{-1.92554} \\
        6 &   &   & 8310  & \num{0.113083} & \num{1.13299e-05}
          & \num{2.29644e-05} & \num{-0.818887} \\\hline
        3 & 5 & 1 & 24726  & \num{1.64894}  & \num{0.00250404}
          & \num{0.00576815}  & --- \\
        4 &   &   & 45836  & \num{1.68264}  & \num{0.0012177}
          & \num{0.00302362}  & --- \\
        5 &   &   & 81777 & \num{1.71372}  & \num{0.000607086}
          & \num{0.00145278}  & \num{0.116762} \\
        6 &   &   & 140265 & \num{1.73903} & \num{0.00029504}
          & \num{0.000657006} & \num{0.296278} \\\hline
        3 & 5 & 5 & 1686   & \num{2.18694}  & \num{0.00253391}
          & \num{0.00508268}  & --- \\
        4 &   &   & 3384   & \num{2.25544}  & \num{0.001236}
          & \num{0.0024606}   & --- \\
        5 &   &   & 6971  & \num{2.33025}  & \num{0.000595522}
          & \num{0.00118686}  & \num{-0.127127} \\
        6 &   &   & 14434  & \num{2.43331}  & \num{0.000279477}
          & \num{0.000567115} & \num{-0.462181} \\\hline
        4\textsuperscript{*} & 1 & \num{0.01} & 6417  & \num{0.0166582}   & \num{0.00320952}
          & \num{0.00392164} & --- \\
        5\textsuperscript{*} &   &   & 11790 & \num{0.0231032}  & \num{0.00148879}
          & \num{0.00194941} & \num{0.43516468826776666} \\
        6\textsuperscript{*} &   &   & 21862 & \num{0.0278454} & \num{0.000693901}
          & \num{0.000955672} & \num{0.4426238495157202} \\
      \end{tabular}
      \caption{Strong Dirichlet anchoring $n_\mathrm{BC}^\mathrm{y}$ (tangential to the clamped boundary).}
      \label{tab:clampedplatey}
    \end{subtable}
    \caption{Experimental results for the clamped plate in Example~\ref{ex:clampedplate} with refinement level $k$ for two different strong Dirichlet anchorings and varying choices of $\bar{r}$ and $\bar{\varepsilon}$. Lines in which the refinement level is marked with an asterisk were obtained using larger step sizes. For the settings $\bar{\e} = 0.01$ we omitted lines for refinement levels $k=3$ due to space restrictions.}
  \end{table}
\end{example}

\begin{example}[Compressive boundary conditions and strong tangential anchoring]\label{ex:compress}
  We let $\mu=1$, $\lambda=1000$ and $\bar{r} = 4$. This example investigates stationary configurations of a rectangular strip with reference configuration $S = (-5,5) \times (-1,1)$ subject to compressive clamped, affine boundary conditions
  \[
  y_\mathrm{BC}(x) = (\alpha x_1, x_2, 0)^\top, \; \Phi_\mathrm{BC} = (\Id, 0)^\top,
  \]
  with compression factor $0 < \alpha < 1$ and $\bar{\e} \in \{1, 4\}$ on $\Gamma_y = \{-5,5\} \times [-1,1]$.
  Throughout the domain we restrict the director to tangential directions via imposing the additional (node-wise) strong tangential anchoring constraint $(\partial_1 y \wedge \partial_2 y) \cdot n = 0$, cf.~Remark~\ref{R:Bridge:Analysis:Numerics}(a).
  The compression forces the plate to buckle up or buckle down, leading to stationary configurations with different local energy minima.
  We consider a triangulation of $S$ into 10'240 halved squares of side-length $\hat{\mathfrak{h}}=2^{-4}$ and use the step size $\tau = \hat{\mathfrak{h}}/4$ in the discrete gradient flow, leading to approximate stationary configurations with deformations and director fields determined by 47'817 and 15'939 degrees of freedom, respectively.
  The stopping criterion $\e_\mathrm{stop} = 2 \times 10^{-3}$ was reached after a few thousand iterations in all cases and the largest observed constraint violations were $\mathrm{Err}_\mathrm{1}(n_\mathfrak{h}^\infty) = 0.0565892$ and $\mathrm{Err}_\mathrm{iso}(y_\mathfrak{h}^\infty) = 0.0185882$, both of which were observed in the most ``extreme'' setting corresponding to $\bar{\e} = 4$ and $\alpha \approx 0.637$.
  In most considered cases the constraint violations were smaller by at least one order of magnitude.
  We use corresponding initial values $y_\mathfrak{h}^0$ in the discrete gradient flow to steer the evolution towards one of the stationary configurations and observe that strong Dirichlet anchoring on $\Gamma_n = \Gamma_y$, namely
  \[
  n^\mathrm{x}_\mathrm{BC} = (1,0,0)^\top \text{ (``perpendicular''), \quad or \quad} n^\mathrm{y}_\mathrm{BC} = (0,1,0)^\top \text{ (``tangential'')},
  \]
  can be used to render either deformation -- ``up'' or ``down'' -- more preferable energy wise.
  Experimental results are collected in Table~\ref{tab:ex2}.
  Figure~\ref{fig:ex2_1} shows the stationary configurations approximating local minima of the energy for $\alpha=3\pi^{-1}\sin(\pi/3)\approx0.827$ and both choices of strong Dirichlet anchorings.
  \begin{figure}
  \centering
  \begin{subfigure}{0.45\textwidth}
    \centering
    \includegraphics[width=\textwidth]{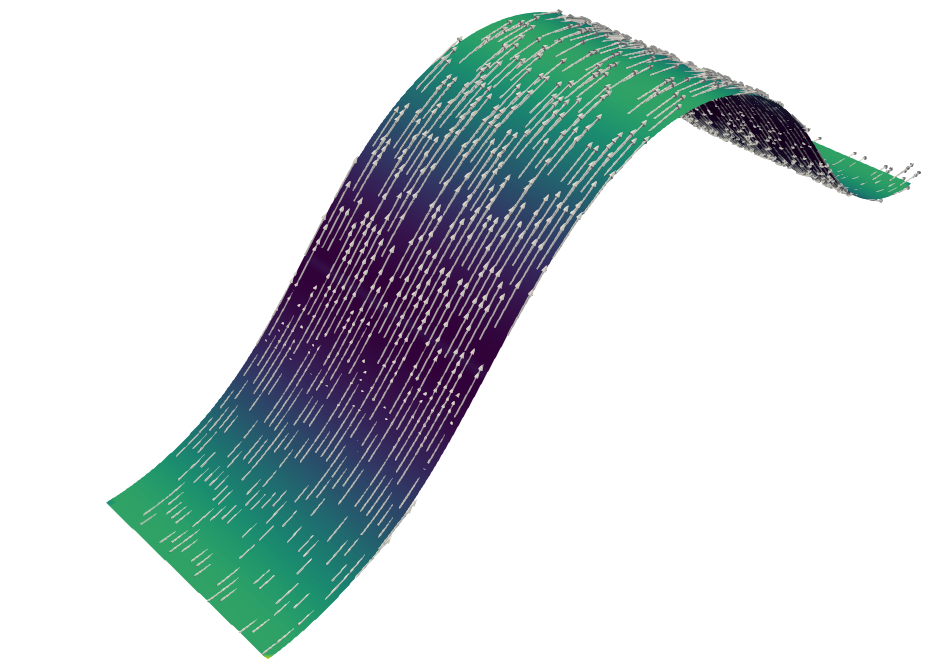}
    \put(-201,135){$\mathcal{E}_\mathfrak{h}=15.23$}
    \put(-201,123){$\mathcal{E}_{\mathfrak{h},\mathrm{OF}}=1.424$}
    \caption{``up'' with perpendicular anchoring $n^\mathrm{x}_\mathrm{BC}$.}
  \end{subfigure}
  \hspace{5mm}
  \begin{subfigure}{0.45\textwidth}
    \centering
    \includegraphics[width=\textwidth]{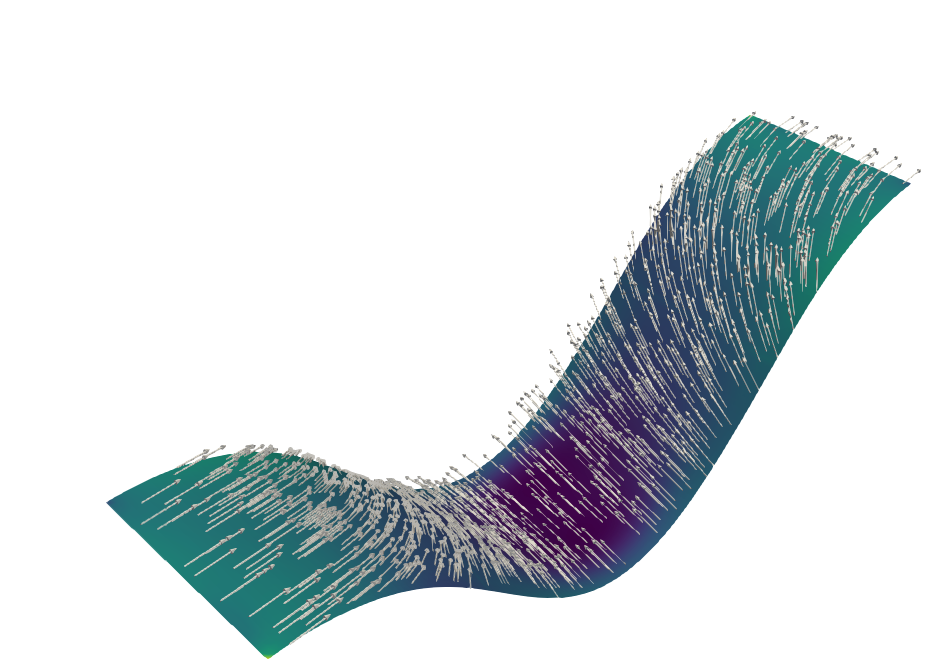}
    \put(-201,135){$\mathcal{E}_\mathfrak{h}=13.73$}
    \put(-201,123){$\mathcal{E}_{\mathfrak{h},\mathrm{OF}}=1.732$}
    \caption{``down'' with perpendicular anchoring $n^\mathrm{x}_\mathrm{BC}$.}
  \end{subfigure}
  \\[5mm]
  \begin{subfigure}{0.45\textwidth}
    \centering
    \includegraphics[width=\textwidth]{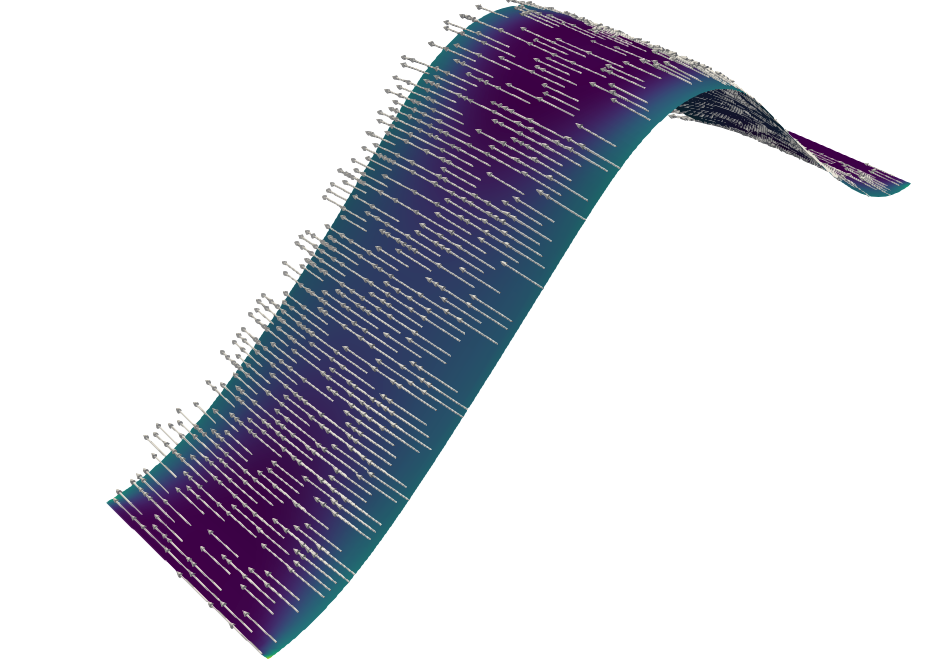}
    \put(-201,135){$\mathcal{E}_\mathfrak{h}=12.74$}
    \put(-201,123){$\mathcal{E}_{\mathfrak{h},\mathrm{OF}}=0.7502$}
    \caption{``up'' with tangential anchoring $n^\mathrm{y}_\mathrm{BC}$.}
  \end{subfigure}
  \hspace{5mm}
  \begin{subfigure}{0.45\textwidth}
    \centering
    \includegraphics[width=\textwidth]{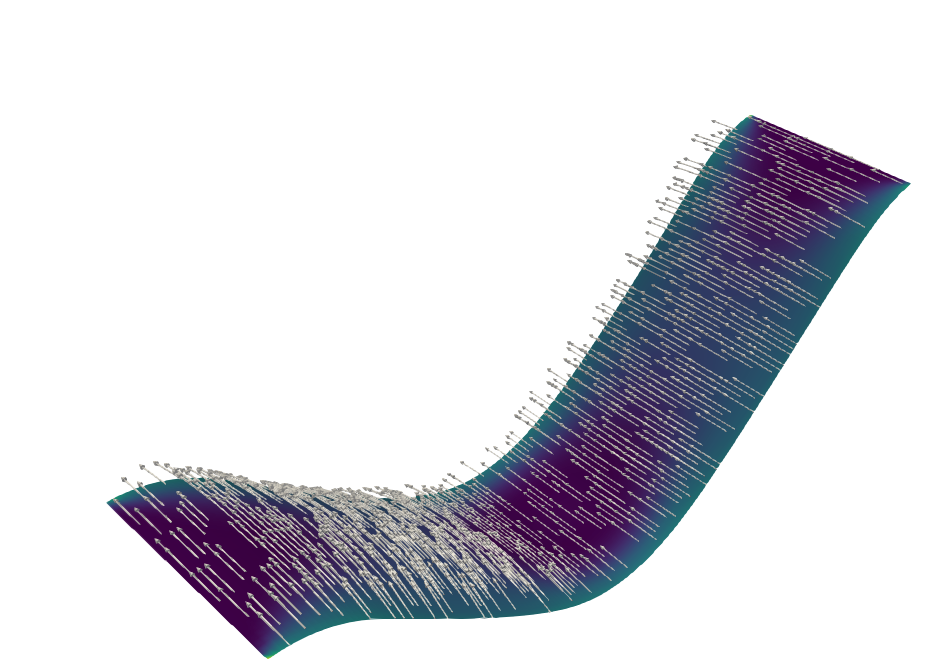}
    \put(-201,135){$\mathcal{E}_\mathfrak{h}=12.74$}
    \put(-201,123){$\mathcal{E}_{\mathfrak{h},\mathrm{OF}}=0.7506$}
    \caption{``down'' with tangential anchoring $n^\mathrm{y}_\mathrm{BC}$.}
  \end{subfigure}
  \vspace{-2mm}
  \caption{Discrete total energies and Oseen-Frank-Energies of the depicted configurations in the case $\bar{r} = 4$, $\bar{\varepsilon} = 1$ of Example~\ref{ex:compress} with strong tangential anchoring with compression factor $\alpha\approx0.827$. The strong Dirichlet anchoring of $n$ on the clamped sides influences which deformation of a compressed strip is preferable energy wise.}
  \label{fig:ex2_1}
  \end{figure}

  Strong Dirichlet anchorings may or may not influence energy response to varying compression factors~$\alpha$, see Figure~\ref{fig:ex2_2} where we compare the total energies of the equilibrium states for three different compression factors.
  With certain combinations of coupling and order parameters $\bar{r}$ and $\bar{\e}$ Dirichlet anchorings can be used to select between a notable energy response or almost no response to compression.
  This is illustrated in Figure~\ref{fig:ex2_2}b for $\bar{r}=4, \bar{\e}=4$, where a director anchored perpendicularly to the clamped boundary corresponds to a significant energy response whereas the energy for all three compression factors stays almost the same in the case of a tangentially anchored director.
  Higher Oseen-Frank energies play an important role in increasing the energy response when using perpendicular anchoring, see Table~\ref{tab:ex2}.
  In this experiment initial configurations in the discrete gradient flow were chosen as interpolants of a continuous deformation consisting of four circular arcs, each of arc length $2.5$, glued together. The three considered compression factors of approximately $0.637$, $0.827$ and $0.900$ correspond to circular arcs with central angles $\pi/2$, $\pi/3$ and $\pi/4$.

  \begin{figure}
  \centering
  \begin{subfigure}{0.49\textwidth}
    \centering
    \includegraphics[width=\textwidth]{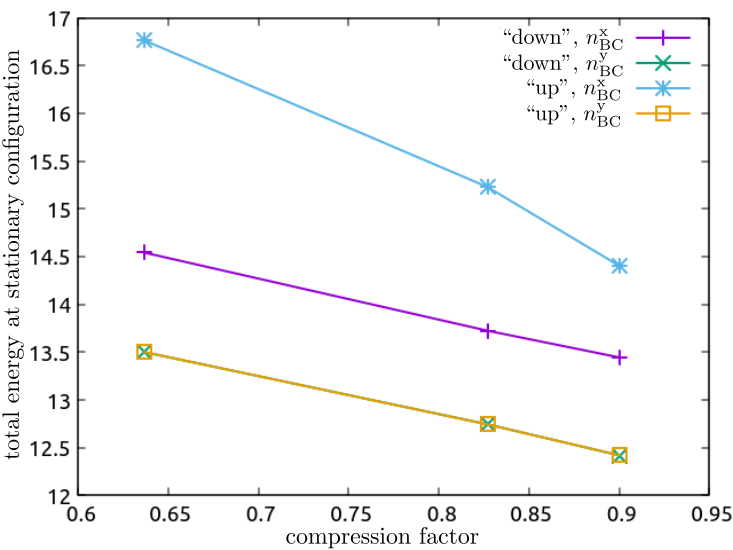}
    \caption{$\bar{r} = 4$, $\bar{\e} = 1$.}
  \end{subfigure}
  \hfill
  \begin{subfigure}{0.49\textwidth}
    \centering
    \includegraphics[width=\textwidth]{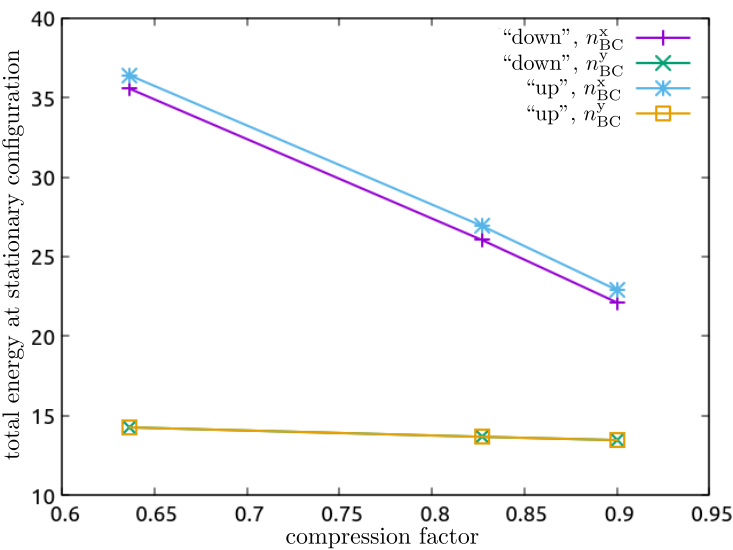}
    \caption{$\bar{r} = 4$, $\bar{\e} = 4$.}
  \end{subfigure}
  \vspace{-2mm}
  \caption{Plots of final energies versus three different compression factors in Example~\ref{ex:compress}. In the case $\bar{r} = 4$, $\bar{\varepsilon} = 4$ strong Dirichlet anchorings may be used to select between almost no energy response (for $n^\mathrm{y}_\mathrm{BC}$) and significant energy response (for $n^\mathrm{x}_\mathrm{BC}$) with respect to varying compression factors.}
  \label{fig:ex2_2}
  \end{figure}

    \begin{table}
    \centering
      \begin{tabular}{ c c c c c c c c c }
        $\alpha$ & $\bar{\e}$ & anchoring & buckling & $\mathcal{E}_\mathfrak{h}$
        & $\mathcal{E}_{\mathfrak{h},\mathrm{OF}}$
        & $\mathrm{Err}_\mathrm{1}(n_\mathfrak{h}^\infty)$
        & $\mathrm{Err}_\mathrm{iso}(y_\mathfrak{h}^\infty)$\\\hline\hline
        $0.637$
        &1&$n^\mathrm{x}_\mathrm{BC}$& ``up''
          & \num{16.7652} & \num{2.25081} & \num{0.000644441} & \num{0.00108599} \\
        &&$n^\mathrm{x}_\mathrm{BC}$&``down''
          & \num{14.5429} & \num{2.43425} & \num{0.0110075} & \num{0.00103477} \\
        &&$n^\mathrm{y}_\mathrm{BC}$& ``up''
          & \num{13.5005} & \num{0.611622} & \num{6.92803e-07} & \num{0.000606317} \\
        &&$n^\mathrm{y}_\mathrm{BC}$&``down''
          & \num{13.5012} & \num{0.611275} & \num{5.77772e-07} & \num{0.00060703} \\\hline
        &4&$n^\mathrm{x}_\mathrm{BC}$& ``up''
          & \num{36.388} & \num{15.4067} & \num{0.0458814} & \num{0.0187701} \\
        &&$n^\mathrm{x}_\mathrm{BC}$&``down''
          & \num{35.549} & \num{15.765} & \num{0.0565892} & \num{0.0185882} \\
        &&$n^\mathrm{y}_\mathrm{BC}$& ``up''
          & \num{14.2776} & \num{0.0907788} & \num{3.64561e-08} & \num{0.000245617} \\
        &&$n^\mathrm{y}_\mathrm{BC}$&``down''
          & \num{14.2776} & \num{0.0907119} & \num{3.65312e-08} & \num{0.000245793} \\\hline\hline
        $0.827$
        &1&$n^\mathrm{x}_\mathrm{BC}$& ``up''
          & \num{15.2285} & \num{1.42447}  & \num{2.55569e-11} & \num{0.000381252} \\
        &&$n^\mathrm{x}_\mathrm{BC}$&``down''
          & \num{13.7253} & \num{1.73153}  & \num{0.00685615} & \num{0.000529905} \\
        &&$n^\mathrm{y}_\mathrm{BC}$& ``up''
          & \num{12.7444} & \num{0.750209} & \num{1.01807e-06} & \num{0.000570308} \\
        &&$n^\mathrm{y}_\mathrm{BC}$&``down''
          & \num{12.7446} & \num{0.750596} & \num{6.49192e-07} & \num{0.000573242} \\\hline
        &4&$n^\mathrm{x}_\mathrm{BC}$& ``up''
          & \num{26.931}  & \num{9.17482}  & \num{0.012309} & \num{0.00783711} \\
        &&$n^\mathrm{x}_\mathrm{BC}$&``down''
          & \num{26.0491} & \num{9.42931}  & \num{0.0188844} & \num{0.00763891} \\
        &&$n^\mathrm{y}_\mathrm{BC}$& ``up''
          & \num{13.6858} & \num{0.1057}   & \num{4.25356e-08} & \num{0.000108281} \\
        &&$n^\mathrm{y}_\mathrm{BC}$&``down''
          & \num{13.6858} & \num{0.105641} & \num{4.27544e-08} & \num{0.000108279} \\\hline\hline
        $0.900$
        &1&$n^\mathrm{x}_\mathrm{BC}$& ``up''
          & \num{14.4009} & \num{0.803111} & \num{1.60889e-11} & \num{0.000208727} \\
        &&$n^\mathrm{x}_\mathrm{BC}$&``down''
          & \num{13.4439} & \num{1.46523} & \num{0.00442234} & \num{0.000385065} \\
        &&$n^\mathrm{y}_\mathrm{BC}$& ``up''
          & \num{12.4196} & \num{0.829226} & \num{1.30263e-06} & \num{0.000625317} \\
        &&$n^\mathrm{y}_\mathrm{BC}$&``down''
          & \num{12.4193} & \num{0.829192} & \num{5.78022e-07} & \num{0.000630694} \\\hline
        &4&$n^\mathrm{x}_\mathrm{BC}$& ``up''
          & \num{22.8881} & \num{6.65747} & \num{0.00173701} & \num{0.00335856} \\
        &&$n^\mathrm{x}_\mathrm{BC}$&``down''
          & \num{22.1056} & \num{6.83343} & \num{0.00473126} & \num{0.00378197} \\
        &&$n^\mathrm{y}_\mathrm{BC}$& ``up''
          & \num{13.4726} & \num{0.114663} & \num{4.2816e-08} & \num{7.23188e-05} \\
        &&$n^\mathrm{y}_\mathrm{BC}$&``down''
          & \num{13.4727} & \num{0.114627} & \num{4.33516e-08} & \num{7.22175e-05} \\
      \end{tabular}
    \caption{Experimental results for the compressed plate in Example~\ref{ex:compress} with $\bar{r}=1$ and strong tangential anchoring with different compression factors $\alpha$. Varying proportions $\mathcal{E}_\mathfrak{h}/\mathcal{E}_{\mathfrak{h},\mathrm{OF}}$ are clearly recognizable.}
    \label{tab:ex2}
  \end{table}
\end{example}

\begin{example}[Strong circular director anchoring]\label{ex:circular}
  We let $\mu = 1$, $\lambda = 1000$ and use parameters $\bar{r} = 4$, $\bar{\e} = 0.5$.
  For cut-out regions
  \[
  K_1 \colonequals  \emptyset,\;\;\; K_2 \colonequals  [1,4]^2,\;\;\; K_3 \colonequals  [2,3]^2,\;\;\; K_4 \colonequals  [2,3] \times [1,4],
  \]
  we consider the two-dimensional reference configurations $S_j \colonequals  (0,5)^2 \setminus K_j$
  with Dirichlet boundary $\Gamma_y = \emptyset$, $\Gamma_n = \{0,5\} \times [0,5] \cup [0,5] \times \{ 0,5\}$.
  There are no clamped boundary conditions in this example. Instead, we fix a reference frame by prescribing values for $y$ and $\nabla y$ at only one vertex, thereby guaranteeing well-definition of the descent steps for $y$ in the discrete gradient flow.
  With the center of the domain $m_\Omega = (2.5,2.5)^\top$, we use strong Dirichlet anchoring
  \[
  \hat{n}_\mathrm{BC}(x) = |x-m_\Omega|^{-1}(-x_2+2.5,\, x_1-2.5,\, 0)^\top
  \]
  for the director $\hat{n} = R_y^\top n$ in local coordinates.

  For triangulations $\mathcal{T}^j_\mathfrak{h} = \mathcal{T}^j_k$, $k=1,2,3,4$, of $S_j$ into halved squares with side lengths $\hat{\mathfrak{h}} = 2^{-k}$ we compute approximations of stationary points of the energy~\eqref{eq:energy-discrete1} using the alternating discrete gradient flow Algorithm~\ref{alg:gradflow} with stopping criterion $\e_\mathrm{stop} = 2 \times 10^{-3}$.
  As initial values we use $y_\mathfrak{h}^0(x) = (x,0)^\top$ for the deformation and a smooth extension of $\hat{n}_\mathrm{BC}$ to $\Omega \setminus B_{\hat{\mathfrak{h}}}(0)$ for the director, and set $n_\mathfrak{h}^0 (0) = (0,0,1)^\top$ for the domain $S_1$.
  The finest triangulations $\mathcal{T}^j_4$ consist of at most 12'800 triangles and corresponding discrete deformations and director fields are determined by 59'049 and 19'683 degrees of freedom, respectively.
  For each triangulation we choose the (pseudo-) time step size $\tau =\hat{\mathfrak{h}}/2$.
  Resulting deformations and director field configurations together with their final energies are illustrated in Figure~\ref{fig:ex3}.
  The resulting iteration numbers, final energies, constraint violations and experimental convergence rates are listed in Table~\ref{tab:ex3}.
  The largest deformation is developed for the largest cut-out region $K_2$.
  In this setting a significantly higher number of iterations was needed in the discrete gradient flow to reach the approximate equilibrium state than in the other settings.
  As in Example~\ref{ex:clampedplate} we observe that the violations of both the nodal unit-length and isometry constraint seem to be bounded linearly in terms of the step size.

  \begin{figure}
  \centering
  \begin{subfigure}{0.45\textwidth}
    \centering
    \includegraphics[width=\textwidth]{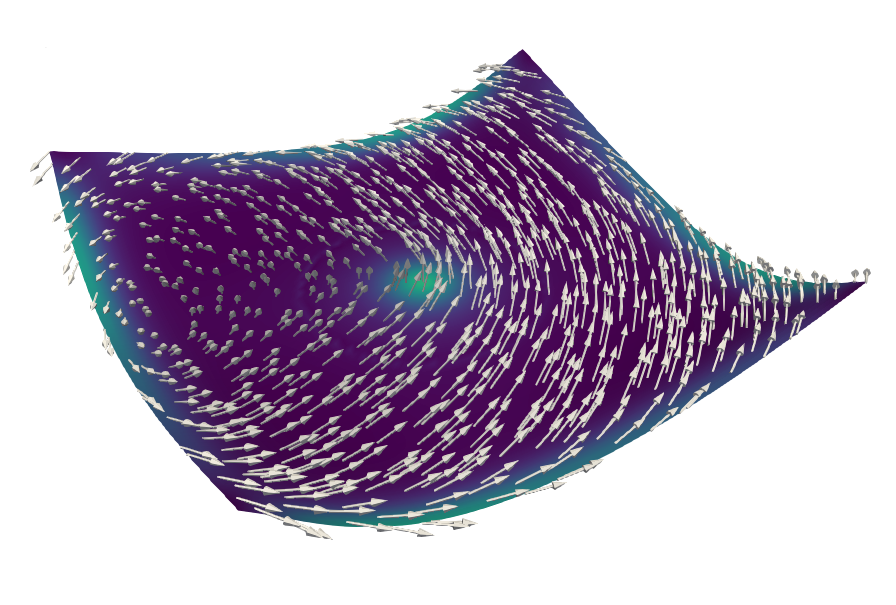}
    \put(-201,128){$\mathcal{E}_\mathfrak{h}=9.153$}
    \put(-201,116){$\mathcal{E}_{\mathfrak{h},\mathrm{OF}}=3.157$}
    \caption{Tangential circular director anchoring on square domain.}
  \end{subfigure}
  \hspace{5mm}
  \begin{subfigure}{0.45\textwidth}
    \centering
    \includegraphics[width=\textwidth]{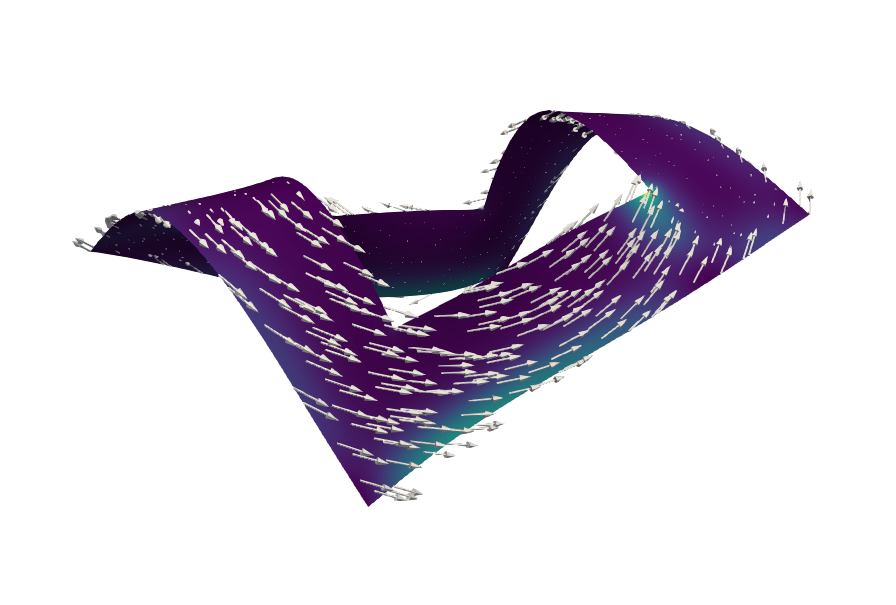}
    \put(-201,128){$\mathcal{E}_\mathfrak{h}=6.726$}
    \put(-201,116){$\mathcal{E}_{\mathfrak{h},\mathrm{OF}}=2.012$}
    \caption{Tangential circular director anchoring on square domain with square cutout.}
  \end{subfigure}
  \\[5mm]
  \begin{subfigure}{0.45\textwidth}
    \centering
    \includegraphics[width=\textwidth]{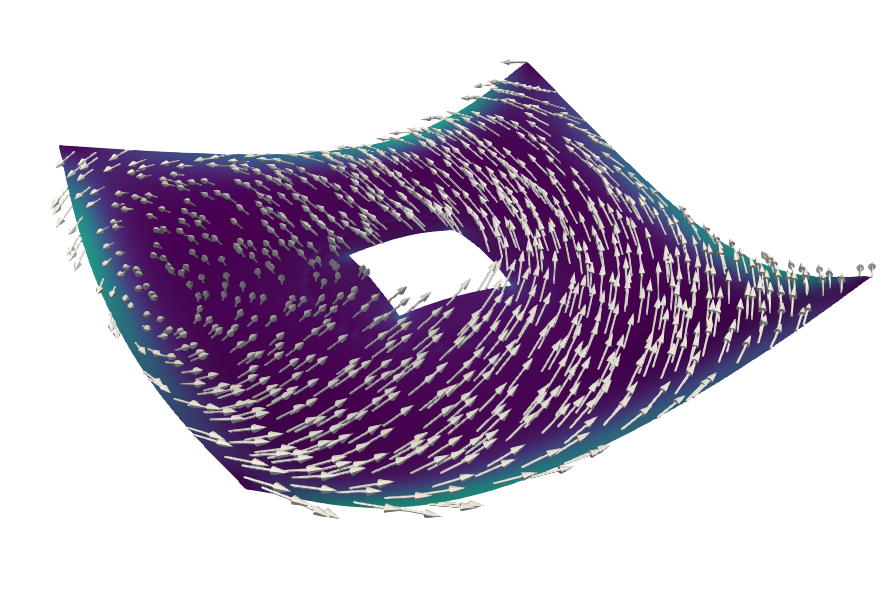}
    \put(-201,128){$\mathcal{E}_\mathfrak{h}=8.221$}
    \put(-201,116){$\mathcal{E}_{\mathfrak{h},\mathrm{OF}}=2.788$}
    \caption{Tangential circular director anchoring on square domain with small square cutout.}
  \end{subfigure}
  \hspace{5mm}
  \begin{subfigure}{0.45\textwidth}
    \centering
    \includegraphics[width=\textwidth]{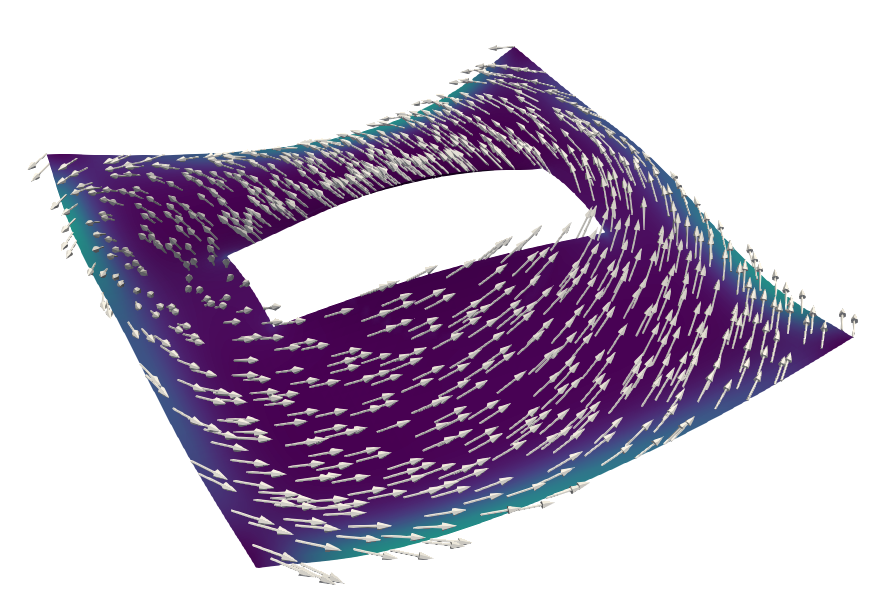}
    \put(-201,128){$\mathcal{E}_\mathfrak{h}=7.437$}
    \put(-201,116){$\mathcal{E}_{\mathfrak{h},\mathrm{OF}}=2.506$}
    \caption{Tangential circular director anchoring on square domain with rectangular cutout.}
  \end{subfigure}
  \vspace{-2mm}
  \caption{Deformations and corresponding final energies for strong circular tangential director anchoring with different cutout regions in Example~\ref{ex:circular}.}
  \label{fig:ex3}
  \end{figure}
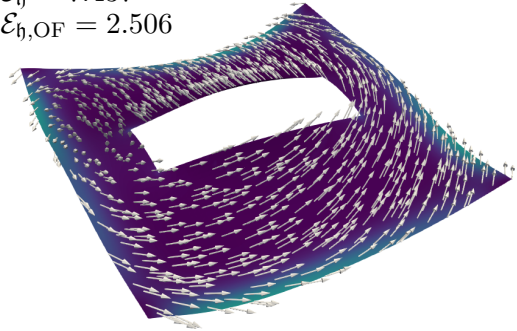

  \begin{table}
    \centering
      \begin{tabular}{ c c r r c c c c c }
        Cutout & $k$ & $\#\mathcal{T}_\mathfrak{h}^j$ & \#iter & $\mathcal{E}_\mathfrak{h}$
        & $\mathcal{E}_{\mathfrak{h},\mathrm{OF}}$
        & $\mathrm{Err}_\mathrm{1}(n_\mathfrak{h}^\infty)$
        & $\mathrm{Err}_\mathrm{iso}(y_\mathfrak{h}^\infty)$
        & $\mathrm{EOC}_\mathfrak{h}^\mathcal{E}$ \\\hline\hline
        $K_1$
        &  1  &    200  &   2516  &  \num{7.35008}  &  \num{4.27803}  &  \num{0.0312128}   &  \num{0.116302}   &  ---  \\
        &  2  &    800  &   1645  &  \num{8.5795}   &  \num{3.24003}  &  \num{0.0149404}   &  \num{0.0431381}  &  ---  \\
        &  3  &   3200  &   2205  &  \num{8.98881}  &  \num{3.17945}  &  \num{0.00767582}  &  \num{0.0147836}  & \num{1.58671}  \\
        &  4  &  12800  &   3648  &  \num{9.1528}   &  \num{3.15671}  &  \num{0.00383612}  &  \num{0.0052758}  & \num{1.31959}  \\\hline
        $K_2$
        &  1  &    128  &   2850  &  \num{4.72833}  &  \num{3.84747}  &  \num{0.00935377}  &  \num{0.323404}   &  ---  \\
        &  2  &    512  &   8983  &  \num{6.41888}  &  \num{2.104}    &  \num{0.0054644}   &  \num{0.111732}   &  ---  \\
        &  3  &   2048  &  20227  &  \num{6.76598}  &  \num{2.08913}  &  \num{0.00341109}  &  \num{0.0355602}  & \num{2.28407}  \\
        &  4  &   8192  &  24782  &  \num{6.72345}  &  \num{2.01209}  &  \num{0.0015231}   &  \num{0.0120369}  & \num{3.0288}  \\\hline
        $K_3$
        &  1  &    192  &   2621  &  \num{6.47214}  &  \num{4.32488}  &  \num{0.0298403}   &  \num{0.143864}   &  ---  \\
        &  2  &    768  &   1767  &  \num{7.68248}  &  \num{2.88389}  &  \num{0.0151286}   &  \num{0.0503269}  &  ---  \\
        &  3  &   3072  &   2470  &  \num{8.04154}  &  \num{2.80975}  &  \num{0.0078652}   &  \num{0.0165502}  & \num{1.75312}  \\
        &  4  &  12288  &   4122  &  \num{8.22126}  &  \num{2.7875}   &  \num{0.00394469}  &  \num{0.0057124}  & \num{0.998474}  \\\hline
        $K_4$
        &  1  &    176  &   2766  &  \num{5.45137}  &  \num{4.31333}  &  \num{0.0286767}   &  \num{0.212539}   &  ---  \\
        &  2  &    704  &   1562  &  \num{6.97267}  &  \num{2.52291}  &  \num{0.0150215}   &  \num{0.0773614}  &  ---  \\
        &  3  &   2816  &   2382  &  \num{7.2839}   &  \num{2.50053}  &  \num{0.00824715}  &  \num{0.0270039}  & \num{2.28925}  \\
        &  4  &  11264  &   4072  &  \num{7.43682}  &  \num{2.50587}  &  \num{0.00403532}  &  \num{0.0145889}  & \num{1.0252}  \\
      \end{tabular}
    \caption{Experimental results for strong circular tangential director anchoring for different refinement levels $k$ and different cutout regions in Example~\ref{ex:circular}.}
    \label{tab:ex3}
  \end{table}
\end{example}

\subsection{Conclusion}\label{sec:conclusion}
In the simulations the delicate interplay of competing energy contributions in the LCE bilayer model results in rich, non-trivial deformations.
Our experiments confirm that the bending and buckling behavior of an LCE plate can be influenced, or even controlled, in a multitude of ways.
Depending on model parameters, anchorings and boundary conditions, either the Oseen-Frank term that corresponds to crystalline order, or the elastic and residual stress terms that correspond to coupling effects, play a dominant role in the energy minimization.
Interesting effects can occur when transitioning from dominant crystalline order to dominant coupling in the energy landscape or vice versa, cf.~the case $\bar{r}=5$, $\bar{\e}=1$ in Figure~\ref{fig:ex1}.
Bearing in mind possible engineering applications, values of $\bar\e \sim1$ correspond to plates with a physical thickness in the range of $1\si{nm}$ -- $100\si{nm}$,~cf.~Section~\ref{S:modeling}.
In theory one could compensate for a greater plate thickness by employing soft materials with very low shear modulus values.
High Oseen-Frank energies can, of course, also result from certain anchoring conditions.
The softening and hardening of a strip with respect to compression in Example~\ref{ex:compress} is closely related to such effects.
In that example the buckling behavior after compression is influenced by the respective anchorings, that either allow for almost perfect director alignment, i.\,e., practically negligible Oseen-Frank energies, or yield Oseen-Frank energies that are equal to approximately half of the total energy.
In the latter case, the effect that comes into play is an incompatibility of the strong anchoring conditions on the Dirichlet boundary and mid-surface with the enforced buckling deformation that results from compressive boundary conditions.

\section{Proofs}\label{S:proofs}
In Section~\ref{Sec:Proof:Compactness} we prove Theorem \ref{T1}~\ref{item:T1:compactness} and Theorem~\ref{T2}~\ref{item:T2:Comp}, which are concerned about compactness.
In Section~\ref{Section:Representation} we prove Lem\-mas~\ref{L:representations:hom:iso}, \ref{L:properties} and \ref{L:reduction}.
These results establish properties of $Q_{\mathrm{el}}$, $E_{\mathrm{res}}$ and $\mathbb B$.
Specifically, Lemma~\ref{L:reduction} yields a relaxation formula, which is used in the proofs of the lower- and upper-bound parts of Theorems~\ref{T1} and \ref{T2}.
These proofs are presented in Sections~\ref{Sec:Proof:Lower:Bound} and \ref{Sec:Proof:Upper:Bound}, respectively.
The proof of Lemma~\ref{L:cont} is presented in Section~\ref{Sec:Proof:Properties:of:Limit}.
In Section \ref{Sec:Proof:Anchoring}, we present the proofs for Lemmas~\ref{P:Gamma:Strong:Anchoring} -- \ref{L:weak:anchoring:with:boundary:slices} about anchorings for the director field.

In the proofs of the upper bound parts of Theorems~\ref{T1} and \ref{T2} (cf. Section~\ref{Sec:Proof:Upper:Bound}) we appeal to two auxiliary results of independent interest: the density result, cf.~Proposition~\ref{P:approx}, and a general construction for 3d-deformations, cf.~Proposition~\ref{P:general_construction}. The proofs of these results are presented in Section~\ref{Sec:Proof:P:approx} and Section~\ref{Sec:Proof:P:general_construction}, respectively.
 
\subsection{Compactness: proof of Theorem~\ref{T1}~\ref{item:T1:compactness} and Theorem~\ref{T2}~\ref{item:T2:Comp}}
\label{Sec:Proof:Compactness}

For the proof of the compactness and lower bound parts of Theorems~\ref{T1} and \ref{T2}, we require the following a priori estimates:
\begin{lemma}[A priori estimates]\label{L:pf:apriori}
  Let Assumption \ref{ass:W} be satisfied and let $S$ be bounded and open.
  Consider a sequence $(y_h,n_h)$ with finite energy in the sense of \eqref{EqFiniteEnergy}. Let $1<p<2$ be defined by the identity $q_W=4\frac{p}{2-p}$. Then,
  \begin{align}
    \label{eq:st:2d}
    &\limsup\limits_{h\to 0}\frac1{h^2}\int_{\Omega}\dist^2(\nabla\!_hy_h,\SO 3)<\infty,\\
    \label{eq:st:2c}
    &\limsup\limits_{h\to 0}\int_{\Omega}|\nabla\!_hy_h|^{q_W}+|\det(\nabla\!_hy_h)|^{-\frac{q_W}{2}}<\infty,\\
    \label{eq:st:2b}
    &\limsup\limits_{h\to0}\int_{\OmegaT}|\nabla\!_hn_h(\nabla\!_hy_h)^{-1}|^2|\det(\nabla\!_h y_h)|<\infty,\\\label{eq:st:2a}
    &\limsup\limits_{h\to0}\int_{\OmegaT}|\nabla\!_hn_h|^p<\infty.
  \end{align}
  Furthermore, there exists $\bar h>0$ (depending on the sequence $(y_h,n_h)$) such that
  \begin{equation}\label{eq:st:0}
    \det(\nabla\!_hy_h)>0\text{ a.e.~in }\Omega\qquad\text{for all }0<h\leq\bar h.
  \end{equation}
\end{lemma}
\begin{proof}
  For convenience we introduce the shorthand notation
  \begin{equation}\label{eq:st:9:Bh}
    B_h(x)\colonequals 
    \begin{cases}
      B_h(n_h(x))&\text{if }x\in\OmegaT,\\
      0&\text{else,}
    \end{cases}
  \end{equation}
  and note that there exists a constant $C=C(\bar r)>0$ such that for all $0<h\leq \frac{1}{C}$,
  \begin{equation}\label{eq:st:1}
    |B_h|\leq C\qquad\text{and} \qquad |\det(I+hB_h)-1|+|(I+hB_h)^{-1}-I|\leq \min\{\tfrac12,C h\}.
  \end{equation}
  In view of \eqref{EqFiniteEnergy} we have $\mathcal E_h(y_h,n_h)<\infty$ for all $0<h\leq\bar h$ and some $\bar h>0$.
  In the rest of the proof we assume that $C\bar h|\bar r|\leq 1$ and $0<h\leq\bar h$.
  
  \step{1 -- Proof of \eqref{eq:st:2d}}\smallskip
  
  From \eqref{EqFiniteEnergy}, \eqref{eq:Bh}, and \ref{item:ass:W:non:degen:nat}, we conclude that
  \begin{equation*}
    \limsup\limits_{h\to 0}\frac{1}{h^2}\int_\Omega\dist^2((I+hB_h)\nabla\!_hy_h,\SO 3)<\infty,
  \end{equation*}
  and thus, \eqref{eq:st:2d} follows with help of \eqref{eq:st:1} and triangle's inequality.

  \step{2 -- Proof of \eqref{eq:st:2c} -- \eqref{eq:st:0}}\smallskip
  
  We start with the argument for \eqref{eq:st:0} by noting that $\det(\nabla\!_hy_h)=\det((I+hB_h)^{-1})\det\big((I+hB_h)\nabla\!_hy_h\big)$.
  The first factor on the right-hand side is positive by \eqref{eq:st:1}, while the second factor is positive thanks to $\mathcal E_h(y_h,n_h)<\infty$ and the growth condition \ref{item:ass:W:growth:cond}.
  Hence, \eqref{eq:st:0} follows.
  In view of the latter, \eqref{eq:st:2b} directly follows from \eqref{EqFiniteEnergy}.
  To see \eqref{eq:st:2c}, we note that \eqref{eq:st:1} yields  $|\nabla\!_h y_h|\leq 3|(I+hB_h)\nabla\!_hy_h|$ and $|\det(\nabla\!_hy_h)|^{-1}\leq 2|\det\big((I+hB_h)\nabla\!_hy_h\big)|^{-1}$,
which together with \ref{item:ass:W:growth:cond} implies
\begin{equation*}
  \int_{\OmegaT}|\nabla\!_hy_h|^{q_W}+|\det(\nabla\!_hy_h)|^{-\frac{q_W}{2}}
  \leq C(q_W,C_W)\int_{\OmegaT}\Big( W(x_3,(I+hB_h)\nabla\!_hy_h)+1\Big),
\end{equation*}
where $C(q_W,C_W)>0$ only depends on $q_W$ and $C_W$.
Hence, \eqref{eq:st:2c} follows, since the right-hand side is controlled by $1+h^2\mathcal E(y_h,n_h)$.

Finally, from $|\nabla\!_h n_h|\leq \det(\nabla\!_hy_h)^{\frac12}|\nabla\!_hn_h(\nabla\!_hy_h)^{-1}|\,|\nabla\!_hy_h|\,\det(\nabla\!_hy_h)^{-\frac12}$,
and a triple H\"older's inequality, we get
\begin{eqnarray*}
  \left(\int_{\OmegaT}|\nabla\!_hn_h|^p\right)^\frac1p&\leq&  \left(\int_{\OmegaT}|\nabla\!_hn_h(\nabla\!_h y_h)^{-1}|^2|\det(\nabla\!_h y_h)|\right)^\frac12\\
  &&\times\left(\int_{\OmegaT}|\nabla\!_hy_h|^{q_W}\right)^{\frac1{q_W}}\left(\int_{\OmegaT}|\det(\nabla\!_hy_h)|^{-\frac{q_W}{2}}\right)^{\frac1{q_W}}.
\end{eqnarray*}
Hence, \eqref{eq:st:2a} follow from \eqref{eq:st:2c} and \eqref{eq:st:2b}.
\end{proof}
\medskip

Next, we recall the celebrated compactness result for sequences with finite bending energies:
\begin{lemma}[\mbox{\cite[Theorem~4.1]{FJM02}}]\label{L:compfinitebend}
Let $S\subseteq \mathbb R^2$ be a bounded Lipschitz domain.
  Let $(y_h)$ be a sequence in $H^1(\Omega;\R^3)$ with finite bending energy in the sense of \eqref{eq:st:2d}. 
Then there exist $y\in H^2_{\iso}(S;\R^3)$, $M\in L^2(S;\R^{2\times 2}_{\sym})$, and $d\in L^2(\Omega;\R^3)$, such that for a subsequence (not relabeled) we have
  \begin{align}\label{eq:st:3a}
    \nabla\!_hy_h\to\,&\,R_y\qquad&&\text{strongly in }L^2(\Omega;\R^{3\times 3}),\\\label{eq:st:3b}
    E_h(y_h)\wto\,&\,\iota(x_3\II_y+M)+\sym(d\otimes e_3)\qquad&&\text{weakly in }L^2(\Omega;\R^{3\times 3}_{\sym}),
  \end{align}
  where we recall the definition $E_h(y_h)\colonequals \frac{\sqrt{(\nabla\!_hy_h)^\top\nabla\!_h y_h}-I}{h}$ of the nonlinear strain.
\end{lemma}
We note that \eqref{eq:st:3b} is not explicitly stated in \cite{FJM02}; yet, it can be established along the lines of the proof of \cite[Theorem 7.1]{FJM02}.
Statement \eqref{eq:st:3b} is also a direct consequence of the two-scale compactness result for the nonlinear strain \cite[Proposition~3.2]{B_HNB14}.

\begin{proof}[Proof of Theorem~\ref{T1}~\ref{item:T1:compactness}]
  \step{1 -- Argument for \eqref{KompaktheitKonvergenzDef}, \eqref{KompaktheitKonvergenzDir}, \eqref{KompaktheitKonvergenzDir2} and $y\in H^2_{\iso}(S;\R^3)$}
  \smallskip
  
  Statement \eqref{KompaktheitKonvergenzDef} and $y\in H^2_{\iso}(S;\R^3)$ directly follow from Lemma~\ref{L:compfinitebend}, which we may apply thanks to \eqref{eq:st:2d}.
  Moreover,  thanks to \eqref{eq:st:2a} and the unit-length constraint, $(n_h)$ is bounded in $W^{1,p}(\OmegaT;\R^3)$.
  Thus, there exists $n\in W^{1,p}(\OmegaT;\R^3)$ such that $n_h\wto n$ weakly in $W^{1,p}(\OmegaT)$ for a subsequence.
  By compact embedding $W^{1,p}(\OmegaT;\R^3)\subseteq L^1(\OmegaT;\R^3)$ we have $n_h\to n$ a.e.~(for a subsequence), and thus $n\in\mathbb S^2$ a.e.~in $\OmegaT$.
  Dominated convergence yields \eqref{KompaktheitKonvergenzDir}.
  By \eqref{eq:st:2a} we obtain \eqref{KompaktheitKonvergenzDir2} (after possibly passing to a further subsequence).
  Moreover, from $\frac{1}{h}\partial_3n_h\wto d$ weakly in $L^p(\OmegaT)$, we deduce that $\partial_3n_h\to 0$ strongly in $L^p(\OmegaT)$.
  Hence, $\partial_3n=0$ and we can identify $n$ with a function $n\in W^{1,p}(S;\R^3)$.

  \step{2 -- Proof of $(\nabla'n,d)\in L^2(\OmegaT;\R^{3\times 3})$}\smallskip

  The claimed integrability is a consequence of the following observation: Consider sequences $(G_h)$, $(F_h)$ such that $G_h\wto G$ weakly in $L^p(\OmegaT;\R^{3\times 3})$, and $F_h\to R$ strongly in $L^2(\OmegaT;\R^{3\times 3})$.
  Suppose that a.e.~in $\Omega$ we have $\det F_h>0$ and $R\in\SO 3$. Then,
  \begin{equation}\label{eq:st:20}
    \int_{\OmegaT}|G|^2\leq \liminf\limits_{h\to 0}\int_{\OmegaT}|G_hF_h^{-1}|^2\det(F_h)\,\mathrm{d}x.
  \end{equation}
  Before we present the argument, we note that an application of \eqref{eq:st:20} with $G_h\colonequals \nabla\!_hn_h$, $G=(\nabla'n,d)$, $F_h\colonequals \nabla\!_hy_h$, $R=R_y$, combined with \eqref{eq:st:2b} (to bound the right-hand side of \eqref{eq:st:20}) yields $(\nabla'n,d)\in L^2(\OmegaT;\R^{3\times 3})$.
  For the proof of \eqref{eq:st:20} consider the map
  \begin{equation*}
    \Phi:\R^{3\times 3}\to \R^{3\times 3},\qquad \Phi(F)\colonequals 
    \begin{cases}
      F^{-1}\sqrt{\det F}&\text{if }\dist(F,\SO 3)\leq\frac12,\\
      0&\text{else,}
    \end{cases}
  \end{equation*}
  and note that $\Phi$ is measurable and bounded, i.e., $C_\Phi\colonequals \sup_{F\in \R^{3\times 3}}|\Phi(F)|<\infty$.
  We may pass to a subsequence such that $F_h\to R$ a.e.~in $\OmegaT$. Since $\Phi$ is continuous in an open neighborhood of $\SO 3$, and since $R\in\SO 3$ a.e., we conclude that $\Phi(F_h)\to \Phi(R)=R^\top$ a.e. 
  Together with the bound $|\Phi(F_h)|\leq C_\Phi$, we conclude that $G_h\Phi(F_h)\wto GR^\top$ weakly in $L^p(\OmegaT)$,
  and thus, by the weak lower semicontinuity of convex integral functionals, we have
  \begin{equation*}
    \int_{\OmegaT}|GR^\top|^2\leq \liminf\limits_{h\to 0}\int_{\OmegaT}|G_h\Phi(F_h)|^2.
  \end{equation*}
  Now, \eqref{eq:st:20} follows, by combining this with the pointwise identity $|GR^\top|=|G|$ and the pointwise upper bound $|G_h\Phi(F_h)|^2\leq |G_h(F_h)^{-1}|^2|\det(F_h)|$ (which is a consequence of the definition of $\Phi$).
\end{proof}

\begin{proof}[Proof of Theorem~\ref{T2}~\ref{item:T2:Comp}]
  Theorem~\ref{T1}~\ref{item:T1:compactness} already yields \eqref{KompaktheitKonvergenzDir}, \eqref{KompaktheitKonvergenzDir2} for some $n\in H^1(S,\mathbb S^2)$.
  
  \step{1 -- Strong convergence of $(y_h)$ in $L^2$}\smallskip
  
  By \eqref{KompaktheitKonvergenzDef} there exists $y\in H^2_{\mathrm{iso}}(S;\R^3)$ such that
  \begin{equation}\label{eq:st:5}
   \nabla\!_h y_h\to R_y\qquad\text{strongly in }L^2(\Omega),
  \end{equation}
  up to a subsequence. With help of the boundary condition \eqref{eq:BC3d} and Poincare's inequality, we conclude that $(y_h)$ is bounded in $H^1(\Omega)$.
  Hence, the compact embedding $H^1\subseteq L^2$ yields in addition to \eqref{eq:st:5} that $y_h\to y+c$ strongly in $L^2(\Omega)$ for some constant vector $c\in\R^3$.
  Since the right-hand side in \eqref{eq:st:5} is invariant w.r.t.~addition of constants, we may assume w.l.o.g. that $c=0$.

  \step{2 -- Argument for $y\in \mathcal A_{BC}$}\smallskip

We consider the line segment $\mathcal L_i$ and show that $y$ and $\nabla' y$ satisfy the required boundary condition on $\mathcal L_i$.  
  Let $x_0'\in\mathcal L_i$.
   Since $\mathcal L_i$ is relatively open, we can choose  $r_0>0$ such that the ball $B(x_0';r_0)\subseteq\R^2$ satisfies
  \begin{align*}
    B(x_0';r_0)\cap\partial S\subseteq\mathcal L_i\qquad\text{and}\qquad B(x_0';r_0)\cap S\text{ is connected}.
  \end{align*}
  Next, we extend $y_h$ affinely to the  3d-domain $B\colonequals B(x_0';r_0)\times(-\frac12,\frac12)$ by the following procedure:
  Since $y_{BC}$ is affine on $\mathcal L_i$, there exist $R_0\in\SO 3$ and $y_0\in\R^3$, such that $y_{BC}(x')=y_0+R_0(x'-x_0',0)^\top$ on $\mathcal L_i$.
  Consider the extension
  \begin{equation*}
    \tilde y_h:B\to\R^3,\qquad \tilde y_h(x',x_3)\colonequals 
    \begin{cases}
      (1-h\bar\delta)(y_0+R_0(x'-x_0',0)^\top)+hx_3R_0e_3&\text{if }x'\in B(x_0';r_0)\setminus S,\\
      y_h(x',x_3)&\text{else.}
    \end{cases}
  \end{equation*}
  By construction, we have $\tilde y_h\vert_D\in H^1(D;\R^3)$ for $D=B\setminus\overline \Omega$ and $D=B\cap \Omega$, respectively.
  Furthermore, the traces of these two restrictions are equal thanks to \eqref{eq:BC3d} and the assumed linearity of $y_{BC}$ on $\mathcal L_i$.
  We conclude that $\tilde y_h\in H^1(B;\R^3)$. Moreover, from
  \begin{equation}\label{eq:st:6}
    \nabla\!_h\tilde y_h=
    \begin{cases}
      R_0-h\delta(R_0e_1,R_0e_2,0)&\text{on }B\setminus \Omega,\\
      \nabla\!_hy_h&\text{on }B\cap \Omega,
    \end{cases}
  \end{equation}
  we conclude that $(\tilde y_h)$ has finite bending energy in the sense that $\limsup\limits_{h\to 0}h^{-2}\int_B\dist^2(\nabla\!_h\tilde y_h,\SO 3)<\infty$. Hence, by Lemma~\ref{L:compfinitebend} and the argument of Step~1 there exists $\tilde y\in H^2_{\iso}(B(x_0';r_0);\R^3)$ such that $(\tilde y_h,\nabla\!_h\tilde y_h)\to (\tilde y,R_{\tilde y})$ strongly in $L^2(B)$. By passing to the limit in \eqref{eq:st:6} with help of \eqref{eq:st:5} and $y_h\to y$ in $L^2(\Omega)$, we obtain the identity
  \begin{equation*}
    (\tilde y,R_{\tilde y})=
    \begin{cases}
      (y_0+R_0\,(\cdot-x_0',0)^\top ,R_0)&\text{in }B(x_0';r_0)\setminus \overline S\\
      (y,R_y)&\text{in }B(x_0';r_0)\cap S.
    \end{cases}
  \end{equation*}
  Since $(\tilde y,R_{\tilde y})\in H^1(B(x_0';r_0))$ and $(y,R_y)\in H^1(B(x_0';r_0)\cap S)$, we deduce
  \begin{equation*}
    (y,R_y)=(y_0+R_0(\cdot-x_0',0)^\top,R_0)=(y_{BC},R_{y_{BC}})\qquad\text{on }\mathcal L_i\cap B(x_0';r_0)
  \end{equation*}
  in the sense of traces.
  Since $x_0'\in\mathcal L_i$ is arbitrary, we conclude that $y$ satisfies the required boundary condition on $\mathcal L_i$.
  Since $i\in\{1,\ldots,k_{BC}\}$ is arbitrary, $y\in\mathcal A_{BC}$ follows.
\end{proof}

\subsection{Representation: a relaxation formula and proof of Lemmas~\ref{L:representations:hom:iso} and \ref{L:properties}}\label{Section:Representation}
We first state and prove the following relaxation result.
It motivates Definition~\ref{def:QEB} and it is used in the proofs of Theorem~\ref{T1} \ref{item:T1:lower:bound} and \ref{item:T1:upper:bound}.
\begin{lemma}[Relaxation formula]\label{L:reduction}
  For all $A\in\R^{2\times 2}_{\sym}$ and $U\in\R^{3\times 3}_{\sym}$ we have
  \begin{align*}
    &\inf\Big\{\int_{-\frac12}^{\frac12}Q\Big(x_3,\iota(x_3 A)+\mathbf 1(x_3>0)\tfrac{\bar r}{2} U+\iota(M)+\sym(d(x_3)\otimes e_3)\Big)\,\mathrm{d}x_3\,:\,\\
     &\qquad\qquad M\in\R^{2\times 2}_{\sym},\,d\in L^2((-\tfrac12,\tfrac12);\R^3)\,\Big\}\\
    &=Q_\mathrm{el}(A+\bar r\mathbb B(U'))+\bar r^2E_\mathrm{res}(U'),
  \end{align*}
  where $U'$ denotes the upper-left $2\!\!\times\!\!2$-submatrix of $U$ and the minimum is attained.
\end{lemma}
\begin{proof}
  We first note that 
  \begin{align*}
    &\inf_{\substack{M\in\R^{2\times 2}_{\sym}\\d\in L^2((-\frac12,\frac12);\R^3)}}\int_{-\frac12}^{\frac12}Q\Big(x_3,\iota(x_3 A)+\mathbf 1(x_3>0)\tfrac{\bar r}{2}U+\iota(M)+\sym(d(x_3)\otimes e_3)\Big)\,\mathrm d x_3\\
    =&\min_{M\in\R^{2\times 2}_{\sym}}\int_{-\frac12}^\frac12Q_2\Big(x_3,x_3A+ M+\mathbf 1(x_3>0)\tfrac{\bar r}{2}U'\Big)\,\mathrm d x_3\equalscolon(\star),
  \end{align*}
  where $Q_2$ is defined as in Definition \ref{def:QEB}.
  We next rewrite $(\star)$ based on a projection scheme that is inspired by \cite{BNS}.
  To that end, we denote by $\mathbf H$ the Hilbert space $L^2((-\frac12,\frac12);\R^{2\times 2}_{\sym})$ with norm $\|G \|_*^2\colonequals \int_{-\frac12}^\frac12 Q_2(x_3,G(x_3))\mathrm dx_3$.
  By $\mathbf H_\mathrm{ aff}\colonequals \{x_3A+M\,:\,A,M\in\R^{2\times 2}_{\sym}\}$ we denote the closed subspace of affine, $\R^{2\times 2}_{\sym}$-valued functions in $\mathbf H$, and by $\mathbf H_{2d}$ the orthogonal complement in $\mathbf H_\mathrm{ aff}$ of $\R^{2\times 2}_{\sym}$.
  We  write $P_\mathrm{ 2d}:\mathbf H\to\mathbf H_{2d}$ (resp. $P_\mathrm{ aff}:\mathbf H\to\mathbf H_\mathrm{ aff}$) for the orthogonal projection onto $\mathbf H_{2d}$ (resp. $\mathbf H_\mathrm{ aff}$).
  By the Pythagorean Theorem and $P_{2d}\circ P_\mathrm{aff}=P_{2d}$, for all $F\in\mathbf H$ we have,
  \begin{align*}
    \min_{M\in\R^{2\times 2}_{\sym}}\|F+M\|^2_*
    &= \min_{M\in\R^{2\times 2}_{\sym}}\|P_\mathrm{aff}F+M\|_*^2+\|(1-P_\mathrm{aff})F\|_*^2
     =\|P_{2d}F\|^2_*+\|(1-P_\mathrm{ aff})F\|^2_*\\
    &=\|P_{2d}F\|^2_*+\min_{A,M\in\R^{2\times 2}}\|x_3A+M+F\|^2_*.
  \end{align*}
  Applied with $F:(-\tfrac12,\tfrac12)\rightarrow \R^{2\times 2}_{\sym}$, $F(x_3)\colonequals x_3A+\mathbf 1(x_3>0)\tfrac{\bar r}{2}U'$, we obtain the identity
  \begin{equation}\label{eq:reduction:1}
    (\star)\;=\;
    		\int_{-\frac12}^\frac12 Q_2\left(x_3,P_{2d}F\right)\,\mathrm d x_3
           +
           \min_{A',M'\in\R^{2\times 2}_{\sym}}\int_{-\frac12}^\frac12 Q_2\Big(x_3,x_3A'+M'+\mathbf 1(x_3>0)\tfrac{\bar r}{2}U'\Big)\,\mathrm dx_3.
  \end{equation}
  Note that the second summand on the right-hand side equals $\bar r^2E_\mathrm{ res}(U')$.
  Thus, it remains to show that the first summand can be rewritten in the claimed form.
  To this end, we consider the map
  \begin{equation*}
    \Phi:\R^{2\times 2}_{\sym}\to\mathbf H_{2d},\qquad \Phi(A)\colonequals P_{2d}(x_3A),
  \end{equation*}
  where $x_3A$ stands short for the map $(-\tfrac12,\tfrac12)\ni x_3\mapsto	x_3A$.
  We note that $\Phi$ is an isomorphism.
  Indeed, $\Phi$ is linear and bounded and satisfies $\|\Phi(A)\|_*^2\geq \frac{1}{12\,C_W}|A|^2$ as a consequence of \eqref{eq:Q:1}. Since the dimensions of $\R^{2\times 2}_{\sym}$ and $\mathbf H_{2d}$ are the same, $\Phi$ must be an isomorphism. We claim that
  \begin{equation}\label{eq:id:B}
    \bar r\mathbb B(U')=(\Phi^{-1}\circ P_{2d})\Big(\mathbf 1(x_3>0)\tfrac{\bar r}{2}U'\Big)\qquad\text{for all }G'\in\R^{2\times 2}.
  \end{equation}
  Before we present the proof of the claim, we note that with help of \eqref{eq:id:B}, the first integral on the right-hand side in \eqref{eq:reduction:1} can be written in the form
  \begin{equation*}
    \int_{-\frac12}^\frac12 Q_2\Big(x_3,P_{2d}\big(x_3A+x_3\bar r\mathbb B(U')\big)\Big)=Q_\mathrm{ el}(A+\bar r\mathbb B(U')),
  \end{equation*}
  which completes the proof of Lemma \ref{L:reduction}.
  It remains to show \eqref{eq:id:B}. For the argument, we recall the notation $G_1,\ldots,G_3$ for the standard basis of $\R^{2\times 2}_{\sym}$ and the definition of $A_i,M_i$, cf.~Definition~\ref{def:QEB}. We claim that
  \begin{equation}\label{eq:reduction:2}
    \Phi(A_i)=P_{2d}\Big(\mathbf 1(x_3>0)\tfrac{1}{2}G_i\Big),\qquad i=1,2,3.
  \end{equation}
  Indeed, in view of the definition of $A_i,M_i$, we have $P_\mathrm{ aff}\Big(\mathbf 1(x_3>0)\tfrac{1}{2}G_i\Big)=x_3A_i+M_i$.
  Hence, in view of $P_{2d}=P_{2d}\circ P_\mathrm{ aff}$, we have $P_{2d}\Big(\mathbf 1(x_3>0)\tfrac{1}{2}G_i\Big)=P_{2d}(x_3A_i+M_i)=P_{2d}(x_3 A_i)=\Phi(A_i)$, and thus \eqref{eq:reduction:2}.
  By applying the isomorphism $\Phi$ to both sides of \eqref{eq:id:B} and by appealing to the definition of $\mathbb B$, we see that \eqref{eq:id:B} is in fact equivalent to \eqref{eq:reduction:2}.
\end{proof}

\begin{proof}[Proof of Lemma~\ref{L:representations:hom:iso}]
In this proof, we write $Q(\cdot)\colonequals Q(x_3,\cdot)$ and $Q_2(\cdot) = Q_2(x_3,\cdot)$, emphasizing the independence from $x_3$.

By definition of $Q_\mathrm{ el}$ and by orthogonality, we indeed have
\begin{equation}\label{eq:Aggenstein}
Q_\mathrm{ el}(A)=\min_{M\in\R^{2\times 2}_{\sym}}\left( \int_{-\frac12}^{\frac12} Q_2(x_3A)+Q_2(M) \right)= 
\frac{1}{12}Q_2(A).
\end{equation}

Next, we note that for each $A,M,T\in \R^{2\times2}_{\sym}$, we have
\begin{equation}\label{eq:Das:Haus:Vom:Nicolaus}
\begin{aligned}
&\int_{-\frac12}^{\frac12}Q_2\left(\mathbf 1(x_3>0)T-(x_3A+M)\right)\,\mathrm dx_3\\
=& Q_2\big(\tfrac12T-M\big)+\frac{1}{12}Q_2\big(\tfrac32T-A\big)+\frac1{16}Q_2\big(T\big),
\end{aligned}
\end{equation}
which can be seen by using the $L^2((-\frac12,\frac12))$-orthogonal decomposition $\mathbf 1(x_3>0)=\frac12+\frac32x_3+(\mathbf 1(x_3>0)-\frac12 -\frac32x_3)$.
Using \eqref{eq:Das:Haus:Vom:Nicolaus} with $T\colonequals \frac12 U$ and \eqref{eq:Aggenstein}, we obtain the claimed form of $E_\mathrm{res}(U)$.
Using \eqref{eq:Das:Haus:Vom:Nicolaus} with $T\colonequals \frac12G_i$ and $G_i$ from Definition \ref{def:QEB}, we obtain that $A_i$ from Definition~\ref{def:QEB} takes the form $A_i=\frac34G_i$.
Now,\eqref{eq:repr:hom:iso:B} directly follows from the definition of $\mathbb B$.

For the second part of the proof, we assume the material to be isotropic.
Then, the minimizer $d$ in \eqref{eq:repr:hom:iso:el} is given by
\begin{equation*}
d=-\frac{\lambda}{2\mu+\lambda}\operatorname{tr} (A)e_3,
\end{equation*}
which leads to the claimed form of $Q_\mathrm{ el}$.
The remaining representation for $E_\mathrm{ res}$ follows immediately with \eqref{eq:repr:hom:iso:Eres}.
\end{proof}

\begin{proof}[Proof of Lemma \ref{L:properties}]
\step{1 -- Proof of (a) and (b)}\smallskip
Note that by \eqref{eq:Q:1} we have $\frac{1}{C_W}\bar Q(G)\leq Q_2(x_3,G)\leq C_W\bar Q(G)$, where $\bar Q(G) \colonequals  |\sym G|^2$.
Hence, the claims of (a) and (b) follow with help of \eqref{eq:Das:Haus:Vom:Nicolaus} applied to $Q_2:=\bar Q$.
 
\step{2 -- Proof of (c)}\smallskip

The linearity of $\mathbb B$ is obvious.
For the boundedness, let $A_i$ and $M_i$ be defined as in Definition~\ref{def:QEB} and note that
\begin{equation*}
  \begin{aligned}
    \frac1{12C_W}|\tfrac34 G_i-A_i|^2\leq   \,&\,\frac{1}{C_W}\int_{-\frac12}^\frac12|\sym(\mathbf 1(x_3>0)\tfrac{1}{2}G_i-(x_3 A_i+M_i))|^2\\
    \leq \,&\,\int_{-\frac12}^{\frac12}Q_2\big(x_3, \mathbf 1(x_3>0)\tfrac{1}{2}G_i-(x_3 A_i+M_i)\big) = E_\mathrm{res}(G_i)\\
    \leq\,& \frac{C_W}{64}|G_i|^2,
\end{aligned}
\end{equation*}
where the first estimate holds by \eqref{eq:Das:Haus:Vom:Nicolaus} applied to $Q_2:=|\sym(\cdot)|$, the second estimate holds thanks to \eqref{eq:Q:1}, and the last estimate holds by (b). Hence, we conclude that $|A_i|\leq\frac{\sqrt 3}{4}(C_W+\sqrt 3)$.
\end{proof}

\subsection{Lower bound: proof of Theorem~\ref{T1}~\ref{item:T1:lower:bound}}
\label{Sec:Proof:Lower:Bound}

\begin{proof}[Proof of Theorem~\ref{T1}~\ref{item:T1:lower:bound}]
  It suffices to consider the case $\liminf_{h\to 0}\mathcal E_h(y_h,n_h)<\infty$. By appealing to Lemma~\ref{L:compfinitebend} and Theorem~\ref{T1}~\ref{item:T1:compactness}, we may pass to a subsequence (that we do not relabel) such that
  \begin{alignat}{5}\label{eq:st:10a}
    \nabla\!_hy_h&\to&&R_y&&\text{ strongly in }L^2(\Omega),\\\label{eq:st:10b}
    E_h(y_h)&\wto&& E\colonequals \iota(x_3\II_y+M)+\sym(d\otimes e_3)&&\text{ weakly in }L^2(\Omega),\\\label{eq:st:10c}
    \nabla\!_hn_h&\wto&& (\nabla'n,d_n)&&\text{ weakly in }L^p(\OmegaT),\\\label{eq:st:10d}
    n_h&\to&&n&&\text{ a.e.~in }\OmegaT,
  \end{alignat}
  for some $M\in L^2(S;\R^{2\times 2}_{\sym})$, $d\in L^2(\Omega;\R^3)$ and $d_n\in L^2(\OmegaT;\R^3)$. We recall the shorthand notation of $B_h$, see \eqref{eq:st:9:Bh}, and note that \eqref{eq:st:10d} yields
  \begin{equation}\label{eq:st:10e}
    B_h\to {\mathbf 1}(x_3>0)B\quad\text{a.e.~in }\Omega,\text{ where }    B\colonequals \tfrac{\bar r}{2}(\tfrac13 I-n\otimes n).
  \end{equation}
  Furthermore, in view of \eqref{eq:st:0} and the polar factorization, for all sufficiently small $h>0$ there exists a measurable $R_h:\Omega\to\SO 3$ such that $\nabla\!_hy_h=R_h(I+h E_h(y_h))$. Hence, 
  \begin{equation*}
    R_h^\top(I+hB_h)\nabla\!_hy_h=I+h\underbrace{\big(E_h(y_h)+R_h^\top B_hR_h+hR_h^\top B_hR_hE_h(y_h)\big)}_{\equalscolon G_h}.
  \end{equation*}
  Since \eqref{eq:st:10a} implies $R_h\to R_y$ a.e.~in $\Omega$, and since $(R_h),(B_h)$ are bounded in $L^\infty(\Omega)$, we conclude from \eqref{eq:st:10b}, \eqref{eq:st:10e}, and the convergence of $R_h$ that
  \begin{equation*}
    G_h\wto \iota(x_3\II_y+M)+\sym(d\otimes e_3)+{\mathbf 1}(x_3>0)R_y^\top B\,R_y\qquad\text{weakly in }L^2(\Omega).
  \end{equation*}
  Hence, arguing as in \cite[Theorem 6.1 (i)]{FJM02} we conclude that
  \begin{align*}
&    \liminf\limits_{h\to 0}\frac{1}{h^2}\int_\Omega W(x_3,(I+hB_h)\nabla\!_hy_h)\,=\,\liminf\limits_{h\to 0}\frac{1}{h^2}\int_\Omega W(x_3,I+hG_h)\\
    \geq& \int_\Omega Q\big(x_3,\iota(x_3\II_y+M)+\sym(d\otimes e_3)+{\mathbf 1}(x_3>0)R_y^\top B\,R_y\big)\\
    \geq& \int_S\Big(\inf_{\tilde M\in\R^{2\times 2}_{\sym}\atop \tilde d\in L^2((-\frac12,\frac12);\R^3)}\int_{-\frac12}^\frac12 Q\big(x_3,\iota(x_3\II_y+\tilde M)+\sym(\tilde d(x_3)\otimes e_3)+{\mathbf 1}(x_3>0)R_y^\top B\,R_y\big)\,\mathrm dx_3\Big)\,\mathrm dx'.
  \end{align*}
  Note that the upper-left $2\!\times\!2$-submatrix of $R_y^\top B \,R_y$ is given by $\frac{\bar r}{2}(\frac{1}{3}I-\hat n'\otimes\hat n')$ where $\hat n'\colonequals \nabla'y^\top n$. Hence, combined with Lemma~\ref{L:reduction} we conclude that
  \begin{equation*}
    \begin{aligned}
    \liminf\limits_{h\to 0}\frac{1}{h^2}\int_\Omega W(x_3,(I+hB_h)\nabla\!_hy_h)\,\geq\,
    &\int_S Q_\mathrm{ el}(\II_y+\bar r\mathbb B(\tfrac13 I_{2\times 2}-\hat n'\otimes\hat n'))+\bar r^2E_\mathrm{res}(\tfrac13 I_{2\times 2}-\hat n'\otimes \hat n')\\
    =&\,\mathcal E_\mathrm{ el}(y,n).
  \end{aligned}
  \end{equation*}
  It remains to show that $\liminf\limits_{h\to 0}\int_{\OmegaT}|\nabla\!_hn_h(\nabla\!_hy_h)^{-1}|^2\det(\nabla\!_hy_h)\geq \mathcal E_{OF}(n)$.
  By the claim of Step~2 in the proof of Theorem~\ref{T1}~\ref{item:T1:compactness}, cf.~\eqref{eq:st:20}, we conclude that
  \begin{equation*}
    \liminf\limits_{h\to 0}\int_{\OmegaT}|\nabla\!_hn_h(\nabla\!_hy_h)^{-1}|^2\det(\nabla\!_hy_h)\geq \int_{\OmegaT}|(\nabla'n,d_n)|^2.
  \end{equation*}
  Since $|(\nabla'n,d_n)|^2=|\nabla'n|^2+|d_n|^2$ and because $n$ is independent of $x_3$, the right-hand side is bounded from below by $\frac12\int_S|\nabla'n|^2\,dx'=\mathcal E_{OF}(n)$,
  and the claimed lower bound follows.
\end{proof}

\subsection{Recovery sequence: proof of Theorem~\ref{T1}~\ref{item:T1:upper:bound} and Theorem~\ref{T2}~\ref{item:T2:Recov}}
\label{Sec:Proof:Upper:Bound}
In this section, we present the construction of recovery sequences.
We only discuss the case with prescribed boundary conditions, i.e., Theorem~\ref{T2}~\ref{item:T2:Recov}, since it requires an additional argument compared to the case without boundary conditions.
In contrast to earlier works in the field, in our situation, constructing the recovery sequences invokes two tasks:
the construction of the deformation $y_h$ and of the director field $n_h$.
Note that the energy invokes a nonlinear coupling between $\nabla\!_hy_h$ and $n_h\otimes n_h$, which is the reason why the recovery sequences for $y_h$ and $n_h$ cannot be constructed independently.
However, by exploiting the good compactness properties of the lower-order term $n_h\otimes n_h$, for the recovery sequence, we can first construct a sequence of 3d-directors $(n_h)$ that recovers a limiting director field $n$,
and in a second step, we can then construct an adapted sequence of 3d-deformations $(y_h)$ (that recovers an $n$-dependent limiting strain).
\medskip

We first discuss the construction of recovery sequences for deformations.
This amounts to finding a sequence of 3d-deformations $(y_h)$ such that the associated sequence of nonlinear strains
\begin{equation*}
  E_h(y_h)=\frac{\sqrt{\nabla\!_hy_h^\top\nabla\!_hy_h}-I}{h}
\end{equation*}
strongly converges in $L^2(\Omega)$ to a limiting strain of the form obtained in the compactness statement Lemma~\ref{L:compfinitebend}, cf.~\eqref{eq:st:3b}, i.e.,
\begin{equation*}
  E\colonequals \iota(x_3\II_y+M)+\sym(d\otimes e_3),
\end{equation*}
for a prescribed isometry $y\in H^2_{\iso}(S;\R^3)$, some $M\in L^2(S;\R^{2\times 2}_{\sym})$ and $d\in L^2(\Omega;\R^3)$.
As in previous works on the derivation of bending theories, we use that smooth isometries are dense in $H^2_{\iso}(S;\R^3)$, see Proposition~\ref{P:approx} for a statement respecting boundary conditions.
The following proposition summarizes the construction for smooth 3d-deformations:
\begin{proposition}\label{P:general_construction}
  Let $S$ satisfy \eqref{ass:domain}. Let $y\in H^2_\mathrm{ iso}(S;\R^3)\cap C^\infty(\overline S;\R^3)$, $M\in L^2(S;\R^{2\times 2}_{\sym})$, $d\in L^2(\Omega;\R^3)$ and assume that there exists $\delta\geq 0$ such that
  \begin{equation}\label{eq:st:40}
    M+\delta I_{2\times 2}\geq 0\text{ a.e.~in }\{\II_y=0\}.
  \end{equation}
  Let $0<\beta<\frac12$.
  Then there exists a sequence $(y_h)\subseteq C^{\infty}(\overline\Omega;\R^3)$ such that
  \begin{equation}\label{eq:recov:onescale}
    \begin{aligned}
      &\limsup\limits_{h\to 0}
      	\,\|y_h-y\|_{L^\infty(\Omega)}=0,\\
      &\limsup\limits_{h\to0}\|E_h(y_h)\,-\,\big(\iota(x_3 \II_y+M)+\sym(d\otimes e_3)\big)\|_{L^2(\Omega)}=0,\\
      &\limsup\limits_{h\to 0} h^{-\beta}\|\nabla\!_h y_h-R_y\|_{L^\infty(\Omega)}=0.
    \end{aligned}
  \end{equation}
  Furthermore, if Assumption~\ref{ass:BC} is satisfied and $y\in\mathcal A_{BC}\cap C^\infty(\overline S;\R^3)$, then we may additionally enforce the clamped, affine boundary condition $y_h\in\mathcal A_{BC,\delta,h}$ (cf.~\eqref{eq:BC3dA}).
\end{proposition}
\textit{(See Section~\ref{Sec:Proof:P:general_construction} for the proof.)}\medskip

Results similar to Proposition~\ref{P:general_construction} have been obtained in the context of dimension reduction earlier:
The special case $M=0$ in \eqref{eq:recov:onescale} is already contained in the seminal paper \cite{FJM02}.
Schmidt~\cite{SCHMIDT07} contains a similar result in the case without boundary conditions and under the additional assumption that $M=0$ in an open neighborhood of the flat region $\{\II_y=0\}$.
Finally, the fourth author introduced in \cite{DPG20} a wrinkling construction (which is inspired by \cite{LP17}) to construct deformations in the case when $M$ may not vanish on flat regions, yet without respecting boundary conditions.
Based on this, to treat boundary conditions, in Section~\ref{Sec:Proof:P:general_construction} we present an argument that for each connected component $V$ of $\{\II_y\neq 0\}$ and $\interieur\{\II_y=0\}$ constructs a recovery sequence with prescribed clamped, affine boundary conditions on $\partial V\cap S$. 
With Propositions \ref{P:approx} and \ref{P:general_construction} at hand, we are in a position to prove Theorem~\ref{T2}~\ref{item:T2:Recov}. In fact, we shall establish a stronger statement:
\begin{lemma}\label{L:stronger}
  Let Assumption \ref{ass:W} and \eqref{ass:domain} be satisfied. 
Let $(y,n)\in \mathcal A_2$ and let $(n_h)\subseteq H^1(\OmegaT;\mathbb S^2)$ denote an arbitrary sequence satisfying
\begin{equation*}
  (n_h,\nabla\!_h n_h)\to (n,\nabla'n,0)\qquad\text{strongly in }L^2(\OmegaT),
\end{equation*}
(e.g.,~$n_h\colonequals n$).
Then there exists $(y_h)\subseteq H^1(\Omega;\R^3)$ such that $(y_h,\nabla\!_hy_h)\to (y,R_y)$ strongly in $L^2(\Omega)$ and $\mathcal E_h(y_h,n_h)\to\mathcal E(y,n)$.
Furthermore, let Assumption~\ref{ass:BC} be satisfied, let $y\in\mathcal A_{BC}$ and $\bar \delta$ satisfy \eqref{eq:bardelta}, then we may additionally enforce $y_h\in\mathcal A_{BC,\bar\delta,h}$.
\end{lemma}

\begin{proof}[Proof of Lemma~\ref{L:stronger} and \mbox{Theorem~\ref{T2}~\ref{item:T2:Recov}}]
  It suffices to prove Lemma~\ref{L:stronger} and for brevity we only consider the case with prescribed boundary conditions. The construction of $y_h$ invokes an approximation argument: By Proposition~\ref{P:approx}, for all $k\in\N$ we find $y^k\in\mathcal A_{BC}\cap C^\infty(\overline S;\R^3)$ with $\|y^k-y\|_{H^2(S;\R^3)}<\frac1k$.
  As an abbreviation, we write $R^k\colonequals R_{y^k}$ and $\II^k\colonequals \II_{y^k}$.
    
  In view of Lemma~\ref{L:reduction} there exist $M^k\in L^2(S;\R^{2\times 2}_{\sym})$, $d^k\in L^2(\Omega;\R^3)$ such that
  \begin{equation}\label{eq:st:41}
    \begin{aligned}
      &\,\int_S Q_\mathrm{ el}(\II^k+\bar r\mathbb B(\tfrac13 I_{2\times 2}-(\hat n^k)'\otimes(\hat n^k)'))+\bar r^2E_\mathrm{ res}(\tfrac13 I_{2\times 2}-(\hat n^k)'\otimes(\hat n^k)')\,\mathrm dx'\\
    =\,&\,
        \int_{\Omega}Q\Big(x_3,\iota(x_3 \II^k)+B^k+\iota(M^k)+\sym(d^k\otimes e_3)\Big)\,\mathrm dx,
      \end{aligned}
  \end{equation}
  where
  \begin{equation*}
    \hat n^k\colonequals (R^k)^\top n\qquad\text{and}\qquad B^k\colonequals \mathbf 1(x_3>0)\tfrac{\bar r}{2}\big(\tfrac13 I_{3\times 3}-\hat n^k\otimes\hat n^k\big).
  \end{equation*}
In particular, Lemma \ref{L:reduction} yields for a.e.~$x'$ with $\II^k(x')=0$,
\begin{equation}\label{eq:Turm:am:verborgenen:Horn}
\begin{aligned}
&\int_{-\frac12}^{\frac12}Q\Big(x_3,B^k(x',x_3)+\iota(M^k(x'))+\sym(d^k(x',x_3)\otimes e_3)\Big)\,\mathrm dx_3
\\
=& \min_{M\in\R^{2\times 2}_{\sym}\atop d\in L^{2}((-\frac12,\frac12);\R^3)}
	\int_{-\frac12}^{\frac12}Q\Big(x_3,B^k(x',x_3)+\iota(M)+\sym(d(x_3)\otimes e_3)\Big)\,\mathrm dx.
\end{aligned}
\end{equation}
We claim that $M^k$ satisfies \eqref{eq:st:40} (up to a null-set) for $\delta\colonequals \bar \delta$.
Indeed, \eqref{eq:Q:1} and a calculation yield
\begin{equation*}
[\mathrm{RHS}\text{ of }\eqref{eq:Turm:am:verborgenen:Horn}]
\leq \min_{M\in\R^{2\times 2}_{\sym}} C_W\int_{-\frac12}^{\frac12}|\mathbf 1(x_3>0)T^k(x')+M|^2\;\mathrm d x_3 = C_W \frac14|T^k(x')|^2,
\end{equation*}
where $T^k(x')\colonequals  \frac{\bar r}{2}\big(\frac13 I_{2\times 2}-(\hat n^k)'(x')\,\otimes\, (\hat n^k)'(x')\big)$ and where $(\hat n^k)'$ contains the first two components of $\hat n^k$.
Also, \eqref{eq:Q:1} yields
\begin{align*}
[\mathrm{LHS}\text{ of }\eqref{eq:Turm:am:verborgenen:Horn}]&\geq
\frac{1}{C_W}\int_{-\frac12}^{\frac12}|\mathbf 1(x_3>0)T^k(x')+M^k(x')|^2\,\mathrm d x_3 \\
&= \frac{1}{C_W} \big(|\tfrac12 T^k(x')+M^k(x')|^2 + \tfrac14 |T^k(x')|^2\big).
\end{align*}
We conclude that
\begin{equation*}
|M^k(x')|\leq \frac{1}{\sqrt 2} C_W|T^k(x')| \leq {\bar r}C_W\frac{1+\sqrt{2}}{4}.
\end{equation*}
Thus, \eqref{eq:st:40} is satisfied for $M:=M^k$ and $\delta \colonequals \bar\delta$.
  Finally, we may apply Proposition~\ref{P:general_construction} and obtain for each $k$ a sequence $(y^k_h)\subseteq\mathcal A_{BC,\bar\delta,h}$ such that
  \begin{equation}\label{eq:st:40:c}
    \begin{aligned}
    \lim\limits_{h\to 0}\Big(\|y^k_h-y^k\|_{L^2(\Omega)}+\|E_h(y^k_h)-\underbrace{(\iota(x_3\II^k+M^k)+\sym(d^k\otimes e_3))}_{\equalscolon E^k}\|_{L^2(\Omega)}\\
    +h^{-\beta}\|\nabla\!_hy^k_h-{R^k}\|_{L^\infty(\Omega)}\Big)=0,
  \end{aligned}
  \end{equation}
  where $0<\beta<\frac12$ denotes an arbitrary exponent, that is fixed from now on. Next, we claim that
  \begin{equation}\label{eq:st:42}
    \lim\limits_{h\to 0}\mathcal E_h(y^k_h,n_h)=\mathcal E(y^k,n).
  \end{equation}
  This can be seen as follows:
  Thanks to the convergence of $\nabla\!_h y^k_h$ with rate $h^\beta$ and the elementary inequality $|E_h(y^k_h)|\leq h^{-1}\dist(\nabla\!_h y^k_h,\SO 3)$, we have $\|hE_h(y^k_h)\|_{L^\infty(\Omega)}\to 0$.
  Furthermore, by polar factorization, we find $R^k_h:\Omega\to\SO 3$ such that $\nabla\!_h y_h^k=R_h^k(I+hE_h(y_h^k))$.
  Let $B_h$ be defined by \eqref{eq:st:9:Bh} and note that $B_h\to {\mathbf 1}(x_3>0)\tfrac{\bar r}{2}(\tfrac13 I_{3\times 3}-n\otimes n)$ a.e.~in $\Omega$ (c.f.~\eqref{eq:st:10e}).
  We obtain
  \begin{equation*}
    (I+hB_h)\nabla\!_h y^k_h=R_h^k(I+h G_h^k),\qquad \text{where } G_h^k\colonequals (R^k_h)^\top B_h(\nabla\!_hy_h^k)+E_h(y_h^k).
  \end{equation*}
  Thanks to \eqref{eq:st:40:c} and $\|hE_h(y^k_h)\|_{L^\infty(\Omega)}\to 0$ we have $R_h^k\rightarrow R^k$ strongly in $L^\infty(\Omega)$ and
  \begin{alignat*}{3}
&    G_h^k\to G^k\colonequals B^k+E^k &&\qquad\text{for $h\to 0$ in $L^2(\Omega)$},\\
&     hG_h^k\rightarrow 0  &&\qquad\text{for $h\to 0$ in $L^\infty(\Omega)$},
  \end{alignat*}
  and we conclude with help of \ref{item:ass:W:frame:ind}, \ref{item:ass:W:quadratic:exp} and \eqref{eq:st:41} that
  \begin{align*}
    &\lim\limits_{h\to 0}\frac{1}{h^2}\int_\Omega W(x_3,    (I+hB_h)\nabla\!_h y^k_h) =    \lim\limits_{h\to 0}\frac{1}{h^2}\int_\Omega W(x_3, I+hG_h^k)=\int_\Omega Q(x_3,G^k)\\
    &=\int_S Q_\mathrm{ el}\big(\II^k+\bar r\mathbb B(\tfrac13 I_{2\times 2}-(\hat n^k)'\otimes(\hat n^k)')\big)+\bar r^2E_\mathrm{ res}\big(\tfrac13 I_{2\times 2}-(\hat n^k)'\otimes(\hat n^k)'\big)\,\mathrm dx'\\
    &=\mathcal E_\mathrm{ el}(y^k,n).
  \end{align*}
  Since $\nabla\!_hn_h\to (\nabla'n,0)$ strongly in $L^2$ by assumption, we have $\int_{\OmegaT}|\nabla\!_hn_h|^2\to\frac12\int_S|\nabla'n|^2=\mathcal E_{OF}(n)$ and thus \eqref{eq:st:42} follows.

  Finally, we obtain the sought for recovery sequence $(y_h,n_h)$ by extraction of a diagonal sequence: For $h>0$ and $k\in\N$ consider
  \begin{equation*}
    c_{k,h}\colonequals \|y_h^k-y\|_{L^2(\Omega)}+\|\nabla\!_hy_h^k-R_y\|_{L^2(\Omega)}+|\mathcal E_h(y_h^k,n_h)-\mathcal E(y,n)|.
  \end{equation*}
  Then, $\limsup\limits_{h\to 0}c_{k,h}\leq \|y^k-y\|_{L^2}+\|R^k-R_y\|_{L^2}+|\mathcal E(y^k,n)-\mathcal E(y,n)|$.
  In view of $y^k\rightarrow y$ strongly in $H^2(S)$, and in view of Lemma~\ref{L:cont}~(d) we conclude that $\limsup\limits_{k\to\infty}\limsup\limits_{h\to 0} c_{k,h}=0$.
  Hence, we find a diagonal sequence $k_h\to\infty$ as $h\to 0$ such that $c_{k_h,h}\to 0$ and we conclude that $y_h\colonequals y_h^{k_h}$ is the sought for recovery sequence.  
\end{proof}

\subsection{Properties of the limit: proof of Lemma \ref{L:cont}}
\label{Sec:Proof:Properties:of:Limit}

\begin{proof}[Proof of Lemma \ref{L:cont}]
\step{1 -- Proof of (a)}\smallskip

In view of Lemma \ref{L:properties} we have $|\II_y|^2\leq C \Big( Q_{\rm el}\big(\II_y+\bar r\mathbb B(\tfrac13 I_{2\times 2}-\nabla'y^\top n\, \otimes\,\nabla'y^\top n)\big)+1\Big)$ a.e.~in $S$, where $C$ can be chosen only depending on $C_W$ and $\bar r$.
Combined with the well-known identity
\begin{equation}\label{eq:II:D2y}
|\II_y|^2=|\nabla'\nabla' y|^2 \text{ a.e.~in }S,
\end{equation}
which in fact holds for all $y\in H^2_{\iso}(S;\R^3)$,
the estimate of (a) follows.

\step{2 -- Compactness properties of $\mathcal A_2$}\smallskip

We claim that
\begin{enumerate}[(i)]
\item $H^2_{\iso}(S;\R^3)$ and $H^1(S;\mathbb S^2)$ are closed w.r.t.~weak convergence in $H^2(S)$ and $H^1(S)$, respectively. 
\item If $(y^k)\subseteq H^2_{\iso}(S;\R^3)$ converges to $y\in H^2_{\iso}(S;\R^3)$ weakly in $H^2(S)$, then $\II_{y^k}\wto \II_y$ weakly in $L^2(S)$.
\item From any sequence $(y^k,n^k)\in\mathcal A_2$ with $\liminf\limits_{k\to\infty}\big(\mathcal E(y^k,n^k)+\|y^k\|_{L^2(S)}\big)<\infty$ we can extract a subsequence that converges weakly in $H^2(S)\times H^1(S)$ to a limit $(y,n)\in\mathcal A_2$.
\end{enumerate}
To see (i), let $(y^k)\subseteq H^2_{\iso}(S;\R^3)$ denote a sequence that weakly converges to some  $y\in H^2(S)$.
Then $\nabla'y^{k} \to \nabla'y$ strongly in $L^2(S)$ (by compact embedding) and thus $(\nabla'y^{k})^\top\nabla'y^{k} \to \nabla'y^\top\nabla'y$ strongly in $L^1(S)$.
Hence, since $(\nabla'y^k)^\top\nabla'y^k=I$ a.e., we get $\nabla'y^\top\nabla'y = I$ a.e., and thus $y\in H^2_{\iso}(S;\R^3)$.
The argument for the weak closedness of $H^1(S;\mathbb S^2)$ is similar and left to the reader.
For the argument of (ii), note that we have $\nabla'y^k\to \nabla'y$ strongly in $L^2(S)$ by compact embedding.
Since $|\partial_1y^k|=|\partial_2y^k|=1$ a.e., we have $|b_{y_k}|\leq 2$ and $|\nabla'b_{y_k}|\leq 2|\nabla'\nabla'y_k|$ a.e.
Hence, we conclude that $b_{y^k}\wto b_y$ weakly in $H^1(S)$.
Combined with the boundedness of $(\II_{y^k})$ in $L^2(S)$ (c.f.~\eqref{eq:II:D2y}), we obtain $\II_{y^k}=\nabla'y_k^\top\nabla' b_{y^k}\wto \nabla' y^\top\nabla' b_y=\II_y$ weakly in $L^2(S)$.
We finally prove (iii).
To that end, note that from  $\liminf\limits_{k\to\infty}(\mathcal E(y^k,n^k)+\|y^k\|_{L^2(S)})<\infty$, we conclude with help of (a) that $\|\nabla'\nabla' y^k\|_{L^2(S)}+\|y^k\|_{L^2(S)}+\|\nabla' n^k\|_{L^2(S)}\leq C$ for a subsequence (not relabeled) and a constant $C$.
Hence, by Poincar\'e's inequality, we infer that the sequences $(y^k)$ and $(n^k)$ are bounded in  $H^2(S)$ and $H^1(S)$, respectively.
Thus, we can extract a subsequence that weakly converges in $H^2(S)\times H^1(S)$ to a limit $(y,n)$. By (i), $(y,n)\in\mathcal A_2$ follows.

\step{3 -- Proof of (c)}\smallskip

Assume w.l.o.g.~that {$\liminf_{k\to\infty} \mathcal E_{OF}(n^k)=\limsup_{k\to\infty} \mathcal E_{OF}(n^k)<\infty$}.
Then $(n^k)$ converges weakly in $H^1(\OmegaT)$, and the liminf-inequality for $\mathcal E_{OF}$ is an immediate consequence of the weak lower semicontinuity of the $L^2$-norm.
For the liminf-inequality of $\mathcal E_{\rm el}$, we assume w.l.o.g.~that 
$
\liminf_{k\to\infty}\mathcal E_{\rm el}(y^k,n^k) < \infty
$.
Then, by the assumption $\limsup_{k\to\infty} \mathcal E_{OF}(n^k)<\infty$ and by Step~2~(iii), we have $(y^k,n^k)\wto (y,n)$ weakly in $H^2(S)\times H^1(S)$.
We conclude that $\nabla'y^k\to \nabla'y$ strongly in $L^2(S)$, and since $|\nabla'y^k| =\sqrt{2}$, we have $\nabla'y^k\to \nabla'y$ strongly in $L^p(S)$ for every $1\leq p<\infty$.
Similarly, we have that $n^k\to n$ strongly in $L^{p}(S)$ for every $1\leq p<\infty$. 
This implies $\mathbb B(\tfrac13 I_{2\times 2}-(\nabla'y^{k})^\top n^{k}\otimes (\nabla'y^{k})^\top n^{k}) \to \mathbb B(\tfrac13 I_{2\times 2}-\nabla'y^\top n\otimes \nabla'y^\top n)$ strongly in $L^2(S)$.
By Step~2 (ii) we further have  $\II_{y^k}\wto \II_y$ weakly in $L^2(S)$.
Thus, the weak lower semicontinuity of convex integral functionals w.r.t.~weak convergence in $L^{2}$ yields
\begin{equation*}
\liminf\limits_{k\to\infty}\mathcal E_{\rm el}(y^k,n^k)\geq \mathcal E_{\rm el}(y,n).
\end{equation*}

\step{4 -- Proof of (b)}\smallskip

Consider the sublevel $\mathcal C\colonequals \{(y,n)\in\mathcal A_2\,:\,\mathcal E(y,n)\leq C\text{ and }\|y\|_{L^2(S)}\leq C\}$, and a sequence $(y^k,n^k)\subseteq\mathcal C$.
Then by Step~2~(iii) we can pass to a subsequence (not relabeled) such that $(y^k,n^k)\wto (y,n)$ weakly in $H^2(S)\times H^1(S)$ with $(y,n)\in\mathcal A_2$.
By (c) we have $\mathcal E(y,n)\leq\liminf\limits_{k\to\infty}\mathcal E(y^k,n^k)\leq C$, and by the lower semicontinuity of the $L^2$-norm, $\|y\|_{L^2(S)}\leq C$. Thus, $(y,n)\in\mathcal C$.
We conclude that $\mathcal C$ is compact w.r.t.~the weak topology of $H^2(S)\times H^1(S)$.
The relative compactness of $\mathcal C$ w.r.t.~the strong topology of $L^2(S)\times L^2(S)$ follows with help of the compact embedding of $H^1(S)\subseteq L^2(S)$.

\step{5 -- Proof of (d)}\smallskip

 Note that if $y^k\to y$ strongly in $H^2(S)$ and $n^k\wto n$ weakly in $H^{1}(S)$ then, proceeding as in Step~3, $\nabla' y^{k} \to \nabla' y$ strongly in $L^{p}(S)$ for every $1\leq p < \infty$, and  $\mathbb B(\tfrac13 I_{2\times 2}-(\nabla'y^{k})^\top n^{k}\otimes (\nabla'y^{k})^\top n^{k}) \to \mathbb B(\tfrac13 I_{2\times 2}-\nabla'y^\top n\otimes \nabla'y^\top n)$ strongly in $L^{2}(S)$.
 Furthermore, by Step 2 (ii) we have that $\II_{y^k}\wto \II_y$  weakly in $L^2(S)$.
 Since in addition $\lim_{k \to \infty} \|\II_{y^k} \|_{L^{2}} = \|\II_{y} \|_{L^{2}}$ (because of \eqref{eq:II:D2y}), we have that $\II_{y^k} \to \II_y$ strongly in $L^2(S)$.
Then the continuity of quadratic integral functionals w.r.t.~strong convergence in $L^{2}$ implies that 
$
\mathcal E_{\rm el}(y^k,n^k)\to\mathcal E_{\rm el}(y,n)
$.
The convergence of $\mathcal E_{OF}(n^k)$ is immediate from the definition of strong convergence.

\end{proof}

\subsection{Anchorings}
\label{Sec:Proof:Anchoring}
\subsubsection{Proof of Lemmas~\ref{P:Gamma:Strong:Anchoring} and \ref{L:volweakanch}}
\label{Sec:Proof:213:214}
\begin{proof}[Proof of Lemma~\ref{P:Gamma:Strong:Anchoring}]
  Part (a) follows from Theorem~\ref{T1} (a), which in particular yields $n_h\wto n$ weakly in $W^{1,p}(\OmegaT)$ for some $p\in(1,2)$, and thus $\Psi_{BC}(x',n_h(x'))\wto \Psi_{BC}(x',n(x'))$ weakly in $W^{1,p}(\OmegaT)$ thanks to the boundedness of the trace operator $T:W^{1,p}(\OmegaT;\R^3)\to L^p(\Gamma_n\times(0,\frac12);\R^3)$ and thanks to $\Psi_{BC}$ being affine. Part (b) follows from Lemma~\ref{L:stronger} with $n_h\colonequals n$.
\end{proof}
\begin{proof}[Proof of Lemma~\ref{L:volweakanch}]
  We note that the compactness statements of part (a) and part (b) directly follow from Theorem~\ref{T1} (a), since we have $\mathcal G_h\geq 0$ and $\mathcal H_h\geq 0$.

  \step{1 -- Convergence of the weak anchoring terms}\smallskip

  We claim that for any sequence $(y_h,n_h)\subseteq H^1(\Omega;\R^3)\times H^1(\OmegaT;\mathbb S^2)$ with
  \begin{equation}\label{eq:st:50c}
  \begin{aligned}
    &\det(\nabla\!_hy_h)>0 \text{ a.e.~in }\Omega\text{, and }\\
    &(\nabla\!_hy_h,n_h)\to (R_y,n)\text{ in }L^2\text{ for some }(y,n)\in\mathcal A_2,
    \end{aligned}
  \end{equation}
  we have
  \begin{alignat}{2}\label{eq:st:50a}
    \lim\limits_{h\to 0}\mathcal G_h(y_h,n_h)=&\mathcal G_{\mathrm{weak}}(y,n),\qquad\lim\limits_{h\to 0}\mathcal H_h(y_h,n_h)=&\mathcal H_{\mathrm{weak}}(y,n).
  \end{alignat}
  We only prove the statement for $\mathcal G_h$, since the argument for $\mathcal H_h$ is similar.
  By passing to a subsequence (not relabeled), we have $ (\nabla\!_hy_h,n_h)\to (R_y,n)$ a.e.~in $\Omega\times\OmegaT$, and thus, $|\frac{\nabla\!_hy_h^\top n_h}{|\nabla\!_hy_h^\top n_h|}-\hat n_0|^2\rho\to |R^\top_y n-\hat n_0|^2\rho$ a.e.
  By dominated convergence the latter holds strongly in $L^1(\OmegaT)$, and thus \eqref{eq:st:50a} follows.
  
  \step{2 -- Proof of part (a)}\smallskip
  
  We only present the proof for $\mathcal G_h$, since the argument for $\mathcal H_h$ is the same.  We start with an argument for the lower bound of part (a): Let $(y_h,n_h)\subseteq H^1(\Omega;\R^3)\times H^1(\OmegaT;\mathbb S^2)$ be a sequence converging to $(y,n)\in\mathcal A_2$ strongly in $L^2$. We need to show
  \begin{equation}\label{eq:st:50}
    \mathcal E(y,n)+\mathcal G_{\mathrm{weak}}(y,n)\leq\liminf\limits_{h\to 0}\big(\mathcal E_h(y_h,n_h)+\mathcal G_h(y_h,n_h)\big).
  \end{equation}
  By passing to a subsequence, we may assume  w.l.o.g.~that $\lim\limits_{h\to 0}\big(\mathcal E_h(y_h,n_h)+\mathcal G_h(y_h,n_h)\big)=\liminf\limits_{h\to 0}\big(\mathcal E_h(y_h,n_h)+\mathcal G_h(y_h,n_h)\big)<\infty$ and (in view of Theorem~\ref{T1} \ref{item:T1:compactness} and \eqref{eq:st:0}) that
  \begin{equation*}
    \det\nabla\!_hy_h>0\text{ a.e.~and }(y_h,\nabla\!_h y_h,n_h)\to (y,R_y,n)\text{ strongly in }L^2.
  \end{equation*}
  Combined with \eqref{eq:st:50a} and Theorem~\ref{T1} \ref{item:T1:lower:bound}, \eqref{eq:st:50} follows.

  Next, we prove the existence of a recovery sequence for $(y,n)\in\mathcal A_2$.
  To that end let $n_h\colonequals n$ and let $y_h$ denote the recovery sequence of Lemma~\ref{L:stronger}.
  In view of \eqref{eq:st:0}, also \eqref{eq:st:50c} holds.
  We conclude with help of \eqref{eq:st:50a} that $\mathcal E_h(y_h,n_h)+\mathcal G_h(y_h,n_h)\to\mathcal E(y,n)+\mathcal G_{\mathrm{weak}}(y,n)$.

  \step{3 -- Proof of part (b)}\smallskip
  
  We only present the proof for $\mathcal H_h$, since the argument for $\mathcal G_h$ is the same. We start with the lower bound.
  Let $(y_h,n_h)\subseteq H^1(\Omega;\R^3)\times H^1(\OmegaT;\mathbb S^2)$ be a sequence converging to $(y,n)\in\mathcal A_2$ strongly in $L^2$ and assume w.l.o.g.~that
$
  \liminf\limits_{h\to 0}\big(\mathcal E_h(y_h,n_h)+h^{-\beta}\mathcal H_h(y_h,n_h)\big)=\lim\limits_{h\to 0}\big(\mathcal E_h(y_h,n_h)+h^{-\beta}\mathcal H_h(y_h,n_h)\big)<\infty.
$
 Since 
$
 \mathcal E_h+h^{-\beta}\mathcal H_h\geq\mathcal E_h
$,
 Theorem~\ref{T1} \ref{item:T1:compactness} and \eqref{eq:st:0} yield \eqref{eq:st:50c} and thus \eqref{eq:st:50a}, which in combination with the bound
$
   \mathcal H_h\leq h^\beta(\mathcal E_h+h^{-\beta}\mathcal H_h)
$
 implies that $\mathcal H_{\mathrm{weak}}(y,n)=0$.
  Hence, $(y,n)$ satisfies the strong anchoring $(R_y^\top n\cdot\hat \nu)\bar\rho=0$ a.e.~in $S$.
  On the other hand, Theorem~\ref{T1} \ref{item:T1:lower:bound} yields  
  \begin{equation*}
    \liminf\limits_{h\to 0}\big(\mathcal E_h(y_h,n_h)+h^{-\beta}\mathcal H_h(y_h,n_h)\big)\geq \mathcal E(y,n).
  \end{equation*}
  Since $(y,n)$ satisfies the strong anchoring condition, the right-hand side is equal to $\mathcal E(y,n)+\mathcal H_{\mathrm{strong}}(y,n)$. This completes the argument for the lower bound.
  \medskip
  
  Next, we construct a recovery sequence for $(y,n)\in\mathcal A_2$ satisfying $\mathcal H_{\mathrm{strong}}(y,n)=0$.
  We start with approximating  $(y,n)$  by a configuration including a smooth isometry: Via \cite{Pakzad,HornungApprox} (or via Proposition \ref{P:approx} in the case of clamped, affine boundary conditions), we choose for each $k\in\N$ some $y^k\in H^2_\mathrm{ iso}(S;\R^3)\cap C^\infty(\overline S;\R^3)$ with
  \begin{equation}\label{eq:st:55}
    \|y^k-y\|_{H^2(S)}\leq \frac1k,
  \end{equation}
  and consider $n^k\colonequals R_{y^k} R_y^\top n$.
  Then $n^k\in H^1(S;\mathbb S^2)$ and $(R_{y^k}^\top n^k\cdot\hat\nu)\bar\rho=(R_y^\top n\cdot\hat\nu)\bar\rho=0$, and thus $\mathcal H_{\mathrm{strong}}(y^k,n^k)=0$. Furthermore, from \eqref{eq:st:55} we conclude that $n_k\to n$ strongly in $H^1(\OmegaT)$, and thus Lemma~\ref{L:cont} yields
  \begin{equation}\label{eq:st:55b}
    \lim\limits_{k\to\infty}|\mathcal E(y^k,n^k)-\mathcal E(y,n)|=0.
  \end{equation}
  Next, we construct the recovery sequence for $(y^k,n^k)$:
  As in the proof of Lemma~\ref{L:stronger}, we find a sequence $(y^k_h)_h\subseteq H^1(\Omega;\R^3)$ such that with $n_h^k\colonequals n^k$ we have
  \begin{align}\label{eq:st:56a}
    &\lim\limits_{h\to0}\mathcal E_h(y^k_h,n^k_h)=\mathcal E(y^k,n^k),\\\label{eq:st:56b}
    &(y^k_h,n^k_h)\to (y^k,n^k)\text{ in }L^2(\Omega)\times L^2(\OmegaT),\\\label{eq:st:56c}
    &\limsup\limits_{h\to 0}h^{-\beta/2}\|\nabla\!_hy^k_h-R_{y^k}\|_{L^\infty(\Omega)}=0.
  \end{align}
  From \eqref{eq:st:56c}, and since $\mathcal H_{\mathrm{strong}}(y^k,n^k)=0$, we infer
  \begin{equation*}
    h^{-\beta}\int_{\OmegaT}    \Big|\frac{(\nabla\!_hy_h^k)^\top n^k_h}{|(\nabla\!_hy_h^k)^\top n^k_h|}\cdot\hat\nu\Big|^2\rho\to 0\qquad(\text{as }h\rightarrow 0 ).
  \end{equation*}
  Thus, combined with \labelcref{eq:st:55,eq:st:56a,eq:st:56b,eq:st:55b}, we conclude that
  \begin{equation*}
    c_{k,h}\colonequals |\mathcal E_h(y^k_h,n^k_h)+h^{-\beta}\mathcal H_h(y^k_h,n^k_h)-(\mathcal E(y,n)+\mathcal H_{\mathrm{strong}}(y,n))|+\|y_h^k-y\|_{L^2}+\|n^k_h-n\|_{L^2}
  \end{equation*}
  satisfies $\limsup\limits_{k\to\infty}\limsup\limits_{h\to 0} c_{k,h}=0$. Hence, the sought for recovery sequence is obtained by extracting a suitable diagonal sequence.
\end{proof}

\subsubsection{Proofs of Lemmas~\ref{L:ext} \& \ref{L:weak:anchoring:with:boundary:slices}}
\label{Sec:Proof:215:217}

The proof of Lemma~\ref{L:ext} directly follows from (a) of the following statement:

\begin{lemma}[Estimates for the canonical extension]\label{L:trace}
	Let $S\subseteq \mathbb R^2$ be a bounded Lip\-schitz do\-main with outer unit normal $\nu_S$.	
	Assume \eqref{eq:gamma_n}.
  Then there exists $C>0$ and $\bar\sigma>0$ (only depending on $\Gamma_n$ and $S$) such that the following properties hold:
  \begin{enumerate}[(a)]
  \item The map $\Psi:\Gamma_n\times(0,\bar\sigma)\to \Gamma^{\bar\sigma}_n$, $(x',s)\mapsto x'-s\nu_S(x')$ is a diffeomorphism and the canonical extension $\mathrm{E}g$ of $g\in L^2(\Gamma_n;\mathbb S^2)$ defined in Lemma~\ref{L:ext} is well-defined.
  \item
  For all $g\in L^2(\Gamma_n;\mathbb S^2)$ and $0<\sigma\leq\bar\sigma$, we have
    \begin{equation*}
      \Big|\,\frac{1}{\sigma}\int_{S\cap\Gamma_n^\sigma}|\mathrm{ E}g|^2\,\mathrm dx'-\int_{\Gamma_n}|g|^2\,\mathrm d\mathcal H^1\,\Big|\leq C\sigma.
    \end{equation*}
  \item For all $0<\sigma\leq\bar\sigma$, $f\in H^1(S\cap\Gamma^\sigma_n;\mathbb S^2)$ and $g\in L^2(\Gamma_n;\mathbb S^2)$, we have
    \begin{equation*}
      \Big|\frac1\sigma\int_{S\cap\Gamma^\sigma_n}|f-\mathrm{ E}g|^2\,\mathrm dx'-\int_{\Gamma_n}|Tf-g|^2\,\mathrm d\mathcal H^1\Big|\\
      \leq\,C\Big(\sqrt\sigma\|\nabla' f\|_{L^2(S\cap\Gamma_n^\sigma)} + \sigma\Big),
    \end{equation*}
    where $T:H^1(\Gamma_n^\sigma)\to L^2(\Gamma_n)$ denotes the trace operator.
  \end{enumerate}
\end{lemma}
\begin{proof}
  We first note that (a) is a direct consequence of the Tubular Neighborhood Theorem (see e.g.,~\cite{oliva2004geometric}) and compactness of $\overline{\Gamma_n}$ in the $C^2$-manifold of \eqref{eq:gamma_n}.
  In the following, we assume w.l.o.g.~that $\Gamma_n$ admits a global parametrization,
  i.e., there exists an arc-length parametrized curve $\gamma\in C^2([a,b];\R^2)$ that homoeomorphically maps $(a,b)$ onto $\Gamma_n$.
  The general case can be reduced to this situation by the usual localization argument with help of a partition of unity and a finite atlas.
  In the following, $C>0$ denotes a constant that may change from line to line, but that can be chosen depending only on $\Gamma_n$ and $S$. 
  
  \step{1 -- Proof of (b)}\smallskip
  
  Using \eqref{eq:gamma_n} and the uniform inner cone condition of the Lipschitz domain $S$, it is possible to obtain $|\Gamma_n^\sigma\setminus S|\leq C\sigma^2$ for all $0<\sigma\leq\bar\sigma$ (for this we possibly have to decrease the constant $\bar\sigma$ of part (a) depending on the Lipschitz constant of the domain).
  In view of this and since $|\mathrm E g|=1$ a.e., for (b) we only need to prove
     \begin{equation} \label{eq:old:estimate:boundedness:E}
      \Big|\,\frac{1}{\sigma}\int_{\Gamma_n^\sigma}|\mathrm{ E}g|^2\,\mathrm dx'-\int_{\Gamma_n}|g|^2\,\mathrm d\mathcal H^1\,\Big|\leq C\sigma.
    \end{equation}
    We start our argument by setting $U\colonequals (a,b)\times(0,1)$ and for $0<\sigma\leq\bar\sigma$,
  \begin{equation*}
    \Phi_\sigma:U\to\Gamma_n^\sigma,\qquad \Phi_\sigma(t,s)\colonequals \gamma(t)-\sigma sR^\perp\gamma'(t),
  \end{equation*}
  where $R^\perp\in\SO 2$ denotes the unique rotation satisfying $\nu_S(\gamma(t))=R^\perp\gamma'(t)$ for all $t\in(a,b)$.
Note that $\Phi_\sigma(t,s)=\Psi(\gamma(t),\sigma s)$  for all $(t,s)\in U$.
Hence, by using the definition of the extension $\mathrm{ E}g$ we conclude that 
 \begin{equation} \label{eq:Identity:g:circ:Phi}
 	(\mathrm{ E}g)(\Phi_\sigma(t,s))=(g\circ\Phi_\sigma)(t,0)=(g\circ\gamma)(t).
 \end{equation}
 Also note that $\Phi_\sigma:U\rightarrow\Gamma_n^\sigma$ is a homeomorphism because of (a).
 Furthermore, $\gamma\in C^2$ yields $\Phi_\sigma\in C^1(U)$ and the Jacobian of $\Phi_\sigma$ satisfies $ |\det(J_{\Phi_\sigma}(t,s))|=\sigma|1-s\sigma\kappa(t)|$ where $\kappa(t)\colonequals R^\perp\gamma''(t)\cdot\gamma'(t)$.
 Now, by a change of variables,
 \begin{alignat*}{2}
   &\frac{1}{\sigma}\int_{\Gamma^\sigma_n}|\mathrm{ E}g|^2\,\mathrm dx' = \frac{1}{\sigma}\int_{U}|(\mathrm{ E}g)\circ\Phi_\sigma|^2|\det J_{\Phi_\sigma}|\,\mathrm ds\mathrm dt\\
   =&\,\int_U|(g\circ\gamma)(t)|^2|1-\sigma s\kappa(t)|\,\mathrm ds\mathrm dt=\int_{\Gamma_n}|g|^2\,\mathrm d\mathcal H^1+\int_U|g\circ\gamma|^2(|1-\sigma s\kappa(t)|-1)\,\mathrm ds\mathrm dt.
 \end{alignat*}
 Since $\|g\|_{L^\infty(\Gamma_n)}\leq 1$ and $\|\kappa\|_{L^\infty(\Gamma_n)}\leq C$, \eqref{eq:old:estimate:boundedness:E} holds and thus (b) follows.

 \step{2 -- Proof of (c)}\smallskip

 First, we extend $f$ to a function $\hat f\in H^1(\Gamma_n^\sigma;\R^3)$ such that $\hat f=f$ a.e.~on $S\cap\Gamma_n^\sigma$, $\|\nabla'\hat f\|_{L^2(\Gamma^\sigma_n)}\leq C\|\nabla' f\|_{L^2(\Gamma^\sigma_n\cap S)}$ and $\|\hat f\|_{L^\infty(\Gamma^\sigma_n)}\leq C\|f\|_{L^\infty(S\cap \Gamma^\sigma_n)}$.
 This is possible since $S\cap\Gamma^\sigma_n$ is a Lipschitz domain.
 In view of this, and since $|\Gamma_n^\sigma\setminus S|\leq C\sigma^2$, we see that for (c), it suffices to show
 \begin{equation}\label{eq:Gruenhorn}
     \Big|\frac1\sigma\int_{\Gamma^\sigma_n}|\hat f-\mathrm{ E}g|^2\,\mathrm dx'-\int_{\Gamma_n}|T\hat f-g|^2\,\mathrm d\mathcal H^1\Big|
     \leq\,C\Big(
     \sqrt\sigma\|\nabla' \hat f\|_{L^2(\Gamma_n^\sigma)}+\sigma\Big).
 \end{equation}
 To prove the latter, we note that
 \begin{equation}
 \begin{aligned}
  \frac1\sigma\int_{\Gamma^\sigma_n}|\hat f-\mathrm{ E}g|^2\,\mathrm dx'
  &=\int_U|(\hat f-\mathrm{ E}g)\circ\Phi_\sigma|^2|1- \sigma s\kappa(t)|\,\mathrm ds\mathrm dt,\\
  \int_{\Gamma_n}|T\hat f-g|^2\,\mathrm d\mathcal H^1 
  &= \int_U|((T\hat f-g)\circ\Phi_\sigma)(t,0)|^2\mathrm ds \mathrm dt.
 \end{aligned}
\end{equation}
Since $|(\hat f-\mathrm{ E}g)\circ\Phi_\sigma|^2\big| |1-\sigma s\kappa(t)|-1 \big|\leq C\sigma$ a.e.~in $U$, we conclude that
\begin{equation*}
\text{[LHS of \eqref{eq:Gruenhorn}]} \leq C\sigma
+
 \int_U\Big||(\hat f-\mathrm{ E}g)\circ\Phi_\sigma|^2 - |((T\hat f-g)\circ\Phi_\sigma)(t,0)|^2 \Big|\,\mathrm ds\mathrm dt.
\end{equation*}
We set $A(t,s):=(\hat f-\mathrm{ E}g)(\Phi_\sigma(t,s))$ and $B(t,s):=((T\hat f-g)(\Phi_\sigma(t,0))$.
Therefore, using the general estimate $\|A^2-B^2\|_{L^1}\leq(\|A\|_{L^2}+\|B\|_{L^2})\|A-B\|_{L^2}$ and $\|A\|_{L^\infty(U)} + \|B\|_{L^\infty(U)}<C$ (which holds thanks to the uniform boundedness of $\hat f$ and $g$), we obtain
\begin{align*}
\text{[LHS of \eqref{eq:Gruenhorn}]} 
&= C\sigma
+ C\big( \int_U |(\hat f \circ\Phi_\sigma)(t,s) - (T\hat f\circ\Phi_\sigma)(t,0)|^2\,\mathrm ds\mathrm dt \big)^{\frac12},
\end{align*}
where we also used \eqref{eq:Identity:g:circ:Phi} in form of the identity $(\mathrm{E}g\circ\Phi_\sigma)(t,s)-(g\circ\Phi_\sigma)(t,0)=0$.
 Moreover, since  $(\hat f\circ\Phi_\sigma)(t,s)-(T\hat f\circ\Phi_\sigma)(t,0)=\int_0^s\partial_{s'}(\hat f\circ\Phi_\sigma)(t,s')\,\mathrm ds'$ and $|\partial_{s'}(\hat f\circ\Phi_\sigma)|\leq C\sigma|\nabla'\hat f\circ\Phi_\sigma|$, we get
\begin{equation*}
\text{[LHS of \eqref{eq:Gruenhorn}]} 
\leq\,
C\sigma + 
C\sigma\big(\int_U\big|\big(\nabla'\hat f\circ\Phi_\sigma\big)(t,s')\big|^2\,\mathrm ds'\mathrm dt\big)^{\frac12}.
\end{equation*}
The claim \eqref{eq:Gruenhorn} then follows from $(1-C\sigma)\|\nabla'\hat f\circ\Phi_\sigma\|_{L^2(U)}^2\leq \sigma^{-1}\|\nabla'\hat f\|_{L^2(\Gamma_n^\sigma)}^2$.
\end{proof}

\begin{proof}[Proof of Lemma~\ref{L:weak:anchoring:with:boundary:slices}]
  We only present the argument for $\mathcal G^\sigma$, since $\mathcal H^\sigma$ can be treated similarly.
  Throughout the proof, $C$ denotes a positive constant that might change from line to line, but that can be chosen independently of $\sigma$.

  \step{1 -- Proof of (a), lower bound}\smallskip
  
   Let $\beta=0$ and consider a sequence $(y_\sigma,n_\sigma)$   that converges to some $(y,n)$ strongly in $L^2(S)$.
  W.l.o.g.~we may assume that 
$
  \liminf\limits_{\sigma\to 0}(\mathcal E(y_\sigma,n_\sigma)+\mathcal G^\sigma(y_\sigma,n_\sigma))=\limsup\limits_{\sigma\to 0}(\mathcal E(y_\sigma,n_\sigma)+\mathcal G^\sigma(y_\sigma,n_\sigma))<\infty
$.
  Thus, Lemma~\ref{L:cont} implies that $(y_\sigma,n_\sigma),(y,n)\in \mathcal A_2$, $y_\sigma\wto y$ weakly in $H^2(S)$, $n_\sigma\wto n$ weakly in $H^1(S)$, and $\liminf_{\sigma\to 0}\mathcal E(y_\sigma,n_\sigma)\geq\mathcal E(y,n)$. It thus remains to prove
  \begin{equation}\label{eq:st:90}
    \liminf\limits_{\sigma\to 0}\frac1\sigma\int_{S\cap\Gamma^\sigma_n}|R_{y_\sigma}^\top n_\sigma-\mathrm{ E}\hat n_{BC}|^2\,\mathrm dx'
    \geq 
    \int_{\Gamma_n}|T(R_y^\top n)-\hat n_{BC}|^2\,\mathrm d\mathcal H^1,
  \end{equation}
  where $T:H^1(S;\R^3)\to L^2(\Gamma_n;\R^3)$ denotes the trace operator.
  For convenience set $f_\sigma\colonequals R_{y_\sigma}^\top n_{\sigma}$.
  We note that $(f_\sigma)$ is a bounded sequence in $H^1(S;\mathbb S^2)$.
  Furthermore, from Lemma~\ref{L:trace} (c), we get
  \begin{equation}
    \label{eq:st:91}
    \Big|\frac1\sigma\int_{S\cap\Gamma^\sigma_n}|f_\sigma-\mathrm{ E}\hat n_{BC}|^2\,\mathrm dx'-\int_{\Gamma_n}|Tf_\sigma-\hat n_{BC}|^2\,\mathrm d\mathcal H^1\Big|\,\leq\,
    C\sigma+C\sqrt{\sigma}\|\nabla'f_\sigma\|_{L^2(S\cap\Gamma_n^\sigma)}.
  \end{equation}
  We infer that
  \begin{equation*}
    \liminf\limits_{\sigma\to 0}\frac{1}{\sigma}\int_{S\cap\Gamma^\sigma_n}|R_{y_\sigma}^\top n_\sigma-\mathrm{ E}\hat n_{BC}|^2\,\mathrm dx'
    =
    \liminf\limits_{\sigma\to0 }\int_{\Gamma_n}|Tf_\sigma-\hat n_{BC}|^2\,\mathrm d\mathcal H^1.
  \end{equation*}
  Now, statement \eqref{eq:st:90} follows from the fact that $f_\sigma\wto R_y^\top n$ weakly in $H^1(S;\R^3)$, and thus $Tf_\sigma-\hat n_{BC}\wto T(R_y^\top n)-\hat n_{BC}$ weakly in $L^2(\Gamma_n;\R^3)$, and the weak lower semicontinuity of $\|\cdot\|_{L^2(\Gamma_n)}$.

  \step{2 -- Proof of (b), lower bound}\smallskip
  
  Let $0<\beta<\frac12$ and consider a sequence $(y_\sigma,n_\sigma)$ that converges to some $(y,n)$ strongly in $L^2(S)$.
  As in Step~1 we may assume w.l.o.g.~that
  $
    \liminf\limits_{\sigma\to 0}(\mathcal E(y_\sigma,n_\sigma)+\sigma^{-\beta}\mathcal G^\sigma(y_\sigma,n_\sigma))=\limsup\limits_{\sigma\to 0}(\mathcal E(y_\sigma,n_\sigma)+\sigma^{-\beta}\mathcal G^\sigma(y_\sigma,n_\sigma))<\infty
  $,
  as well as
  $(y_\sigma,n_\sigma),(y,n)\in \mathcal A_2$, $y_\sigma\wto y$ weakly in $H^2(S)$, $n_\sigma\wto n$ weakly in $H^1(S)$, and
$
   \liminf_{\sigma\to 0}\mathcal E(y_\sigma,n_\sigma)\geq\mathcal E(y,n)
$.
Thanks to \eqref{eq:st:90} and $\frac1\sigma\int_{S\cap \Gamma^\sigma_n}|R_{y_\sigma}^\top n_\sigma-\mathrm{ E}\hat n_{BC}|^2= \sigma^\beta\sigma^{-\beta}\mathcal G^\sigma(y_\sigma,n_\sigma)\rightarrow 0 $, we conclude that the strong anchoring $R_y^\top n=\hat n_{BC}$ on $\Gamma_n$ is satisfied, and thus $\mathcal E(y,n)+\mathcal G_{\mathrm{strong},\Gamma_n}(y,n)=\mathcal E(y,n)\leq \liminf_{\sigma\to 0}\mathcal E(y_\sigma,n_\sigma)\leq \liminf_{\sigma\to 0}(\mathcal E(y_\sigma,n_\sigma)+\sigma^{-\beta}\mathcal G^\sigma(y_\sigma,n_\sigma))$.

  \step{3 -- Proof of (a), upper bound}\smallskip
  
  It suffices to construct a recovery sequence for $(y,n)$ with $\mathcal E (y,n)+\mathcal G_{\mathrm{weak},\Gamma_n}(y,n)<\infty$, which implies that $(y,n)\in\mathcal A_2$. We claim that the constant sequence is a recovery sequence, i.e., we need to show that $
 \frac1\sigma\int_{S\cap\Gamma_n^\sigma}|R^\top_yn-\mathrm{ E}\hat n_{BC}|^2\,\mathrm dx'\to\int_{\Gamma_n}|T(R^\top_yn)-\hat n_{BC}|^2\,\mathrm d\mathcal H^1
$.
  The latter follows from Lemma \ref{L:trace} (c).
  
  \step{4 -- Proof of (b), upper bound} \smallskip
  
  Let $(y,n)$ with $\mathcal E(y,n)+\mathcal G_{\mathrm{strong},\Gamma_n}(y,n)<\infty$. The latter implies that $(y,n)\in\mathcal A_2$ and $R^\top_yn=\hat n_{BC}$ on $\Gamma_n$. We claim that the constant sequence is a recovery sequence.
  Indeed, from Lemma \ref{L:trace} (c) applied with $f:=R^\top_yn$ and thanks to $Tf=\hat n_{BC}$ on $\Gamma_n$, we learn that
  \begin{align*}
    \sigma^{-(1+\beta)}\int_{S\cap\Gamma_n^\sigma}|R^\top_yn-\mathrm{ E}\hat n_{BC}|^2\,\mathrm dx'
    &\leq \sigma^{-\beta}\left(\int_{\Gamma_n}|Tf-\hat n_{BC}|^2\,\mathrm d\mathcal H^1+C(\sqrt{\sigma}\|\nabla'f\|_{L^2(S\cap\Gamma_n^\sigma)}+\sigma) \right).
  \end{align*}
  Since $\frac12-\beta>0$, we conclude that $\sigma^{-(1+\beta)}\int_{S\cap\Gamma_n^\sigma}|R^\top_yn-\mathrm{ E}\hat n_{BC}|^2\,\mathrm dx'\to 0$ and
  \begin{equation*}
    \lim\limits_{\sigma\to 0}(\mathcal E(y,n)+\sigma^{-\beta}\mathcal G^\sigma(y,n))=\mathcal E(y,n)=\mathcal E(y,n)+\mathcal G_{\mathrm{strong},\Gamma_n}(y,n).
  \end{equation*}
\end{proof}
\subsection{Approximation of the nonlinear strain: proof of Proposition \ref{P:general_construction}} \label{Sec:Proof:P:general_construction}
The structure of the proof of Proposition \ref{P:general_construction} is as follows.
\begin{enumerate}[(i)]
\item
  For each smooth isometry $y:S\to\R^3$, the domain $S$ can be decomposed into a flat part $\{\II_y=0\}$ and a non-flat part $\{\II_y\neq 0\}$, and, as we shall see, connected components $V$ of $\{\II_y\neq 0\}$ and $\interieur\{\II_y=0\}$ have a special geometry:
  The relative boundary $\partial V\cap S$ is a disjoint union of line segments whose endpoints are contained in $\partial S$, and $y$ is affine on $\partial V\cap S$.
  Furthermore, in the non-flat case, $\partial V\cap S$ consists of (at most) two such line segments.
  Furthermore, we see that $\partial\{\II_y\neq 0\}\cap S$ is a null set, see Lemma~\ref{L:BoundaryMeasureVanishes}.
  This allows us to construct the recovery sequence on each connected component independently.

\item For each connected component $V\subseteq\{\II_y\neq 0\}$, we construct the recovery sequence based on ideas in \cite{SCHMIDT07, HNV}.
The construction invokes a solution to the system $\nabla'_{\sym}g+\alpha\II_y= M$.
The existence of such a solution has already been shown in \cite{SCHMIDT07}. In order to recover the prescribed boundary conditions, we shall slightly upgrade that existence result by showing that solutions exist that vanish on the relative boundary $\partial V\cap S$ and the line segments $\mathcal L_i$ that appear in our boundary condition.
\item The construction on a connected component $V\subseteq\interieur\{\II_y=0\}$ is similar to \cite{DPG20} and relies on an approximation result for the Monge-Amp\'ere equation introduced in \cite{LP17}.
Nevertheless, we shall give a self-contained presentation of the construction, since we need to take care of the additional boundary conditions.
\end{enumerate}
The above discussion is made precise in the following two lemmas.
\begin{lemma}\label{L:BoundaryMeasureVanishes}
    Let $S$ satisfy \eqref{ass:domain}. Let $y\in H^2_{\rm iso}(S;\R^3)\cap C^\infty(\overline S;\mathbb{R}^3)$. Then $\partial \{\II_y=0\}\cap S=\partial\{\II_y\neq 0\}\cap S$ is a null set in $\R^2$.
\end{lemma}
\textit{(See Section \ref{Sec:Proof:L:BoundaryMeasureVanishes} for the proof.)}\medskip
\begin{lemma}[Recovery sequence on curved and flat connected components]\label{L:recov_general}
  Let $S$ satisfy \eqref{ass:domain}. 
  Let $y\in H^2_{\iso}(S;\R^3)\cap C^\infty(\overline S;\R^3)$, $M\in L^2(S;\R^{2\times 2}_{\sym})$, $d\in L^2(\Omega;\R^3)$, $\delta\geq 0$, and $V\subseteq S$. Suppose that one of the following cases holds:
  \begin{itemize}
  \item (curved case). $V$ is a connected component of $\{\II_y\neq 0\}$.
  \item (flat case). $V$ is a connected component of $\interieur\{\II_y=0\}$ and $M+\delta I_{2\times 2}$ is positive semi-definite a.e.~in $V$.
  \end{itemize}
  Let $0<\beta<\frac12$. Then there exists a sequence $(y_h)\subseteq C^{\infty}(\overline V\times[-\tfrac12,\tfrac12];\R^3)$ that satisfies
  \begin{equation}\label{eq:recov:onescalea}
    \begin{aligned}
      &y_h\,\to\, y\qquad\text{uniformly in } V\times(-\tfrac12,\tfrac12),\\
      &E_h(y_h)\,\to\,\iota(x_3 \II_y+M)+\sym(d\otimes e_3)	\quad\text{strongly in }L^2(V\times(-\tfrac12,\tfrac12)),\\
    &\limsup_{h\to 0} h^{-\beta}\|\nabla\!_hy_h-R_y\|_{L^\infty(V\times(-\frac12,\frac12))}=0.
  \end{aligned}
  \end{equation}
  Furthermore, we have
  \begin{equation}
    \label{eq:BC3daa}
    y_h=(1-h\delta)y+hx_3b_y\qquad\text{ in a uniform neighborhood of }(\partial V\cap S)\times(-\frac12,\frac12).
  \end{equation}
  Additionally, if Assumption~\ref{ass:BC} is satisfied and $y\in\mathcal A_{BC}$, then
  \begin{equation}
    \label{eq:BC3dac}
    \begin{aligned}
    &y_h=(1-h\delta)y+hx_3b_y\\
    &\text{in a uniform neighborhood of  }\big(\partial V\cap (S\cup\overline{\mathcal L_1}\cup\cdots\cup\overline{\mathcal L_{k_{BC}}})\Big)\times(-\frac12,\frac12).
    \end{aligned}
  \end{equation}
\end{lemma}
\textit{(See Sections \ref{Sec:Proof:L:recov_general_curved} and \ref{Sec:Proof:L:recov_general_flat} for the proof.)}\medskip
\smallskip

With these results at hand we are in position to prove Proposition~\ref{P:general_construction}:
\begin{proof}[Proof of Proposition~\ref{P:general_construction}]
  Let $(V_i)_{i\in I}$ denote an enumeration of the connected components of the two sets $\{\II_y\neq 0\}$ and $\interieur\{\II_y=0\}$,
  and denote by $(y^i_h)_h\subseteq C^\infty(\overline{V_i}\times[-\tfrac12,\tfrac12];\R^3)$ the sequences constructed in Lemma~\ref{L:recov_general}, respectively. Furthermore, we set
  \begin{equation*}
      \bar y_h(x',x_3)\colonequals (1-h\delta)y(x')+hx_3b_y(x').
    \end{equation*}
    Thanks to Lemma~\ref{L:BoundaryMeasureVanishes} we have $|S\setminus\bigcup_{i\in I}V_i|=0$ and thus for any $n\in\N$ there exists $N_n\in \N$ such that
    \begin{equation}\label{eq:def:choicedelta}
      \int_{\Omega\setminus\Omega_n}|\iota(M)+\sym(d\otimes e_3)|^2\,\mathrm dx'<\frac1n\text{ where }\Omega_n\colonequals \bigcup_{i=1}^{N_n}V_i\times\left(-\tfrac12,\tfrac12\right).
    \end{equation}
    We define
    \begin{equation*}
      y_{n,h}(x)\colonequals 
      \begin{cases}
        y^i_h(x)&\text{if }x\in V_i\times(-\tfrac12,\tfrac12)\text{ for some }i\leq N_n,\\
        \bar y_h(x)&\text{else}.
      \end{cases}
    \end{equation*}
    Thanks to \eqref{eq:BC3daa}, $y_{n,h}\in C^\infty(\overline{\Omega};\R^3)$. Furthermore, if $y\in\mathcal A_{BC}$, then $y_{n,h}\in\mathcal A_{BC,\delta,h}$  thanks to \eqref{eq:BC3dac}.
    Also set
    \begin{equation*}
      E_n(x)\colonequals 
      \begin{cases}
        \iota(x_3\II_y+M)+\sym(d\otimes e_3)&x\in \Omega_n,\\
        \iota(x_3\II_y)&\text{else.}
      \end{cases}
    \end{equation*}
    Then, direct calculations show
    \begin{equation*}
      \begin{aligned}
        \lim_{h\rightarrow 0} 
        \|E_h(y_{n,h})- E_n\|_{L^2(\Omega)}
        +
  		h^{-\beta}\|\nabla\!_hy_{n,h}-R_y\|_{L^\infty(\Omega)}
  		+
        \|y_{n,h}-y\|_{L^\infty(\Omega)}
        =0.
      \end{aligned}
    \end{equation*}
    In view of \eqref{eq:def:choicedelta}, we have $\|E_n-\big(\iota(x_3\II_y+M)+\sym(d\otimes e_3)\big)\|_{L^2(\Omega)}\to 0$ as $n\to\infty$, and thus, the error term
    \begin{equation*}
      C_{n,h}\colonequals \|E_h(y_{n,h})-\big(\iota(x_3\II_y+M)+\sym(d\otimes e_3)\big)\|_{L^2(\Omega)}+h^{-\beta}\|\nabla\!_hy_{n,h}-R_y\|_{L^\infty(\Omega)}+\|y_{n,h}-y\|_{L^\infty(\Omega)}
    \end{equation*}
    satisfies $\limsup_{n\to\infty}\limsup_{h\to 0}C_{n,h}=0$. This allows us to pass to a diagonal sequence $h\mapsto n_h$ such that $\lim\limits_{h\to 0}C_{n_h,h}=0$ and we conclude that $y_h\colonequals y_{n_h,h}$ defines the sought for sequence.
\end{proof}

\subsubsection{Auxiliary results for isometries}
For the proof of Lemma~\ref{L:BoundaryMeasureVanishes} and the recovery sequence construction in the curved case, we appeal to standard properties of isometric immersions, which we recall in the following. We refer to \cite{Pakzad,HornungApprox,HornungFine,mullerPakzad2005regularity} for details. Since $H^2$-regular isometric immersions are continuously differentiable (see \cite[Proposition~5]{mullerPakzad2005regularity}) the following definition is meaningful:
\begin{definition}
\label{DefinitionBodiesAndArms}
Let $S\subseteq \mathbb R^2$ be open.
For $y\in H^2_\mathrm{iso} (S;\mathbb{R}^3)$ we denote by
\begin{equation*}
  C_{\nabla'y}\colonequals  \left\lbrace
    x'\in S: \nabla' y \text{ is constant in some neighborhood of } x'\right\rbrace
\end{equation*}
the \emph{flat part} of $y$.
Furthermore, we denote by $\widehat C_{\nabla'y}$ the union of all connected components $U\subseteq C_{\nabla'y}$ with $\partial U\cap S$ consisting of more than two connected components.
\end{definition}
An important observation is that at each point $x'$ in the non-flat part $S\setminus C_{\nabla'y}$, the isometry $y$ is affine on a unique line segment that contains $x'$ and whose endpoints are contained in $\partial S$.
In fact, this property extends to the larger set $S\setminus \widehat C_{\nabla'y}$ as the following lemma shows.
For arbitrary $x'\in S$ and $N\in \mathbb S^1$, we use the notion of $[x';N]$ as the connected component of $S\cap\{x'+sN:\, s\in\mathbb{R}\}$ that contains $x'$ itself.

\begin{lemma}[\cite{Pakzad},\cite{HornungApprox}]
  \label{LemmaExistenceAsymptoticsField}
  Let $S$ be a bounded Lipschitz domain.
  Let $y\in  H^2_\mathrm{iso} (S;\mathbb{R}^3)$.
  Then, there is a locally Lipschitz continuous vector field $N:S\setminus \widehat C_{\nabla' y}\rightarrow \mathbb S^1$ such that for all $x_1',x_2'\in S\setminus \widehat C_{\nabla' y}$,
  \begin{align}
    \label{EqFourteenA}
    &\nabla' y \text{ is constant on } [x_1';N(x_1')],
    \\
    \label{EqFourteenB}
    &[x_1';N(x_1')]\cap [x_2';N(x_2')]\neq \emptyset \quad\Rightarrow\quad [x_1';N(x_1')] = [x_2';N(x_2')].
  \end{align}
  Moreover, if $x_1'\in S\setminus C_{\nabla' y}$, then property \eqref{EqFourteenA} defines the value $N(x_1')$ uniquely up to a sign. We call $N$ a \emph{field of asymptotic directions} of $y$.
\end{lemma}
Based on $N$ one can introduce natural parametrizations of $y$:
\begin{definition}[Patch]
  \label{DefinitionPatch}
  Let $S$ be a bounded Lipschitz domain.
  Let $y\in  H^2_\mathrm{iso}(S;\mathbb{R}^3) $ and $N$ as in Lemma \ref{LemmaExistenceAsymptoticsField}.
A curve $$\Gamma\in W^{2,\infty}([0,T];{S\setminus\widehat C_{\nabla'y}}),\qquad T>0,$$
  is called a \emph{ruling curve} of $y$, if for all $t,t'\in [0,T]$,
  \begin{align}
    & \label{eq:Def:ODE:ruling:curve}
      \Gamma'(t) = R^\perp N(\Gamma(t)),\qquad 0<\Gamma'(t)\cdot \Gamma'(t'),\qquad R^\perp\colonequals 
  \begin{pmatrix}
    0 &-1\\
    1&0
  \end{pmatrix}.
  \end{align}
  A set $P\subseteq S\setminus \widehat C_{\nabla'y}$ is called a \emph{patch} (with respect to $y$ and $N$), if there exists a ruling curve $\Gamma$ such that
  \begin{equation}\label{eq:patchdom}
    P = \bigcup_{t\in(0,T)}[\Gamma(t);(N\circ\Gamma)(t)].
  \end{equation}
\end{definition}
\begin{remark}[Line of curvature parametrization in the smooth case]\label{R:linecurvparam}
  On patches the line of cur\-va\-ture para\-metri\-zat\-ion yields a simple representation of the surface and its second fundamental form, e.g.,~see \cite{Pakzad,SCHMIDT07}, and in particular \cite[Proposition 1]{HornungApprox}.
  In the case of a smooth isometry (which is the situation that we consider in this section) the following additional properties hold:
  Let $y\in H^2_{\mathrm{iso}}(S;\mathbb R^3)\cap C^\infty(\overline{S};\R^3)$ with $N$ as in Lemma \ref{LemmaExistenceAsymptoticsField}.
  Then we have $C_{\nabla'y} = \interieur(\{\II_y= 0\})$ and $S\setminus C_{\nabla'y}=S\cap\overline{\{\II_y\neq 0\}}$ by continuity of $\II_y$.
  Moreover, we note that $N$ is smooth on $\{\II_y\neq 0\}$, as follows from $\II_y\,N=0$ (cf.~\ref{EqFourteenA}) and the smoothness of $\II_y$.

  Now, let $P$ be a patch for $y$ with ruling curve $\Gamma$ such that $\Gamma([0,T])\subseteq \{\II_y\neq 0\}$.
  Then $\Gamma$ is smooth on the closed interval $[0,T]$ as a consequence of \eqref{eq:Def:ODE:ruling:curve} and the smoothness of $N$.
  If $S$ is convex, then the set $S_\Gamma\colonequals \{(t,s)\in (0,T)\times\R\,:\,\Gamma(t)+s N(\Gamma(t))\in S\}$ is open and convex (see \cite[Lemma 3.6]{Pakzad}), and $S_\Gamma\to P$, $(t,s)\mapsto \Gamma(t)+sN(\Gamma(t))$ defines a smooth diffeomorphism.
  Furthermore, for all $(t,s)\in[0,T]\times\R$ with $\Gamma(t)+sN(\Gamma(t))\in S$ we have $1-s\kappa(t)>0$, and
  \begin{align}\label{eq:st:1121}
    \II_y(\Gamma(t)+sN(\Gamma(t)))=\,&\,\frac{\kappa_{\mathrm{n}}(t)}{1-s\kappa(t)}\Gamma'(t)\otimes\Gamma'(t),
  \end{align}
  where $\kappa\colonequals \Gamma''\cdot (N\circ \Gamma)$ and $\kappa_{\mathrm{n}}:[0,T]\to\R$ denotes a smooth function related to the curvature of $y\circ\Gamma$.
\end{remark}
\subsubsection{Proof of Lemma~\ref{L:BoundaryMeasureVanishes}}\label{Sec:Proof:L:BoundaryMeasureVanishes}
  By continuity of $\II_y$, the set $\{\II_y\neq 0\}$ is open; thus, it is the union of at most countably many connected components $\{V_i\}_{i\in I}$ of $\{\II_y\neq 0\}$.
  In view of $\partial\{\II_y\neq 0\}\cap S\subseteq \bigcup_{i\in I}(\partial V_i\cap S)$, it suffices to show that $\partial V_i\cap S$  is a null-set.
  To conclude the latter, we show that $\partial V_i\cap S$ consists of at most two line segments.
  (Note that this is only true since $S$ is convex.)
  The latter follows from the following two properties:
  \begin{enumerate}[(i)]
  \item The connected components of $\partial V_i\cap S$ are line segments of the form $[x';N(x')]$.
  \item There are at most two disjoint line segments of the form $[x';N(x')]$ contained in $\partial V_i\cap S$.
  \end{enumerate}
  Property (ii) is true, since otherwise $V_i$ contains a small triangle on which $y$ is affine -- in contradiction to $\II_y\neq 0$ on $V_i$.
  To prove (i), let $x'\in\partial V_i\cap S$. We show that $[x';N(x')]\subseteq\partial V_i\cap S$.
  We first note that $\II_y(x')=0$ and thus $\II_y=0$ on $[x';N(x')]$ (by \eqref{eq:st:1121} and a contradiction argument).
  Hence, $[x';N(x')]\cap V_i=\emptyset$.
  To conclude $[x';N(x')]\subseteq\partial V_i\cap S$, it suffices to show that any $\bar x'\in[x';N(x')]$ is the limit of a sequence in $V_i$.
  To that end, choose $s\in\R$ with $\bar x'=x'+sN(x')$ and consider the sequence $\bar x_n'\colonequals x_n'+sN(x_n')$ where $(x_n')\subseteq V_i$ with $x_n'\to x'$.
  By continuity of $N$, we have $\bar x_n'\to \bar x'$.
  Furthermore, since $S$ is open and $\bar x'\in S$, we have $\bar x_n'\in S$ for all $n$ large enough, and thus $\bar x_n'\in[x_n';N(x_n')]$.
  With \eqref{eq:st:1121}, we conclude that $\II_y(\bar x_n')\neq 0$ and hence $\bar x_n'\in V_i$.\qed

  \subsubsection{Proof  of  Lemma~\ref{L:recov_general} - curved case}\label{Sec:Proof:L:recov_general_curved}
We start with a lemma showing that connected components of $\{\II_y\neq 0\}$ can be conveniently covered (up to a small set) by a finite number of patches whose boundaries do not intersect the boundary segments $\mathcal L_1,\ldots,\mathcal L_{k_{BC}}$:
\begin{lemma}[Approximate covering by patches]
\label{LemmaFinePropertiesOfSmoothisometries}
    Let $S$ satisfy \eqref{ass:domain}. Let $y\in H^2_{\rm iso}(S;\R^3)\cap C^\infty(\overline S;\mathbb{R}^3)$, let $N$ denote the field of asymptotic directions from Lemma \ref{LemmaExistenceAsymptoticsField}, let $V$ be a connected component of $\{\II_y\neq 0\}$, and let $\varepsilon>0$.
Then there exist finitely many, pair-wise disjoint patches $P_1,\ldots,P_k\subseteq V$ associated with ruling curves $\Gamma_i$ such that $|V\setminus (P_1\cup\ldots\cup P_k)|<\varepsilon$, and $\Gamma_i:[0,T_i]\to V$ for all $i=1,\ldots,k$.
If additionally Assumption~\ref{ass:BC} is satisfied and $y\in\mathcal A_{BC}$, then
\begin{equation*}
  \overline{\mathcal L_j}\cap \partial P_i=\emptyset\qquad\text{for all }j=1,\ldots,k_{BC}\text{ and }i=1,\ldots,k.
\end{equation*}
\end{lemma}
\textit{(We postpone the proof to the end of this section.)}\medskip

For the construction of the recovery sequence on a single patch $P\subseteq\{\II_y\neq 0\}$, we need to solve the equation
\begin{equation}\label{eq:monge}
  \nabla'_{\sym}g+\alpha\II_y= M\qquad\text{in }P,
\end{equation}
such that $\alpha$ and $g$ vanish close to $\partial P\cap S$ provided that $M$ is compactly supported in $P$.
This can be achieved by modifying \cite[Lemma~3.3]{SCHMIDT07}:
\begin{lemma}[Modification of \mbox{\cite[Lemma~3.3]{SCHMIDT07}}]\label{L:monge1}
  Let $S$ satisfy \eqref{ass:domain}. 
  Let $y\in H^2_{\iso}(S;\R^3)\cap C^\infty(\overline S;\R^3)$ and $M\in C^\infty_c(P;\R^{2\times 2}_{\sym})$.
  Suppose that $P\subseteq\{\II_y\neq 0\}$ is a patch with its ruling curve $\Gamma$ satisfying $\Gamma:[0,T]\to \{\II_y\neq 0\}$.
  Then there exist $\alpha\in C^\infty(\overline P;\R)$ and $g\in C^\infty(\overline P;\R^2)$ solving \eqref{eq:monge} such that $g=\alpha=0$ in a uniform neighborhood of $\partial P\cap S$.
\end{lemma}
\textit{(We postpone the proof to the end of this section.)}\medskip
We are now in a position to present the construction of the recovery sequence.
\begin{proof}[Proof of Lemma \ref{L:recov_general} -- curved case]
  \step{1 -- Construction of a displacement field on a single patch}\smallskip

  Let $P\subseteq V$ be a patch such that the endpoints of the associated ruling curve are contained in $V$.
  Let $M_P\in C^\infty_c(P;\R^{2\times 2})$.
  We claim that there exists a sequence $w_h\in C^\infty(\overline P\times [-\frac12,\frac12];\R^3)$ that vanishes in a uniform neighborhood of $(\partial P\cap S)\cap(-\frac12,\frac12)$ such that
  \begin{equation}\label{eq:recov_general:1}
    \limsup\limits_{h\to 0}\|\nabla\!_hw_h\|_{L^\infty}<\infty,
    \qquad
    \limsup\limits_{h\to 0}h^{-1}\|\sym\big(R_y^\top\nabla\!_hw_h\big)-\iota (M_P)\|_{L^\infty}<\infty.
  \end{equation}
  For the argument let $g,\alpha\in C^\infty(\overline P)$ be a solution of \eqref{eq:monge} (with $M$ replaced by $M_P$) such that $g,\alpha=0$ in a uniform neighborhood of $\partial P\cap S$.
  The existence of $g,\alpha$ follows from Lemma~\ref{L:monge1}.
  As in \cite{HNV} we consider
  \begin{equation*}
    w_h\colonequals V+hx_3\mu,\qquad V=R_y
    \begin{pmatrix}
      g\\\alpha
    \end{pmatrix},\qquad \mu\colonequals (I_{3\times 3}-b_y\otimes b_y)(\partial_1V\wedge\partial_2y+\partial_1y\wedge\partial_2V),
  \end{equation*}
  We note that
  \begin{equation*}
    \sym\Big((\nabla'y)^\top(\nabla'V)\Big)=M_P,\qquad \sym
    \begin{pmatrix}
      0_{2\times 2}&(\nabla'y)^\top\mu\\
      b_y^\top\nabla'V&0
    \end{pmatrix}=0.
  \end{equation*}
  Indeed, the first identity follows from \eqref{eq:monge} by a direct calculation and the second identity is a consequence of the definition of $\mu$.
  Now, \eqref{eq:recov_general:1} follows from the smoothness of $R_y$,$b_y$, $g$ and $\alpha$.

  \step{2 -- Construction of the recovery sequence} \smallskip

  It suffices to provide a construction in the case $d\in C^\infty_c(V;\R^3)$ and $M\in C^\infty_c(V;\R^{2\times 2}_{\sym})$, since then the general case follows by a diagonal sequence argument. Furthermore, we only present the argument in the case of prescribed boundary conditions $y\in\mathcal A_{BC}$ (since this adds an additional layer of difficulty).
  
  Let $n\in\N$. By Lemma \ref{LemmaFinePropertiesOfSmoothisometries} there exists $k_n\in\N$, and pair-wise disjoint patches $P_1,\ldots,P_{k_n}\subseteq V$, whose ruling curves have endpoints in $V$, such that $|V\setminus ({P_1}\cup\ldots\cup {P_{k_n}})|<\frac1n$.
  Furthermore, there exists $M_n\in C^\infty_c(P_1\cup\cdots\cup P_{k_n})$ such that $\|M-(M_n-\delta I_{2\times 2})\|_{L^2(V)}^2\leq \frac{2}{n}\|M+\delta I_{2\times 2}\|_{L^\infty(V)}^2$.
  By Step~1 (applied to each patch $P_i$ with $M_P:=M_n\vert_{P_i}$), there exists a sequence $(w_{n,h})_h\subseteq C^\infty(\overline V\times[-\frac12,\frac12];\R^3)$ such that \eqref{eq:recov_general:1} is satisfied on each $P_i\times(-\frac12,\frac12)$ and such that $w_{n,h}=0$ in a uniform neighborhood of $(V\setminus (P_1\cup\ldots\cup P_{k_n}))\times(-\frac12,\frac12)$.
  Thus, $w_{n,h}$ also vanishes in a uniform neighborhood of $\big(\partial V\cap (S\cup \overline{\mathcal L_1}\cup\cdots\cup \overline{\mathcal L_{k_{BC}}})\Big)\times(-\frac12,\frac12)$; for the latter we used that $\overline{\mathcal L_j}\cap \partial P_i=\emptyset$, cf.~Lemma~\ref{LemmaFinePropertiesOfSmoothisometries}.
  We consider the sequence
  \begin{equation*}
    y_{n,h}\colonequals (1-h\delta)y+hx_3 b_y+h^2R_y\int_0^{x_3}d(\cdot,t)\,\mathrm dt+hw_{n,h}.
  \end{equation*}
  It is a standard Kirchhoff-Love ansatz modified by the displacement $w_{n,h}$ constructed via Step~1. By construction $y_{n,h}$ satisfies the required boundary conditions, i.e., we have
  $
    y_{n,h}=(1-h\delta)y+hx_3 b_y
  $
  in a uniform neighborhood of $\big(\partial V\cap (S\cup \overline{\mathcal L_1}\cup\cdots\cup \overline{\mathcal L_{k_{BC}}})\Big)\times(-\frac12,\frac12)$. A direct calculation shows that
  \begin{align*}
    \limsup\limits_{h\to 0}\Big(
    &\|y_{n,h}-y\|_{L^\infty(V\times(-\frac12,\frac12))}
    +
    h^{-\beta}\|\nabla\!_hy_{n,h}-R_y\|_{L^\infty(V\times(-\frac12,\frac12))}
    \\
    &\; 
    +
    \|E_h(y_{n,h})-\big(\iota(x_3\II+M_n-\delta I_{2\times 2})+\sym(d\otimes e_3)\big)\|_{L^2(V\times(-\frac12,\frac12))}
    \Big)
    =0.
  \end{align*}
  Since $M_n-\delta I_{2\times 2}\to M$ in $L^2(V)$, we obtain the sought for sequence by passing to a suitable diagonal sequence $y_h\colonequals y_{n_h,h}$.
\end{proof}
\begin{proof}[Proof of Lemma~\ref{LemmaFinePropertiesOfSmoothisometries}]
  \step{1 -- Construction of $P_1,\ldots,P_k$}\smallskip

  Since $V$ is open and $f(x')\colonequals (-N_2(x'),N_1(x'))$ is bounded and locally Lipschitz continuous in $V$, for each $x'\in V$, we can find $T_{x'}>0$ and a ruling curve $\Gamma_{x'}$ such that $\Gamma_{x'}:[0,T_{x'}]\to V$ and $\Gamma(T_{x'}/2)=x'$;
  indeed, we only need to solve the differential equation $\Gamma'_{x'}(t)=f(\Gamma_{x'}(t))$ with initial value given by ${x'}$.
  We denote by $P_{x'}\subseteq V$ the associated patch.

  Since $V$ is open and bounded, we can choose $\delta>0$ such that $|V\setminus V_\delta|<\varepsilon$, where $V_\delta\colonequals \{{x'}\in V\,:\,B({x'};\delta)\subseteq V\}$.
  Then $\overline {V_\delta}\subseteq\bigcup_{{x'}\in V_\delta}P_{x'}$.
  Since $\overline {V_\delta}$ is compact, we can find $x_1',\ldots,x_k'\in V_\delta$ such that with $P_i\colonequals P_{x_i'}$ we have $V_\delta\subseteq P_1\cup\ldots\cup P_k$ and thus $|V\setminus (P_1\cup\ldots\cup P_k)|\leq |V\setminus V_\delta| <\varepsilon$.
  So far, the patches are not disjoint.
  However, it is easy to check that intersections $P_i\cap P_j$ or differences $P_i\setminus \overline{P_j}$ of two patches are either empty or patches as well.
  Therefore, w.l.o.g., we may assume that the patches $P_1,\ldots, P_k$ are disjoint.
  This proves the first part of the lemma.

  \step{2 -- The case $y\in\mathcal A_{BC}$}\smallskip

  Fix $i,j$.
  To shorten the notation, we simply write $\mathcal L$, $P$, $\Gamma$ and $T$ instead of $\mathcal L_j$, $P_i$, $\Gamma_i$ and $T_i$.
  For the following argument it is useful to note that
  \begin{equation}\label{eq:patchbdry}
    \partial P\cap\partial S=\bigcup_{t\in[0,T]}\overline{[\Gamma(t);N(\Gamma(t))]}\cap\partial S.
  \end{equation}
  Indeed, this follow from \eqref{eq:patchdom} and convexity of $S$.

  We show that $\overline{\mathcal L}\cap \partial P=\emptyset$ by contraposition and thus assume that  $\overline{\mathcal L}\cap \partial P\neq\emptyset$.
  Then, by \eqref{eq:patchbdry}, there is $x_0'\in V\subseteq \{\II_y\neq 0\}$ such that $\overline{[x_0';N(x_0')]}\cap \overline{\mathcal L}\neq \emptyset$.
  But $S$ is convex, $N$ is continuous, and $V$ and $\mathcal L$ are (relatively) open.
  So we even find $x_1'\in V\subseteq \{\II_y\neq 0\}$ such that $\overline{[x_1';N(x_1')]}\cap \mathcal L\neq \emptyset$.
  For the same reasons, there is $\varepsilon>0$ such that 
  \begin{equation*}
  \forall x_2'\in B(x_1';\varepsilon):\quad \overline{[x_2';N(x_2')]}\cap{\mathcal L}\neq \emptyset.
  \end{equation*}
  By $y\in\mathcal A_{BC}\cap C^\infty(\overline{S};\R^3)$ and \eqref{EqFourteenA}, we conclude that $\nabla'y$ is constant on $B(x_1';\varepsilon)$ -- a contradiction to $x_1'\in\{\II_y\neq 0\}$.
\end{proof}

\begin{proof}[Proof of Lemma~\ref{L:monge1}]$ $

  \step{1 -- Preliminaries}\smallskip

For convenience set $N(t)\colonequals N(\Gamma(t))$.
  Let $\kappa,\kappa_{\mathrm{n}},S_\Gamma$ be defined as in Remark~\ref{R:linecurvparam} and note that by convexity of $S$ (c.f.~\cite{Pakzad,HornungApprox})
  \begin{equation*} 
    \partial P=\{\Gamma(t)+sN(t)\,:\,(t,s)\in\partial S_\Gamma\},\qquad \partial P\cap S=[\Gamma(0);N(0)]\cup[\Gamma(T);N(T)].
  \end{equation*}
By Remark \ref{R:linecurvparam}, $\Gamma$, $\kappa$ and $\kappa_\mathrm{n}$ are smooth on $[0,T]$.
Furthermore, since $|\II_y|$ is bounded on $\overline S$ and $\min_{[0,T]}|\kappa_{\mathrm{n}}|>0$, it follows by \eqref{eq:st:1121}, that
  $ 
    \inf_{(t,s)\in S_\Gamma}(1-s\kappa(t))>0
  $.

  \step{2 -- Construction of a solution to \eqref{eq:monge}}\smallskip

  Following \cite[Lemma~3.3]{SCHMIDT07} we consider the orthonormal basis of $\R^{2\times 2}_{\sym}$ given by
  \begin{equation*}
    F_1(t)\colonequals N(t)\otimes N(t),\qquad F_2(t)\colonequals \sqrt{2}\sym(N(t)\otimes\Gamma'(t)),\qquad F_3(t)\colonequals \Gamma'(t)\otimes\Gamma'(t),
  \end{equation*}
  and represent $M$ with help of the coefficient functions $\bar m_1,\bar m_2,\bar m_3:S_\Gamma\to\R$ defined by the identity $M(\Phi_\Gamma(t,s))=\sum_{i=1}^3\bar m_i(t,s)F_i(t)$,  where $\Phi_\Gamma(t,s):=\Gamma(t)+sN(t)$ denotes the diffeomorphism of Remark~\ref{R:linecurvparam}.
  Since $M$ is smooth and compactly supported in $P$, we conclude via Step 1 that $\bar m_i\in C^\infty_c(S_\Gamma)$.
  Next, we consider the system
  \begin{equation}\label{eq:st:12312123:sys}
    \begin{aligned}
    \partial_s\bar g_2=\,&\bar m_1,\\
    -\partial_s\bar g_1-\frac{\kappa_{\mathrm{n}}(t)}{1-s\kappa(t)}\bar g_1-\frac{1}{1-s\kappa(t)}\partial_t\bar g_2=\,&\sqrt{2}\bar m_2.
  \end{aligned}
  \end{equation}
  We shall see below, that it admits a solution $\bar g=(\bar g_1,\bar g_2)$ in $X$, where
  \begin{equation*}
    X\colonequals \{\bar u\in C^\infty(\overline{ S_{\Gamma}})\,:\,\exists \e>0\text{ such that }\bar u(t,\cdot)=0\text{ for all $t\in[0,\e)\cup(T-\e,T]$}\,\}.
  \end{equation*}
  Furthermore, for $x'\colonequals \Gamma(t)+sN(t)\in \overline P$, we set $ g(x')\colonequals \bar g_1(t,s)\Gamma'(t)+\bar g_2(t,s)N(t)$ and $\bar G(t,s)\colonequals \nabla'_{\sym}g(x)$.
  We let $\bar \alpha$ be defined by 
  \begin{equation*}
    \bar G\cdot F_3(t)+\bar\alpha\frac{\kappa_{\mathrm{n}}(t)}{1-s\kappa(t)}=\,\bar m_3.
  \end{equation*}
  Since $\frac{\kappa_{\mathrm{n}}(t)}{1-s\kappa(t)}$ is smooth and non-negative, $\bar\alpha\in X$ is uniquely defined.
  A direct calculation (which is detailed in the proof of \cite[Lemma~3.3]{SCHMIDT07}) then shows that the functions $g$ and $\alpha\colonequals \bar\alpha\circ \Phi_\Gamma^{-1}$ solve \eqref{eq:monge}.
  Furthermore, from $\bar g_1,\bar g_2,\bar\alpha\in X$ we infer that $\alpha$ and $g$ are smooth on $\overline P$ and vanish in a uniform neighborhood of $\partial P\cap S$.
  
  To solve the system \eqref{eq:st:12312123:sys}, we note that by Step~1 the coefficients and the right-hand side of the system are smooth and bounded. Since $\bar m_1\in X$ we can easily find $\bar g_2\in X$ solving the first equation. Likewise, since $\sqrt{2}\bar m_2+\frac{1}{1-s\kappa(t)}\partial_t\bar g_2\in X$, integration of the second equation yields a solution $\bar g_1\in X$.

\end{proof}

\subsubsection{Proof of Lemma \ref{L:recov_general} -- flat case}\label{Sec:Proof:L:recov_general_flat}
A key ingredient in the construction is the following approximation result for the Monge-Amp\'ere equation, which we recall from \cite{LP17}:
\begin{lemma}[\mbox{\cite[Proposition~3.2]{LP17}}]\label{L:convex}
  Let $V\subseteq\R^2$ be open and bounded, $A\in C^\infty(\overline V;\R^{2\times 2}_{\sym})$, and $\delta>0$.
  Suppose that $A(x')-\delta I_{2\times 2}$ is positive definite for all $x'\in \overline V$. Then for all $n\in\N$ there exists $g_n\in C^\infty(\overline V;\R^2)$ and $\alpha_n\in C^\infty(\overline V;\R)$ such that
  \begin{equation}\label{L:convex:eq1}
    \lim\limits_{n\to\infty}\Big(\|A-\big(\sym\nabla'g_n+\frac12\nabla'\alpha_n\otimes\nabla'\alpha_n\big)\|_{L^\infty(V)}+\|\alpha_n\|_{L^\infty(V)}+\|g_n\|_{L^\infty(V)}\Big)=0,
  \end{equation}
  and
  \begin{equation}\label{L:convex:eq2}
    \|\nabla'g_n\|_{L^\infty(V)}+\|\nabla'\alpha_n\|_{L^\infty(V)}\leq C\|A\|_\infty,
  \end{equation}
  for some constant $C>0$ independent of $n$ and $A$.
\end{lemma}

\begin{proof}[Proof of Lemma \ref{L:recov_general} -- flat case]
  For convenience, we set $\Omega_V\colonequals V\times(-\frac12,\frac12)$.
  W.l.o.g. we may assume that $d\in C^\infty_c(\Omega_V;\R^3)$, $M\in C^\infty(\overline V;\R^{2\times 2}_{\sym})$, and $M+\delta I_{2\times 2}$ is uniformly positive definite on $\overline V$.
  The general case then follows by an approximation argument and the extraction of a diagonal sequence.

  \step{1 -- Modification of Lemma~\ref{L:convex}}\smallskip

  We claim that there exist sequences $(g_n)\subseteq C^\infty_c(V;\R^3)$, $(\alpha_n)\subseteq C^\infty_c(V)$ satisfying
  \begin{equation*}
    \|\nabla'g_n\|_{L^\infty(V)}+\|\nabla'\alpha_n\|_{L^\infty(V)}\leq C(\|M\|_\infty+|\delta|)
  \end{equation*}
  (with $C>0$ independent of $M$ and $n$) and
  \begin{equation}\label{eq:st:1112323aa}
    \lim_{n\to\infty}\Big(\|(M+\delta I_{2\times 2})-\big(\sym\nabla'g_n+\frac12\nabla'\alpha_n\otimes\nabla'\alpha_n\big)\|_{L^2(V)}+\|\alpha_n\|_{L^\infty(V)}+\|g_n\|_{L^\infty(V)}\Big)=0.
  \end{equation}
  This can be seen as follows: Let $(\tilde g_n,\tilde \alpha_n)$ denote the sequence of Lemma~\ref{L:convex} applied with $A\colonequals M+\delta I_{2\times 2}$ (since $M+\delta I_{2\times 2}$ is uniformly positive definite, we can find $\delta'>0$ such that $A-\delta' I_{2\times 2}$ is positive definite on $\overline V$).
  Since $V$ is open and bounded, there exist cut-off functions $(\eta_k)\subseteq C^\infty_c(V;\R)$ satisfying $0\leq\eta_k\leq 1$ and $\lim_{k\to\infty}\|1-\eta_k\|_{L^4(V)}\to 0$.
  Set $g_{k,n}\colonequals \eta_k\tilde g_n$ and $\alpha_{k,n}\colonequals \eta_k\tilde \alpha_n$, and note that
  \begin{equation*}
    |\nabla'g_{k,n}|\leq |\nabla'\eta_k||\tilde g_n|+|\nabla'\tilde g_n|,\qquad     |\nabla'\alpha_{k,n}|\leq |\nabla'\eta_k||\tilde \alpha_n|+|\nabla'\tilde \alpha_n|.
  \end{equation*}
  Then,
  \begin{align*}
    C_{k,n}\,\colonequals \,&\|\tilde g_n-g_{k,n}\|_{L^\infty(V)}+    \|\tilde \alpha_n-\alpha_{k,n}\|_{L^\infty(V)} + \|\nabla'(\tilde g_n-g_{k,n})\|_{L^4(V)}+\|\nabla'(\tilde \alpha_n-\alpha_{k,n})\|_{L^4(V)}\\
    \leq\,&(1 + \|\nabla'\eta_k\|_{L^4(V)})\Big(\|\tilde g_n\|_{L^\infty(V)}+\|\tilde \alpha_n\|_{L^\infty(V)}\Big)\\
    &
    +\|1-\eta_k\|_{L^4(V)}\Big(\|\nabla'\tilde g_n\|_{L^\infty(V)}+\|\nabla'\tilde \alpha_n\|_{L^\infty(V)}\Big).
  \end{align*}
  Consider
  \begin{equation*}
    C'_{k,n}\colonequals C_{k,n}+
    \begin{cases}
      1&\text{if }\max\{\|\nabla'g_{k,n}\|_{L^\infty(V)},\|\nabla'\alpha_{k,n}\|_{L^\infty(V)}\}>2C\|M+\delta I_{2\times 2}\|_{L^\infty(V)},\\
      0&\text{else,}
    \end{cases}.
  \end{equation*}
  where $C$ is from \eqref{L:convex:eq2} (applied for $A:=M+\delta I$).
  By \eqref{L:convex:eq1} and \eqref{L:convex:eq2}, we have $\limsup_{k\to\infty}\limsup_{n\to\infty}C'_{k,n}=0$, and thus $C'_{k_n,n}\to 0$ for a diagonal sequence $(k_n)$. We conclude that $(g_{k_n,n},\alpha_{k_n,n})_n$ satisfies the claimed properties.

  \step{2 -- Construction of a sequence of von-K\'arm\'an-displacements}\smallskip

  Let $(g_n,\alpha_n)$ be as in Step~1 and consider the (von-K\'arm\'an-type) displacement
  \begin{align*}
    w_{n,h}\colonequals 
    \begin{pmatrix}
      g_n-h^{\frac12} x_3(\nabla'\alpha_n)^\top\\
      h^{-\frac12}\alpha_n-h\frac{x_3}{2}|\nabla'\alpha_n|^2
    \end{pmatrix}.
  \end{align*}
  By construction we have $w_{n,h}\in C^\infty(\overline{\Omega_V};\R^3)$ with
  \begin{equation}\label{eq:Ukraine}
  w_{n,h}=0 \text{ in a uniform neighbourhood of }\partial V\times(-\tfrac12,\tfrac12).
  \end{equation}
  Furthermore, a direct computation and the estimates of Step~1 show that there exists a constant $C>0$ such that for all $n\in\N$ we have
  \begin{align*}
    &\limsup_{h\to0}h^\frac12\Big(\|w_{n,h}\|_{L^\infty(\Omega_V)}+\|\nabla\!_hw_{n,h}\|_{L^\infty(\Omega_V)}\Big)\leq C,\\
    &\limsup_{h\to0}h^{-1}\Big\|(I_{3\times 3}+h\nabla\!_hw_{n,h})^\top(I_{3\times 3}+h\nabla\!_hw_{n,h})-I_{3\times 3}\Big\|_{L^\infty(\Omega_V)}<\infty,\\
    &\limsup_{h\to 0}h^{-\frac32}\Big\|(I_{3\times 3}+h\nabla\!_hw_{n,h})^\top(I_{3\times 3}+h\nabla\!_hw_{n,h})\\
    	&\qquad\qquad-\Big(I_{3\times 3}+2h\iota\Big(\sym\nabla'g_n+\frac12\nabla'\alpha_n\otimes\nabla'\alpha_n\Big)\Big)\Big\|_{L^\infty(\Omega_V)}\leq C.
  \end{align*}
  Hence, we may pass to a diagonal sequence $w_h\colonequals w_{n_h,h}$ with $\lim_{h\to 0}n_h=\infty$ such that
  \begin{align}\label{eq:prop:wh1}
    &\limsup_{h\to0}h^\frac12\Big(\|w_{h}\|_{L^\infty(\Omega_V)}+\|\nabla\!_hw_{h}\|_{L^\infty(\Omega_V)}\Big)<\infty,\\\label{eq:prop:wh2}
    &\limsup_{h\to0}h^{-1}\Big\|(I_{3\times 3}+h\nabla\!_hw_{h})^\top(I_{3\times 3}+h\nabla\!_hw_{h})-I_{3\times 3}\Big\|_{L^\infty(\Omega_V)}<\infty,\\\label{eq:prop:wh3}
    &\limsup_{h\to 0}h^{-\frac32}\Big\|(I_{3\times 3}+h\nabla\!_hw_{h})^\top(I_{3\times 3}+h\nabla\!_hw_{h})-\Big(I_{3\times 3}+2h\iota(M+\delta I_{2\times 2})\Big)\Big\|_{L^2(\Omega_V)}<\infty,
  \end{align}
  where the last estimate is a consequence of \eqref{eq:st:1112323aa}.

  \step{3 -- Construction of the recovery sequence}\smallskip

  Let $(w_h)\subseteq C^\infty(\overline{\Omega_V})$ be the sequence constructed in Step~2. We consider the sequence
  \begin{align*}
    y_{h}\colonequals \,&(1-h\delta)y+hx_3 b_y+h^2R_y\int_0^{x_3}d(\cdot,t)\,\mathrm dt+hR_yw_h.
  \end{align*}
  By $d\in C_c^\infty(\Omega_V)$ and \eqref{eq:Ukraine}, we have $y_h\in C^\infty(\overline{\Omega_V};\R^3)$ and
  \begin{equation*}
    y_h=(1-h\delta)y+hx_3 b_y\qquad\text{in a uniform neighborhood of }\partial V\times(-\frac12,\frac12).
  \end{equation*}
  Thanks to  \eqref{eq:prop:wh1} we have
  \begin{equation*}
    \|y_h-y\|_{L^\infty(\Omega_V)}\to 0.
  \end{equation*}
  Since $y$ is affine in $V$ by assumption, $R_y$ and $b_y$ are constant, and thus,
  \begin{equation*}
    \nabla\!_h y_h=R_y\Big(I_{3\times 3}-h\delta\iota(I_{2\times 2})+h\nabla\!_hw_h+h \,d\otimes e_3+h^2\int_0^{x_3}\nabla'd(\cdot,t)\,\mathrm dt\Big).
  \end{equation*}
  With help of \eqref{eq:prop:wh1}, we obtain the claimed convergence of $(\nabla\!_hy_h)$ in \eqref{eq:recov:onescalea}.
  To prove convergence of $E_h(y_h)$, we consider
  \begin{equation*}
    G_h\colonequals \frac{1}{h}((\nabla\!_h y_h)^\top\nabla\!_hy_h-I_{3\times 3}).
  \end{equation*}
  A direct calculation that uses \eqref{eq:prop:wh2}, \eqref{eq:prop:wh3} yields
  \begin{equation*}
    \limsup\limits_{h\to 0}\|G_h\|_{L^\infty(\Omega_V)}<\infty,\qquad \lim\limits_{h\to 0}\|G_h-2\big(\iota(M)+\sym(d\otimes e_3)\big)\|_{L^2(\Omega_V)}=0.
  \end{equation*}
  Hence, by a Taylor expansion of the matrix square root around $I_{3\times 3}$ in $\sqrt{(\nabla\!_hy_h)^\top\nabla\!_hy_h}=\sqrt{I_{3\times 3}+hG_h}$, we conclude
  \begin{equation*}
    \limsup\limits_{h\to 0}\|E_h(y_h)-\big(\iota(M)+\sym(d\otimes e_3)\big)\|_{L^2(\Omega_V)}=0.
  \end{equation*}
\end{proof}

\subsection{Density of smooth isometries: proof of Proposition~\ref{P:approx}} \label{Sec:Proof:P:approx}
\newcommand{\yext}{{\hat y}}
\newcommand{\Sext}{{\widehat S}}

In this section, we prove Proposition~\ref{P:approx}.
In fact, we prove a slightly stronger statement, where instead of \eqref{ass:domain}, we only assume the following condition introduced in \cite{HornungApprox}:

\begin{equation}
  \label{ConditionStar}
  \begin{aligned}
    &\text{$S\subseteq\R^2$ is a bounded Lipschitz domain,}\\
      &\text{and there is a closed subset }\Sigma\subseteq \partial S\text{ with }\mathcal{H}^1(\Sigma)=0, \\
      &\text{such that the unit outer normal exists on }\partial S\setminus \Sigma\text{ and is continuous on that set.}
\end{aligned}
\end{equation}
Our proof is a consequence of the theory developed in \cite{Kirchheim, Pakzad, HornungApprox, HornungFine} and builds on ideas of \cite{Olb}.
In particular, we shall appeal to the following result, which follows from the proof of \cite[Theorem~1]{HornungApprox}:

\begin{proposition}[cf.~\mbox{\cite[Theorem~1]{HornungApprox}}]\label{P:aux}
  Let $S\subseteq\R^2$ satisfy \eqref{ConditionStar} and let $y\in H^2_{\iso}(S;\R^3)$ be finitely developable in the sense that $\widehat C_{\nabla'y}$ consists of finitely many connected components. Then for all $\delta>0$ there exists $y_\delta\in C^\infty(\overline S;\R^3)\cap H^2_{\iso}(S;\R^3)$ such that
  \begin{align*}
    \|y-y_\delta\|_{H^2(S;\R^3)}<\delta\qquad\text{and}\qquad y=y_\delta\text{ on each connected component }U\subseteq\widehat C_{\nabla'y}.
  \end{align*}
\end{proposition}
We shall deduce Proposition~\ref{P:approx} from Proposition~\ref{P:aux} by first extending $y\in\mathcal A_{BC}$ to an isometry $\yext$ defined on a larger domain $\Sext$ such that each line segment $\mathcal L_i$ is contained in the closure of a connected component  $U_i\subseteq\widehat C_{\nabla'\yext}$, see Lemma~\ref{LemmaExtendTheIsometricImmersion} below.
Then, in a second step, we check that $\yext$ can be approximated by an isometry that is finitely developable and that is equal to $\yext$ on the components $U_i$.
An application of Proposition~\ref{P:aux} then yields a smooth isometry $\yext^\delta$ that approximates $\yext$ and satisfies $\yext^\delta=\hat y$ on the components $U_i$.
We conclude that $\yext^\delta\big\vert_S\in\mathcal A_{BC}$.

\begin{lemma} \label{LemmaExtendTheIsometricImmersion}
  Let $S\subseteq \R^2$ satisfy \eqref{ConditionStar} and suppose Assumption~\ref{ass:BC}. Then there exists $\Sext\subseteq \mathbb{R}^2$ such that:
  \begin{enumerate}[(a)]
  \item $\Sext$ satisfies \eqref{ConditionStar}, $S\subseteq\Sext$, and $\Sext\setminus\overline S$ is the union of mutually disjoint, open, connected sets $D_1,\ldots,D_{k_{BC}}$ such that $\overline{\mathcal L_i}=\overline D_i\cap \overline S$ for $i=1,\ldots,k_{BC}$.
  \item  For all non-affine $y\in\mathcal A_{BC}$ there exists an extension $\yext\in H^2_\mathrm{iso}(\Sext,\mathbb{R}^3)$ such that $\yext=y$ on $S$, and for all $i=1,\ldots,k_{BC}$ there exists exactly one connected component $U_i\subseteq\widehat C_{\nabla'\yext}$ with $D_i\cap U_i\neq\emptyset$.
  Furthermore, $\mathcal L_i\subseteq\overline{U_i}$ and $\partial U_i\cap D_i$ consists of exactly two connected components.

  \end{enumerate}
\end{lemma}
\begin{figure}
\center
\definecolor{ttqqqq}{rgb}{0,0,0}
\definecolor{ududff}{rgb}{0.0,0.0,0}
\begin{tikzpicture}[line cap=round,line join=round,>=triangle 45,x=1cm,y=1cm]
\clip(-5.4,-2.3) rectangle (5.0,2);
\fill[line width=.2pt,color=ttqqqq,fill=ttqqqq,fill opacity=0.06] (4,0) -- (-4,0) -- (-3.429305223439793,1.9292572089754076) -- (-7.372384946814471,-3.146546589895281) -- (-0.8452414666395512,-4.424274896007801) -- (5.566271445426404,-4.207720977796457) -- cycle;
\fill[line width=.2pt,color=ttqqqq,fill=ttqqqq,fill opacity=0.2] (-4,0) -- (4,0) -- (2,1.5432311678355755) -- cycle;
\fill[line width=.2pt,color=ttqqqq,fill=ttqqqq,fill opacity=0.2] (4,0) -- (-2,1.5432311678355755) -- (-4,0) -- cycle;
\draw [line width=.2pt,color=ttqqqq] (4,0)-- (-4,0);
\draw [line width=.2pt,color=ttqqqq] (-4,0)-- (-3.429305223439793,1.9292572089754076);
\draw [line width=.2pt,color=ttqqqq] (-3.429305223439793,1.9292572089754076)-- (-7.372384946814471,-3.146546589895281);
\draw [line width=.2pt,color=ttqqqq] (-7.372384946814471,-3.146546589895281)-- (-0.8452414666395512,-4.424274896007801);
\draw [line width=.2pt,color=ttqqqq] (-0.8452414666395512,-4.424274896007801)-- (5.566271445426404,-4.207720977796457);
\draw [line width=.2pt,color=ttqqqq] (5.566271445426404,-4.207720977796457)-- (4,0);
\draw (-0.32551533363693336,0.05) node[anchor=north west] {$\mathcal L_1$};
\draw (-1.0708348814852182,-1.2) node[anchor=north west] {\huge $S$};
\draw (-4.2,-0.05) node[anchor=north west] {$x_1'$};
\draw (3.4,-0.05) node[anchor=north west] {$x_2'$};
\draw [line width=.2pt] (-4,0)-- (-2,1.5432311678355755);
\draw [line width=.2pt] (2,1.5432311678355755)-- (4,0);
\draw [line width=.2pt] (2,1.5432311678355755)-- (0,1.0288207785570502);
\draw [line width=.2pt] (-2,1.5432311678355755)-- (0,1.0288207785570502);
\draw [line width=.2pt, dashed] (0,1.0288207785570502)-- (4,0);
\draw [line width=.2pt, dashed] (0,1.0288207785570502)-- (-4,0);
\draw [line width=.2pt,color=ttqqqq] (-4,0)-- (4,0);
\draw [line width=.2pt,color=ttqqqq] (4,0)-- (2,1.5432311678355755);
\draw [line width=0.2pt,color=ttqqqq] (2,1.5432311678355755)-- (-4,0);
\draw [line width=0.2pt,color=ttqqqq] (4,0)-- (-2,1.5432311678355755);
\draw [line width=.2pt,color=ttqqqq] (-2,1.5432311678355755)-- (-4,0);
\draw [line width=.2pt,color=ttqqqq] (-4,0)-- (4,0);
\draw (1.1,1.3) node[anchor=north west] {$\Delta_1\setminus\Delta_2$};
\draw (-2.7,1.3) node[anchor=north west] {$\Delta_2\setminus\Delta_1$};
\draw (-1,0.8) node[anchor=north west] {$\Delta_1\cap\Delta_2$};

\draw [line width=.2pt, loosely dashed] (-2.3333333333333335,1.2860259731963126)-- (-1.666666666666666,1.457496102955821);
\draw [line width=.2pt, loosely dashed] (-1.3333333333333335,1.371761038076067)-- (-2.666666666666666,1.0288207785570507);
\draw [line width=.2pt, loosely dashed] (-3,0.7716155839177877)-- (-1,1.2860259731963126);
\draw [line width=.2pt, loosely dashed] (-0.6666666666666665,1.2002909083165585)-- (-3.333333333333333,0.5144103892785253);
\draw [line width=.2pt, loosely dashed] (-3.666666666666666,0.2572051946392629)-- (-0.3333333333333335,1.1145558434368044);

\draw [line width=.2pt, loosely dashed] (0.3333333333333335,1.1145558434368044)-- (3.666666666666666,0.2572051946392629);
\draw [line width=.2pt, loosely dashed] (3.333333333333333,0.5144103892785253)-- (0.6666666666666665,1.2002909083165585);
\draw [line width=.2pt, loosely dashed] (1,1.2860259731963126)-- (3,0.7716155839177877);
\draw [line width=.2pt, loosely dashed] (2.666666666666666,1.0288207785570507)-- (1.3333333333333335,1.371761038076067);
\draw [line width=.2pt, loosely dashed] (1.666666666666666,1.457496102955821)-- (2.3333333333333335,1.2860259731963126);
\draw (-2.0,2.1) node[anchor=north west] {\scalebox{0.7}{$\frac{3x_1'+x_2'}{4}+\rho\nu$}};
\draw (2,2.1) node[anchor=north west] {\scalebox{0.7}{$\frac{x_1'+3x_2'}{4}+\rho\nu$}};
\begin{scriptsize}
\draw [fill=ududff] (-7.372384946814471,-3.146546589895281) circle (2.5pt);
\draw [fill=ududff] (-0.8452414666395512,-4.424274896007801) circle (2.5pt);
\draw [fill=ududff] (5.566271445426404,-4.207720977796457) circle (2.5pt);
\end{scriptsize}
\end{tikzpicture}
\caption{In the constructions of Lemma \ref{LemmaExtendTheIsometricImmersion}, the domain $S$ is extended by the triangles $\Delta_1$ and $\Delta_2$.
On the dashed regions, the corresponding extended deformation is curved (with the dashed lines indicating the field of asymptotic directions).
On $\Delta_1\cap\Delta_2$, the extended isometric deformation is flat.}
\label{Fig:Extend:Triangles}
\end{figure}
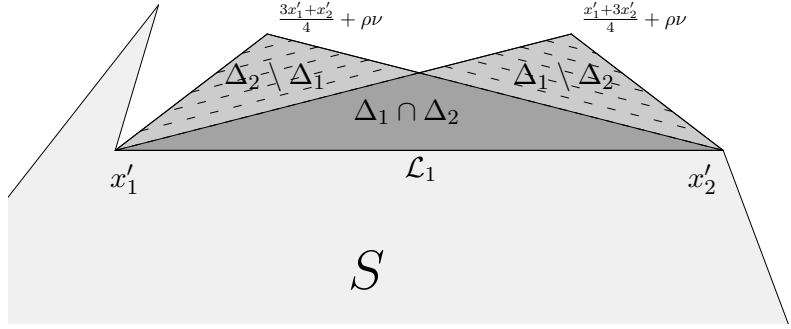
\begin{proof}
  It suffices to consider the case $ k_{BC} = 1$ since the argument below can be iterated finitely many times.
  The idea of the proof is sketched in Figure \ref{Fig:Extend:Triangles}.

  \step{1 - Proof of (a)}\smallskip

  Since $\mathcal L_1\subseteq\partial S$ is a relatively open line segment, there exists $x_1',x_2'\in\partial S$ such that $\mathcal L_1=\mbox{conv}\{x_1',x_2'\}$.
  Let $\nu\in\mathbb S^1$ denote the outer unit normal on $\mathcal L_1$.
  Since $\partial S$ is Lipschitz, there exists $\rho>0$ such that the triangles $\Delta_1\colonequals \mbox{conv}\{x_1',x_2',\frac{x_1'+3x_2'}{4}+\rho\nu\}$ and $\Delta_2\colonequals \mbox{conv}\{x_1',x_2',\frac{3x_1'+x_2'}{4}+\rho\nu\}$ have empty intersections with $S$.
  Set $D_1\colonequals \Delta_1\cup\Delta_2$.
  We conclude that $\Sext\colonequals \interieur(\overline{S\cup D_1})$ satisfies the claimed properties.

  \step{2 - Proof of (b)}\smallskip
  
  Let $\Sext$ as in (a) and let $y\in\mathcal A_{BC}$ be non-affine.
  Since $\nabla'y$ is constant on $\mathcal L_1$, there exists $\yext\in H^{2}_{\iso}(\Sext;\R^3)$ such that $\yext=y$ in $S$, and $\II_{\yext}=0$ in $\Delta_1\cap\Delta_2$, and $\II_{\yext}$ is constant and non-zero on each of the two sets $\Delta_1\setminus\overline{\Delta_2}$ and $\Delta_2\setminus\overline{\Delta_1}$.
  Let $U_1$ denote the connected component of $C_{\nabla'\yext}$ that contains $\Delta_1\cap\Delta_2$. Since $y$ is assumed to be non-affine, $\partial U_1\cap\Sext$ must consist of at least three connected components.
  We conclude that $\yext$ satisfies the sought for properties.
\end{proof}

We are now in the position to prove Proposition \ref{P:approx}.

\begin{proof}[Proof of Proposition \ref{P:approx}]
  Let $y\in\mathcal A_{BC}$ and assume w.l.o.g. that $y$ is non-affine.
  Let $\Sext$ and $\yext$ be as in Lemma~\ref{LemmaExtendTheIsometricImmersion}. For $i=1,\ldots,k_{BC}$ let $D_i$ and $U_i$  be as in Lemma~\ref{LemmaExtendTheIsometricImmersion}. We apply \cite[Proposition~5]{HornungApprox} and obtain a sequence $(u^\delta)$ of finitely developable isometries in $H^2_{\iso}(\Sext;\R^3)$ such that:
  \begin{itemize}
  \item   $\|u^\delta-\yext\|_{H^2(\Sext;\R^3)}\to 0$,
  \item for all $\delta>0$ (sufficiently small) we have $u^\delta=\yext$ on each connected component $U\subseteq \widehat C_{\nabla'\yext}$ satisfying
    \begin{equation}\label{eq:approx:crit}
      U\cap\{x'\in\Sext\,:\,\dist(x',\partial\Sext)>\delta\}\neq\emptyset.
    \end{equation}
  \end{itemize}
  Since the line segments $\mathcal L_i$ are contained in the open set $\Sext$, for all $\delta>0$ (sufficiently small), the components $U_i$, $i=1,\ldots,k_{BC}$, satisfy \eqref{eq:approx:crit}.
  Thus, we conclude that
  \begin{equation}\label{eq:st:112sss12ll}
    u^\delta=\yext\text{ in }U_i\text{ for all }i=1,\ldots,k_{BC}.
  \end{equation}
  Note that this implies $(u^\delta,\nabla'u^\delta)=(\yext,\nabla'\yext)$ on $\mathcal L_i$, since we have $\mathcal L_i\subseteq\overline{ U_i}$ by construction (and $H^2_{\iso}(S;\R^3)\subseteq C^1(S;\R^3)$).
  In particular, we have $u^\delta\vert_S\in\mathcal A_{BC}$.
  If $u^\delta\vert_S$ is affine, then it is smooth and thus the desired approximation.
  If $u^\delta\vert_S$ is not affine, we proceed as follows. 
  From \eqref{eq:st:112sss12ll} we obtain that $(u^\delta,\nabla'u^\delta)=(\yext,\nabla'\yext)$ on $\partial U_i\cap D_i$, and thus,
  \begin{equation*}
    \yext^\delta(x')\colonequals 
    \begin{cases}
      \yext(x')&x\in \bigcup_{i=1}^{k_{BC}}(D_i\setminus U_i),\\
      u^\delta(x')&\text{else.}
    \end{cases}
  \end{equation*}
  defines an isometry in $H^2_{\iso}(\Sext;\R^3)$.
  Since $\yext$ is curved in $D_i\setminus U_i$, we conclude that $\yext^\delta$ is finitely developable.
  Furthermore, since $\partial U_i\cap D_i$ consists of at least two connected components, and since $u^\delta\vert_S$ is non-affine, there exist connected components $U^\delta_i\subseteq\widehat C_{\nabla'\yext^\delta}$ with $U_i\subseteq U^\delta_i$.
  An application of Proposition~\ref{P:aux} to $\yext^\delta$ thus yields an isometry $y^\delta\in H^2_{\iso}(S;\R^3)\cap C^\infty(\overline\Sext;\R^3)$ with $\|y^\delta-\yext^\delta\|_{H^2(\Sext;\R^3)}<\delta$ such that $y^\delta=\yext^\delta=\yext$ on $U_i$ for $i=1,\ldots,k_{BC}$.
  Combined with the property $\mathcal L_i\subseteq\overline {U_i}$ we conclude that $y^\delta\vert_S\in\mathcal A_{BC}\cap C^\infty(\overline S;\R^3)$ and $y^\delta\to y$ in $H^2(S;\R^3)$.
\end{proof}

\subsection*{Acknowledgments}
The authors acknowledge support by the German Research Foundation (DFG) via the
re\-search unit FOR 3013 ``Vector- and tensor-valued surface PDEs'' (grant no.~BA2268/6--1 and \mbox{NE2138/4-1}).

\bibliographystyle{plain}
\bibliography{bibliography}

\end{document}